\documentclass[times]{nmeauth}
\pdfoutput=1
\newcommand\BibTeX{{\rmfamily B\kern-.05em \textsc{i\kern-.025em b}\kern-.08em
T\kern-.1667em\lower.7ex\hbox{E}\kern-.125emX}}

\def\volumeyear{2014}

\usepackage {amsmath}
\usepackage {amssymb}
\usepackage{graphicx}
\usepackage{empheq}
\usepackage{dsfont}
\usepackage{subfig}
\usepackage{booktabs}
\usepackage{algorithm}
\usepackage{algpseudocode}
\floatname{algorithm}{Procedure}

\usepackage{cite}
\usepackage[noprefix]{nomencl}
\usepackage{multirow}

\usepackage[noprefix]{nomencl}
\makenomenclature

\newcommand{\subeq}[2]{\eqref{#1}$_{#2}$}
\newcommand{\rhof}{\rho_\text{f} }
\newcommand{\rhos}{\rho_\text{s} }
\newcommand{\kh}{K_\text{h} }
\newcommand{\Nn}{N_\text{n}}
\newcommand{\Nfc}{N_{\text{eg}}}
\newcommand{\Nel}{N_{\text{el}}}
\newcommand{\prsc}{\text{prsc}}

\newcommand{\bld}{\boldsymbol}
\newcommand{\Mf}{M_\text{f}}

\newcommand{\Nfn}{N_\text{fn}}
\newcommand{\Nffc}{N_{\text{feg}}}
\newcommand{\divg}[1]{{\nabla} \cdot #1} 
\newcommand{\grad}[1]{{\nabla} #1} 
\newcommand{\ms}{\frac{\text{m}}{\text{s}}}
\newcommand{\linprt}{\text{P}_1-\text{RT}_0}
\newcommand{\quadprt}{\text{P}_2-\text{RT}_0}
\newcommand{\quadElem}{\text{P}_2-\text{P}_0}
\newcommand{\linElem}{\text{P}_1-\text{P}_0}

\newcommand{\nf}{n_\text{f}}

\newcommand{\rt}{\text{RT}_0}

\newcommand{\sigmat}{\sigma^{\text{t}} } 
\newcommand{\tprsc}{T_{\text{prsc}}}

\begin{document}
\runningheads{Lotfian and Sivaselvan}{Dynamic Poroelasticity}

\title{A topology-motivated mixed finite element method for dynamic response of porous media}

\author{Z.~Lotfian  and M.~V.~Sivaselvan} 

\address{Department of Civil, Structural and Environmental Engineering, University at Buffalo, Buffalo, NY 14221, U.S.A.}

\corraddr{zahrasad@buffalo.edu}

\begin{abstract}
In this paper, we propose a numerical method for computing solutions to Biot's fully dynamic model of incompressible saturated porous media \cite{Biot_lowfrequency_1956}. Our spatial discretization scheme is based on the three-field formulation ($u$-$w$-$p$) and the coupling of a lowest order Raviart-Thomas mixed element \cite{Raviart_Thomas_1977} for fluid variable fields ($w$, $p$ ) and a nodal Galerkin finite element for skeleton variable field ($u$). These mixed spaces are constructed based on the natural topology of the variables; hence, are physically compatible and able to exactly model the kind of continuity which is expected. The method automatically satisfies the well known LBB (inf-sup) stability condition and avoids locking that usually occurs in the numerical computations in the incompressible limit and very low hydraulic conductivity. In contrast to the majority of approaches, our three-field formulation can fully capture dynamic behavior of porous media even in high frequency loading phenomena with considerable fluid acceleration such as liquefaction and biomechanics of porous tissues under rapid external loading. Moreover, we address the importance of consistent initial conditions for poroelasticity equations with the incompressibility constraint, which represent a system of differential algebraic equations. The energy balance equation is derived for the full porous medium and used to assess the stability and accuracy of our time integration. To highlight the capabilities of our method, a variety of numerical studies are provided including verification with analytical and boundary element solutions, wave propagation analyses, hydraulic conductivity effects on damping and frequency content, energy balance analyses, mass lumping considerations, effects of mesh pattern and size, and stability analyses. We also explain some discrepancies commonly found in dynamic poroelasticity results in the literature.
\end{abstract}


\maketitle
\section{Introduction}
Dynamics of porous media are of interest in a variety of fields including mechanics of saturated and partially saturated soil, porous biological materials and petroleum reservoirs. The dynamic response of saturated soil under earthquake and the subsequent soil liquefaction process is subject of many studies \cite{ Zienkiewicz_Shiomi_1984,Jeremić_etal_2008, Li_etal_2004}. Additionally, the ability to accurately model the behavior of rock/soil and fluids in petroleum reservoirs is essential for maximizing oil recovery and assessing the performance of processes such as wellbore stability and subsidence. For this also, coupled modeling of geomechanics and fluid flow is required \cite{Wan_2002, Jha_Juanes_2007}. Similarly in biomechanics, mechanical response of hydrated soft tissues such as brain, skin, and articular cartilage \cite{ Simon_etal_1985,Suh_etal_1991, Levenston_etal_1998, Zhu_ Suh_2001, Yang_2006}, bones \cite{Cowin_1999}, and even individual cells \cite{Moeendarbary2013,SachsSivaselvan2015} involves coupled response of a skeleton and pore fluid. Furthermore, all the results of dynamic poromechanics can be applied to dynamic thermomechanics based on the mathematical analogy between the formulations \cite{Biot_thermo_1956,Biot_1962}.

In this paper, we propose a numerical method for dynamic poroelasticity, starting from Biot's formulation \cite{Biot_lowfrequency_1956}. The porous medium is considered as a solid skeleton saturated with fluid. While the skeleton itself is deformable and elastic, the material that makes up the skeleton and the fluid are modeled as incompressible. This is reasonable for many of the physical problems described above. We also restrict the formulation to linearized kinematics. Alternate formulations based on the Theory of Porous Media \cite{Bowen_1980} have been shown to be equivalent in this realm of incompressible solid and fluid materials and linear kinematics \cite{Ehlers_Kubik_1994}. Our focus is on fully dynamic problems, i.e., inertia effects in both the fluid and the skeleton have been taken into account. We note that this is in contrast to most previous work, where fluid acceleration terms are fully or partially neglected \cite{Zienkiewicz_Shiomi_1984,Li_etal_2004,Zhu_Suh_2001,Oka_etal_1994}.

Some of the specific features of our approach are as follows: (a) In contrast to much of the literature, we use a full three-field formulation in which the skeleton displacement ($u$), pore pressure ($p$) and pore-fluid relative velocity ($w$) fields are all approximated independently and calculated directly. As a result (see section \ref{ThreeFieldFormulation}), the formulation is applicable to problems with considerable fluid acceleration, for instance in liquefaction analysis \cite{Jeremić_etal_2008} and biomechanics of porous tissues under rapid external loading \cite{Yang_2006,Levenston_et al_1998}. Moreover, greater accuracy is attained compared to formulations where some fields are computed by post processing techniques as reported in \cite{Simon_etal_higherOrder_1986,Gajo_etal_1994}. (b) We employ a topology-motivated spatial discretization based on the lowest order Raviart-Thomas ($\rt$) mixed element \cite{Raviart_Thomas_1977} to approximate the flow variables and a linear/quadratic standard Galerkin finite element for the skeleton displacement. This mixed approach automatically enforces the precise kind of continuity (local mass conservation) that is required. Additionally, in the incompressible limit (when material is modeled as incompressible and under low hydraulic conductivity) and the rigid skeleton limit, the proposed method satisfies the well known LBB (inf-sup) condition \cite{bathe2006finite,boffi2013mixed}, avoiding associated locking effects and the need for any additional stabilization techniques. (c) Taking a differential algebraic equations (DAE) point of view, we demonstrate the importance of consistent initial conditions for poroelasticity equations with the incompressibility constraint. (d) The energy balance equation for the whole porous medium is also derived and used to confirm the stability and accuracy of the Newmark time integration method with particular emphasis on consistent initial condition.

In many problems, such as soil mechanics and biomechanics of soft tissues, the fluid and skeleton material could be reasonably modeled as being incompressible. However when developing a finite element formulation, the incompressibility condition imposes a certain constraint (see section \ref{LimitingCases}), which in turn gives rise to stability concerns. The discrete form of the  Lady\v{z}enskaja-Babu\v{s}ka-Brezzi (LBB) condition provides the necessary and sufficient criterion for a stable finite element formulation. To satisfy this condition, several approaches including stabilized methods \cite{Hughes_etal_1986, Hughes_etal_1989,Zienkiewicz_Wu_1991,Murad_etal_1996, Wan_2002,Hughes_etal_2006,White_Borja_2008}, discontinuous Galerkin \cite{Baumann_1997, Oden_etal_1998,Rivière_etal_2000, Bastian_2003, Chen_etal_2008, Scovazzi_etal_2013,Liu_2004}, special elements such as Taylor-Hood elements \cite{Taylor_Hood_1973} or MINI elements \cite{Arnold_etal_1984}, and mixed finite elements \cite{Raviart_Thomas_1977,Nédélec_1980,Chen_FiniteElement} have been proposed in the literature.
As we discuss in section \ref{LimitingCases}, in a finite element formulation for poroelasticity, the LBB condition must be satisfied in two limiting cases --- low hydraulic conductivity and rigid skeleton \cite{Wan_2002,White_2009}.
In this paper, we develop a finite element formulation using a two fold approach to satisfy the LBB condition --- (i) We discretize the pressure ($p$) and velocity ($w$) fields of the pore-fluid flow using Raviart-Thomas (RT) elements \cite{Raviart_Thomas_1977}. These elements are topology-inspired, associating velocity degrees of freedom (DOF) with element edges and pressure DOF with the element interior, and result in automatically satisfying the LBB condition in the rigid skeleton limit. (ii) We discretize the skeleton displacement field ($u$) using standard Galerkin finite elements in such a way that together with the fluid pressure field, the LBB condition in the low hydraulic conductivity limit is satisfied in the same manner as in incompressible elasticity problems \cite{bathe2006finite,strang2007computational}. In the context of RT elements, we note other topology-inspired finite elements \cite{BossavitWhitney,DiscreteExteriorCalculus,Arnold_etal_1984}.

Most finite elements of the three-field category for porous media are obtained using standard nodal Galerkin discretization of all three fields \cite{Sandhu_etal_1990, Gajo_etal_1994, Diebels_Ehlers_1996, Levenston_etal_1998, Breuer_1999,Breuer_Jägering_1999,Arduino_Macari_I_2001, Yang_Smolinski_2006,Jeremić_etal_2008,tasiopoulou2014solution}. However, these finite elements usually fail to fulfill the LBB condition. As described above, in this paper, we exploit a mixed finite element approach to deal with the challenges that conventional node-based Galerkin methods are not able to handle easily. The main advantage of mixed elements is their ability to exactly model the kind of continuity which is expected i.e., continuity of the normal velocity across element edges which results in the desirable element-wise (local) mass conservation. Mixed elements form a structure that imitates, at the discrete level, all the operators that exist at the continuous level PDEs (such as divergence and gradient). In addition, boundary conditions can be applied more naturally. For these reasons, mixed elements such as RT elements \cite{Raviart_Thomas_1977}, N\'ed\'elec elements \cite{Nédélec_1980} or a variety of others summarized in \cite{Chen_FiniteElement} are widely used in solid and fluid mechanics problems of incompressible elasticity and fluid flow. However, their application in coupled problems (simultaneous solution) has been limited. Phillips and Wheeler \cite{Phillips_Wheeler_I_2007, Phillips_Wheeler_II_2007} and Lipnikov \cite{Lipnikov_2002} developed a formulation based on a combination of the lowest order Raviart-Thomas ($\rt$) mixed finite element for the $w$ and $p$ fields and the standard nodal linear Galerkin finite element for the $u$ field for quasistatic poroelasticity problems. We call this element $\linprt$. They elaborated on the solvability and error estimates of this finite element for quasistatic consolidation problems. In this paper, we extend the application of this $\linprt$ element to dynamic problems. In doing so, we are able to effectively capture wave propagation effects (section \ref{Example1_halfSpace}). We also observe that careful consideration of the finite element mesh is required to obtain the proper relationship between hydraulic conductivity and the damping of the mechanical response (section \ref{Example2_soilBlock}). This relationship has not been captured correctly by many past results in the literature. Furthermore, we find that while stability in terms of the LBB condition is achieved for typical ranges of hydraulic conductivity, in the limiting case of very small hydraulic conductivity, the $\linprt$ element exhibits instability in the form of checkerboard pattern  in the $p$ field. For reasons described in section \ref{PoroWLowPermeabil}, this effect (also referred to as locking in this context), is more pronounced in dynamic problems than previously observed in quasistatic problems \cite{Phillips_Wheeler_2009}. To resolve this issue, we propose a $\quadprt$ element which employs a quadratic Galerkin finite element to estimate the space of $u$. The element is shown to pass the numerical inf-sup test and therefore the stability criterion (see section \ref{PoroWLowPermeabil}).

This paper is organized as follows. First we summarize Biot's mathematical model for porous media in section \ref{Governingequations}. In section \ref{FiniteElement}, we motivate and present our mixed finite element formulation; show how the poroelasticity equations reduce to two classical models in the limiting case, and discuss the associated stability considerations; we develop our spatial discretization scheme based upon the three-field formulation ($u$-$w$-$p$) and describe the proposed finite element formulation using RT elements for the $w$-$p$ fields and nodal Galerikin element for the $u$ field. In section \ref{TimeDiscritization}, we discuss time integration, in particular addressing the need for consistent initial conditions; we also present the energy balance equation for the full porous medium, which we later use to assess the stability and accuracy of time integration.
We present the details of implementation in section \ref{Implementation}, including aspects involving unconventional edge degrees of freedom. To demonstrate the performance of this approach, we consider detailed numerical examples in section \ref{NumericalExamples}. These include verification with analytical and boundary element solutions, wave propagation analysis, energy balance analyses, mass lumping considerations, effects of mesh pattern and size, hydraulic conductivity effects on damping and frequency, spurious pressure modes and locking, and stability analysis by means of inf-sup tests. We also comment on some discrepancies found in the results in the literature studies.
\section{Governing equations} \label{Governingequations}
\begin{figure}
    \centering
    \subfloat[]
    {
	 \resizebox{0.45\textwidth}{!}{\includegraphics {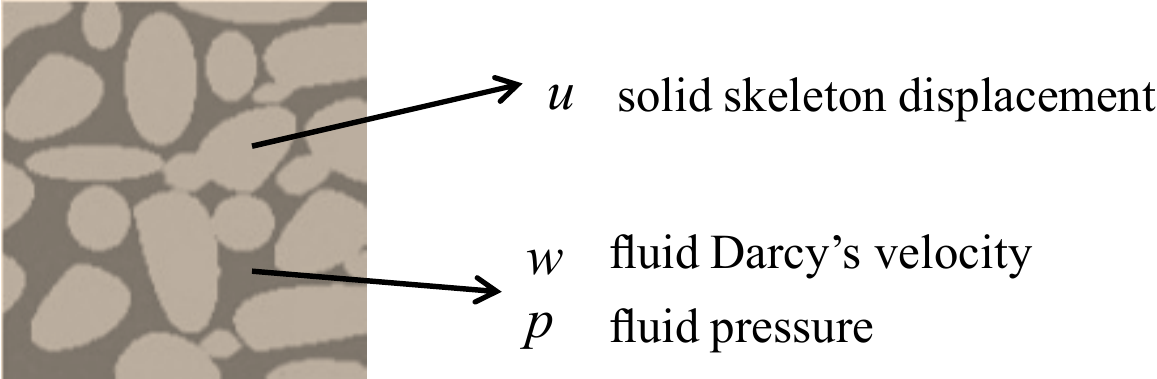}}        
 \label{fig:porousMedia}
    }
    \hspace{0.1\textwidth}
    \subfloat[]
    {
     \resizebox{0.45\textwidth}{!}{\includegraphics {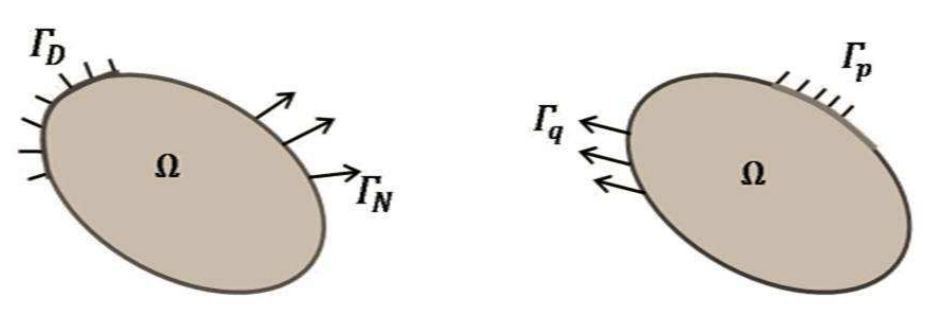}}
        \label{fig:bc}
    }
    \caption{ (a) Saturated porous medium with primary variables. (b) Domain $\Omega$ corresponds to the skeleton boundaries is shown on the left while at the same time,  $\Omega$ corresponds to the fluid boundaries is shown on the right. In each figure boundary $\Gamma$ is divided into two disjoint subsets which cover the entire boundary of the domain; i.e., $\Gamma_D \cup \Gamma_N= \Gamma$ and $\Gamma_p \cup \Gamma_q= \Gamma$. }
\end{figure}

We describe Biot's equations \cite{Biot_lowfrequency_1956} under the following conditions: (i) fluid and skeleton materials are incompressible; (ii) there is a single fluid phase; (iii) skeleton has a linear elastic behavior; (iv) a linear form of Darcy's law (i.e,. constant hydraulic conductivity) is used; (v) a linearized kinematics formulation is applied which also results in constant porosity;  (vi) the additional apparent mass\footnote{also called the added mass; in porous medium, when the fluid (skeleton) is moving, a force must be exerted on the skeleton (fluid) to prevent an average displacement
of the latter which is defined through a mass coupling coefficient.} used in Biot's equations is neglected; (vii) to focus on the essential features of our formulation, we do not consider fluid sources or sinks, body forces or gravity effects, inclusion of which is straightforward. Accordingly, the linearized governing equations read \cite{Chen_Dargush_1995}
\begin{equation} \label{eq:continuum_system}
	\begin{alignedat}{2}
		& \rho\, \ddot u + \rhof\, \dot w  - \divg (\sigma\,-\,p\, I)&=&\,0 \\
		&\rhof \,  \ddot u + \frac{\rhof}{\nf} \, \dot w+  \frac{\rhof \,  g}{ \kh} \,  w+\grad (p ) \, &=&\,0  \\
		&\divg (\dot u ) + \divg ( w )&=&\,0
	\end{alignedat}
\end{equation}
These represent the momentum balance of the porous medium, the dynamic extension of Darcy's law, and the mass balance equation (continuity equation) respectively. The three basic fields are skeleton displacement, $u$, pore fluid pressure, $p$, and fluid Darcy velocity, $w$, as shown in Figure \ref{fig:porousMedia}. Darcy velocity is defined as the relative rate of discharge per unit area of the medium, $w=\nf (\dot u_{\text f}-\dot{u}̇)$, where $\dot u_{\text f}$ refers to the fluid absolute velocity. $\divg$, $\grad$, and dot denote the gradient operator, the divergence operator, and the derivative with respect to time, $t$, respectively. $\rhof$ is the fluid density, while $\rho=\nf \, \rhof+(1-\nf) \rhos$ is the average density of the medium, with $\rhos$ and $\nf$ being solid phase density and porosity respectively. $I$, $g$, and $\kh$ represent the identity matrix, acceleration due to gravity, and the hydraulic conductivity (its associated matrix is assumed to be isotropic).

The concept of effective stress is used as a principle that controls the constitutive behavior of a porous medium: $\sigmat=\sigma\,- \alpha \,p\, I$  where the total stress $\sigmat$ is composed of effective stress of skeleton $\sigma$ and the pore pressure $p$. Here, $\alpha$ denotes the Biot-Willis coefficient which is commonly taken as 1 for soil material. For an isotropic elastic skeleton, the effective stress is given by the usual constitutive equation
\begin{equation} \label{elasticity} 
\sigma= 2 G \, \varepsilon+ \lambda I \, \text{trace}(\varepsilon)
\end{equation}
where $\lambda$ and $G$ are the Lam{\'e}'s first and second constant (shear modulus) of the skeleton under drained condition, respectively. $\varepsilon=\frac{1}{2}(\nabla u+\nabla u ^ \top)$ is the strain.

Equation \eqref{eq:continuum_system} is a system of coupled partial differential equations with $u$, $w$, and $p$ as the three unknown fields on the domain $\Omega$. The following are appropriate boundary conditions
\begin{equation}
\begin{alignedat}{6}
u=u_\prsc  & \quad \text{over  } \Gamma_D & \quad \text{and }&\quad \sigmat \cdot \hat n=\tprsc &\quad \text{over  } \Gamma_N & \quad \text{(skeleton boundary condition)}\\
\quad  w \cdot \hat n=q_\prsc &\quad \text{over  } \Gamma_q & \quad  \text{and } &\quad p=p_\prsc  &\quad  \text{over  } \Gamma_p   & \quad \text{(fluid boundary condition)}
\end{alignedat}
\end{equation}
where $u_\prsc, \tprsc, p_\prsc$ and $q_\prsc$ are the prescribed skeleton displacement, traction, pressure and volumetric flux, and $\hat n$ is the outward unit normal to the boundary $\Gamma$. We note that for each of the flow and skeleton, the boundary is decomposed into two disjoint subsets which cover the entire boundary $\Gamma$ as shown in Figure \ref{fig:bc}.
\nomenclature{$u$, $w$, $p$}{skeleton displacement, Darcy velocity, pressure}
\nomenclature{$\rho$, $\rhof$}{average density of the medium, fluid density}
\nomenclature{$\kh$, $\nf$}{hydraulic conductivity, porosity}
\nomenclature{$\sigma$}{effective stress}
\nomenclature{$\varepsilon$}{strain}
\nomenclature{$\tprsc$}{prescribed traction} 
\nomenclature{$p_\prsc$}{prescribed pressure} 
\section{Finite element formulation} \label{FiniteElement}
The coupled nature of the poroelastic governing equations \eqref{eq:continuum_system} precludes analytical solutions, except in some special cases like those presented in the seminal paper of Rice and Cleary \cite{rice1976some}. In most cases, a numerical approximation must be computed. For this purpose, typically a finite element approach is used for spatial discretization, resulting in a system of differential algebraic equations (DAE) that can be solved by a time integration scheme. In this section, we discuss spatial discretization, specifically our approach of combining RT elements for the $w$-$p$ fields with standard nodal Galerkin elements for the $u$ field. First we motivate our choice of a three-field formulation. Then we discuss the stability concerns that arise. Finally, we present the finite element formulation.
\subsection{Three-field formulation} \label{ThreeFieldFormulation}
Zienkiewicz and Shiomi \cite{Zienkiewicz_Shiomi_1984} first discussed various possible finite element formulations for Biot's dynamic poroelastcity model with the notion of incompressibility. These can be classified into three-field ($u$-$w$-$p$) formulations and two-field ($u$-$p$ or $u$-$w$) formulations. Zienkiewicz and Shiomi \cite{Zienkiewicz_Shiomi_1984} did not present an implementation of a three-field formulation.
A majority of current implementations are based on two-field formulations obtained by eliminating either the fluid Darcy velocity $w$, resulting in a $u$-$p$ formulation \cite{Zienkiewicz_etal_1980,Zienkiewicz_1982,Zienkiewicz_etal_1990,Murad_Loula_1994,Oka_et al_1994,Zhu_Suh _2001,Li_et al_2004,Liu_2004}
or the pressure, leading to a $u$-$w$ formulation \cite{Prevost_1982, Engelman_et al_1982, Zienkiewicz_Shiomi_1984, Prevost_1985,Spilker_Maxian_1990,Suh_etal_1991}. These methods are efficient in terms of the computational cost as they reduce the number of DOFs. However, we pursue here a three-field formulation motivated by the following considerations.

\noindent \textit{Generality} -- When the fluid and the solid material making up the skeleton are modeled as incompressible, the pressure term does not appear in the mass conservation equation \subeq{eq:continuum_system}{3}. Consequently, the fluid Darcy velocity, $w$, cannot be directly eliminated from the Darcy's law \subeq{eq:continuum_system}{2} \cite{Zienkiewicz_Shiomi_1984,Gajo_etal_1994}. Therefore, to obtain a $u$-$p$ formulation, either (a) penalty procedures must be used in the mass conservation equation \cite{Prevost_1982, Engelman_et al_1982, Zienkiewicz_Shiomi_1984, Prevost_1985}, or (b) the fluid relative acceleration, $\dot{w}$, (and hence inertia) must be neglected. The latter is a simplification, first introduced by Zienkiewicz and coworkers \cite{Zienkiewicz_etal_1980,Zienkiewicz_Shiomi_1984}, has been the most common means of obtaining a $u$-$p$ formulation \cite{Zienkiewicz_etal_1980, Zienkiewicz_1982,Zienkiewicz_etal_1990, Murad_Loula_1994, Oka_et al_1994, Zhu_Suh _2001, Li_et al_2004, Liu_2004}. This approach has been found to be typically valid for the slow and medium frequency level phenomena \cite{Zienkiewicz_et al_1980, Zienkiewicz_1982}.  
However, in a variety of problems, it is desirable to take into account all acceleration terms; which is only possible by exploiting the three-field formulation, without using a penalty formulation. For instance using verification and validation procedures, Jeremi{\'c} and coworkers \cite{Jeremić_etal_2008,tasiopoulou2014solution,tasiopoulou2014Validation} show that a three-field formulation is required to realistically model the saturated soil-foundation-structure interaction in earthquake as well as the liquefaction and cyclic mobility phenomena. They conclude that a full-field method can appropriately model the physical damping in the dynamic poroelasticity problem during seismic events by directly taking into account the fluid and skeleton interaction through Darcy's equation. Accordingly, we could model and compute the velocity proportional damping energy as shown in equation \eqref{eq:energy4} whereas in other studies an artificial damping has to be used \cite{Li_et al_2004}. Another area of interest is soft biological tissues that are often modeled mechanically as poroelastic media, and experience very rapid external forces \cite{Moeendarbary2013}. This situation creates significant inertia forces in the fluid phase that is not negligible and requires use of a full formulation \cite{Yang_2006}. An improved performance of the full-field mixed finite element formulation over the two-field formulation is also reported in \cite{Levenston_et al_1998, Zhu_Suh_2001} for large deformation simulations and in \cite{Turan_Arbenz_2013} for modeling of bones subjected to Osteoporosis disease.
Presence of the fluid velocity field in an explicit fashion will also allow the formulation to be more readily extended for coupled transport phenomena.

\noindent \textit{Accuracy} -- The $u$-$p$ formulation amounts to solving for pressure and skeleton displacement as primary variable fields. Typically post-processing techniques are required to compute the velocity field from the gradient of the pressure field, and is therefore of lower accuracy \cite{Simon_etal_higherOrder_1986,Gajo_etal_1994}. Furthermore, both pressure and flux boundary conditions can be applied directly with a three-field formulation, but require special treatment if a two-field formulation is used \cite{gresho1987pressure}.
\subsection{Nature of stability considerations} \label{LimitingCases} 
To understand the nature of stability issues that are relevant in a three-field finite element formulation, we describe two limiting conditions in terms of hydraulic conductivity and rigidity of skeleton similar to \cite{Wan_2002,White_2009}. Both of these limiting cases give rise to saddle point problems. Without loss of generality, for these limiting cases we only consider the static version of the poroelasticity problem where all the inertia terms are zero\footnote{The dynamic version affects only on ellipticity condition.}.

Case (i) -- Low hydraulic conductivity: In the impermeable limit $(\kh \to 0)$, the static version of \eqref{eq:continuum_system} reduces to
\begin{equation} \label{Eq:ElasticitySaddle}
\begin{alignedat}{3}
&2\, G \, \nabla ^2 u &\,- \, \nabla p&\,=\,0 \\
& \divg u&&\,=\,0
\end{alignedat}
\end{equation}
where $\nabla ^2$ is the Laplacian operator. Note $\divg \dot u=0$ is  replaced by $\divg  u=0$ as they are equivalent if the initial conditions satisfy. We recognize that these are the equations of the classical incompressible elasticity problem. Thus in the impermeable limit, there is no flow and the skeleton becomes incompressible but deformable.

Case (ii) -- Rigid skeleton: The other limiting case is when the skeleton matrix become rigid (which implies $ u \to 0 $) and results in
\begin{equation} \label{Eq:DarcySaddle}
\begin{alignedat}{3}
&w&\,+\,\frac{\kh}{\rhof g} \, \nabla p&\,=\,0 \\
& \divg w&&\,=\,0
\end{alignedat}
\end{equation}
These are the classical Darcy's flow equations. We note that this second limiting case applies only to a three-field formulation.

Both these limit cases are saddle-point problems \cite{braess2001finite}. The coefficient matrix associated with the discrete form of the saddle point problem has a zero-block structure on lower diagonal (see equation \eqref{incompressElasticity} for example) that creates stability challenges. There are two stability criteria that provide the necessary and sufficient condition for well-posedness of such problems, namely the ellipticity condition and the LBB condition (also called inf-sup condition). We refer to \cite{braess2001finite, Brezzi_Fortin_1991} for a thorough description on this subject.
The discrete version of the ellipticity condition is satisfied for any choice of finite element with positive definite mass and stiffness matrices upon application of boundary conditions. The LBB stability condition is more important and challenging. Often, standard finite element approaches presented in the literature fail to satisfy this condition, leading to instability manifested as checkerboard patterns in the pressure field poor convergence behavior. We provide a comparison of stable and unstable finite element results in section \ref{sec:NumericalExmpl3}.

Since the poroelasticity model degenerates to the saddle point problems \eqref{Eq:ElasticitySaddle} and \eqref{Eq:DarcySaddle}, a finite element formulation must satisfy the stability conditions associated with both these limiting cases. Thus the poroelasticity problem has the additional requirement that \emph{both} the $u$-$p$ pair and the $w$-$p$ pair must be able to satisfy the stability conditions simultaneously. We show in section \ref{sec:inf-sup} that the proposed $\quadprt$ element qualifies.

\subsection{Finite element discretization} \label{FEDescription}
%
%
We start from the weak-form representation of \subeq{eq:continuum_system},
\begin{equation} \label{eq:weakform1}
	\begin{aligned}
		&\rho\ \int_\Omega { \delta u}^\top \ddot {u} \, d\Omega\,+ \,\rhof \int_\Omega { \delta u}^\top \dot {w}   \, d\Omega\,+\,  \int_\Omega ({ \grad \delta u})^\top (\sigma\,-\,p\, I)   \, d\Omega
		 \, -\,    \int_{\Gamma_N} { \delta u}^\top \, \tprsc \, d\Gamma\ +\\
		&\rho_f  \int_\Omega { \delta w}^\top  {\ddot u} \, d\Omega \, +\, \frac{\rho_f}{\nf} \int_\Omega { \delta w}^\top  {\dot w} \, d\Omega \,+\,  \frac{\rhof \, g}{\kh}
		\int_\Omega { \delta w}^\top w \, d\Omega \,- \, \int_\Omega  (\divg w)^\top p \, d\Omega \, + \,  \int_{\Gamma_p}  (w^\top\cdot\hat n) p_\prsc  \, d\Gamma +\\
		&\int_\Omega { \delta p}^\top \divg (\dot u)  \, d\Omega \,+\,   \int_\Omega { \delta p}^\top \divg (w) \, d\Omega \, =\,0
	\end{aligned}
\end{equation}
where $\delta u$, $\delta w$, and $\delta p$ are the skeleton displacement, fluid relative velocity, and pressure test functions respectively. We recognize now that the functions $u$, $w$, and $p$ must belong to the spaces
\begin{equation}\label{eq:hilbertSpaces}
\begin{aligned}
	& \mathcal U=\{ u \in (H^1 (\Omega))^2 : u= u_\prsc  \text{ over } \Gamma_D \} \\
	&  \mathcal W=\{ w \in H(\text{div}; \Omega) : w \cdot \hat n= q_\prsc  \text{ over } \Gamma_q \}  \\
	& \mathcal P=\{ p \in L^2(\Omega) \}
\end{aligned}
\end{equation}
respectively. Here $L^2(\Omega)$ is the space of square integrable functions on $\Omega$, $H^1 (\Omega)$ is the usual Sobolev space in which each function and its first order derivative belong to $L^2(\Omega)$, and $H(\text{div}; \Omega)$ represents the space of vector-valued functions that belong to $L^2(\Omega)$ and whose divergences belong to $L^2(\Omega)$ as well. Further elaboration on these spaces can be found in \cite{Brezzi_Fortin_1991} and also in \cite{Raviart_Thomas_1977, Lipnikov_2002,Phillips_Wheeler_I_2007}.
\nomenclature{$L^2 (\Omega)$}{square integrable functions in $\Omega$}
\nomenclature{$H^1 (\Omega)$}{$H^1 (\Omega)=\{ u: \, u \in   L^2(\Omega),\, Du \in  L^2(\Omega) \} $  }
\nomenclature{$H(\text{div}; \Omega)$}{$H(\text{div}; \Omega)=\{ w: \, w \in   (L^2(\Omega))^2,\, \divg w \in  L^2(\Omega) \} $ }
Similarly, the test functions of $\delta u$, $\delta w$, and $\delta p$ belong to the same spaces as their corresponding variable fields, except that instead of $u_\prsc$ and $q_\prsc$ a zero is placed. We call these test spaces $\mathcal U_0$, $\mathcal W_0$, and $\mathcal P_0$. 

As discussed above, with the goal of satisfying the appropriate stability conditions, we simultaneously apply the lowest order Raviart-Thomas mixed finite element ($\rt$) \cite{Raviart_Thomas_1977} for flow variable fields and the standard linear/quadratic Galerkin finite elements for skeleton variable fields. Let the domain $\Omega$ be partitioned with $\mathcal {\hat{T}}  =   \{{\mathcal {T}}_1,\cdots,{\mathcal {T}}_{\Nel}\}$ into $\Nel$ non-overlapping triangular elements. We denote by $\Nn$ and $\Nfc$ number of nodes and edges. Also position vector is $\bld x=\begin{bmatrix} x& y \end{bmatrix}^\top$. 

\begin{figure}
    \centering  
        \resizebox{0.8\textwidth}{!}{\includegraphics{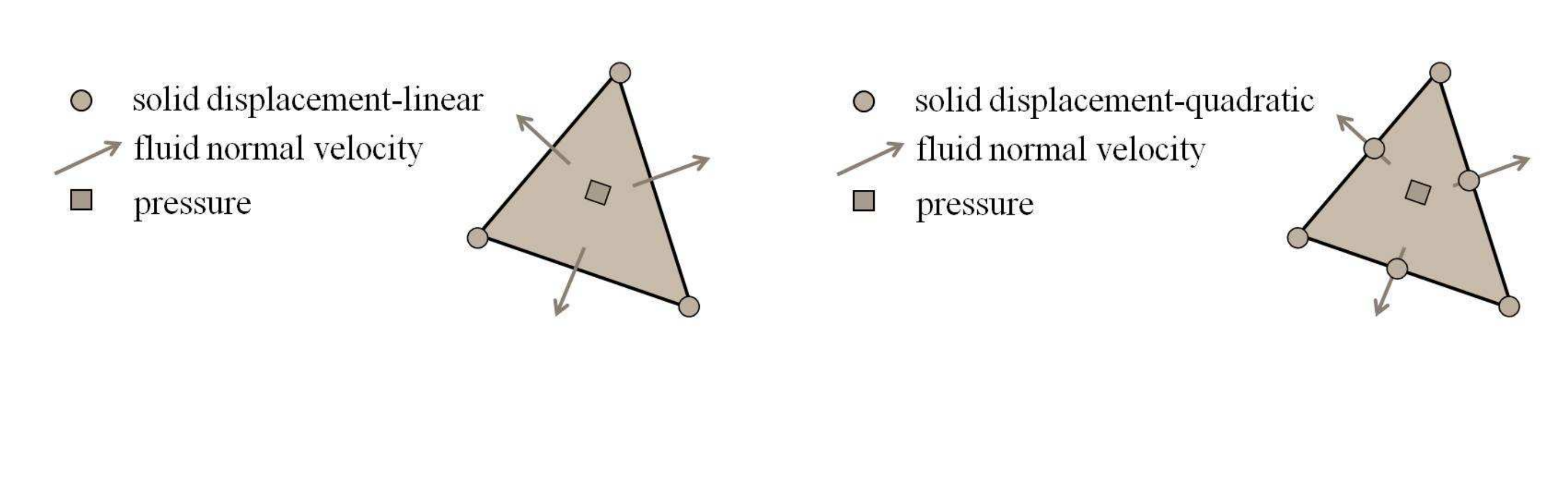}}
    \caption{Fluid and skeleton DOFs in the coupled Galerkin/$\rt$ elements -- (a) linear Galerkin/$\rt$ element, referred as $\linprt$, (b) quadratic Galerkin/$\rt$ element referred as $\quadprt$. Both elements are constructed on the natural topology of the variable fields, i.e., skeleton displacement fields are defined on the nodes using a standard Galerkin approximation, while fluid variable fields are defined on the element edges (velocity) and volumes (pressure) using $\rt$ element. We note that Darcy’s velocity DOFs are the values of the normal velocity at the midpoints of edges of elements and the pressure DOFs are set to be the constant elemental pressure.}
   \label{fig:mixedElement}
\end{figure}

\textbf{Skeleton displacement discretization}: We discretize the skeleton displacement field using standard linear/quadratic shape functions over triangles, $\tilde N_i(\bld x)$,
\begin{equation} \label{eq:interpFunc1}
	\hat u(\bld x,t)= \sum_{i=1}^{\Nn} \tilde N_i(\bld x) \tilde{u}_i(t)= N(\bld x)  \bld{u}(t) 
\end{equation}
We note that the components of vector $\bld{u}(t)$ are skeleton displacement DOFs, $\tilde{u}_i$, that are displacements at the element nodes. The skeleton and fluid DOFs are illustrated in Figure \ref{fig:mixedElement}. We call the combination of linear Galerkin finite element with $\rt$ finite element, $\linprt$, and the combination of quadratic Galerkin finite element with $\rt$ finite element, $\quadprt$. We propose the latter finite element to resolve the stability issues that the linear element creates in the limiting case of very small hydraulic conductivity as elaborated in section \ref{sec:NumericalExmpl3}.

\begin{figure}
    \centering  
    \subfloat[]
    {
        \resizebox{0.185\textwidth}{!}{\includegraphics{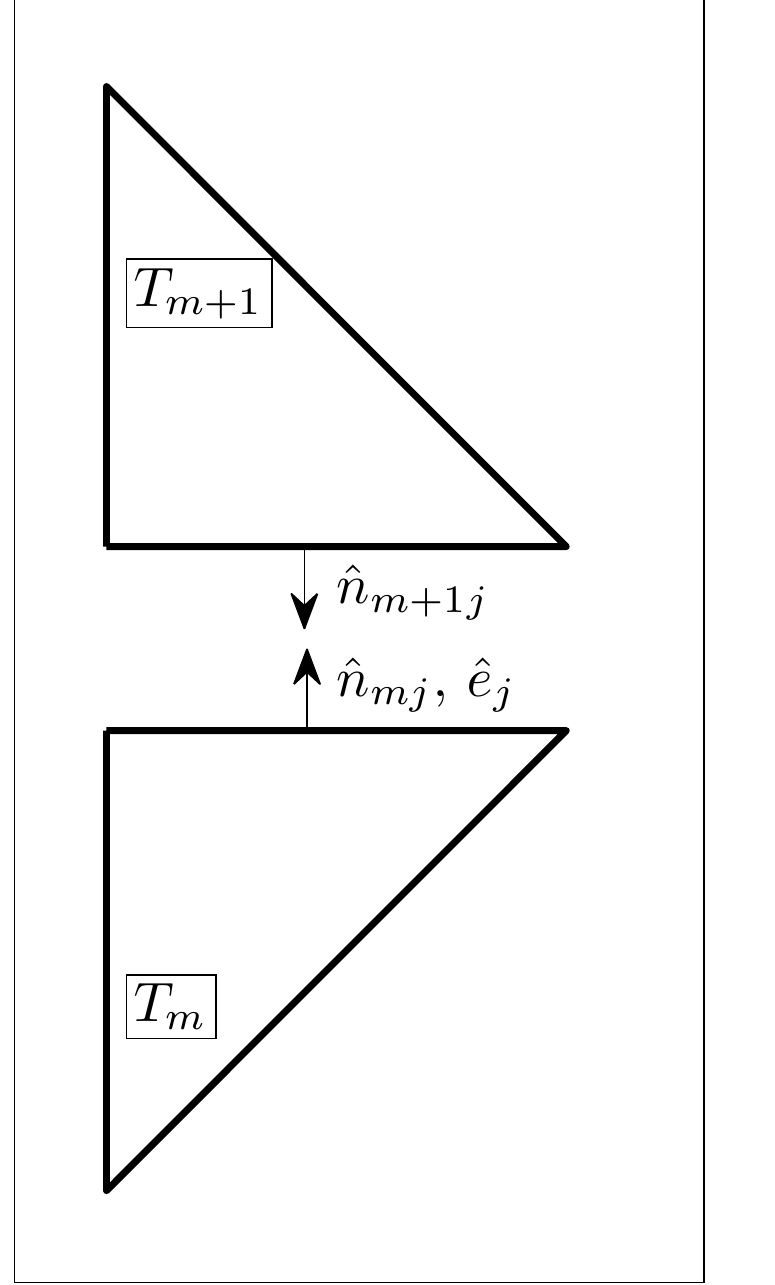}}
        \label{fig:neighberTrinagles}
    }
  \hspace{0.05\textwidth}
    \subfloat[]
    {
        \resizebox{0.17\textwidth}{!}{\includegraphics{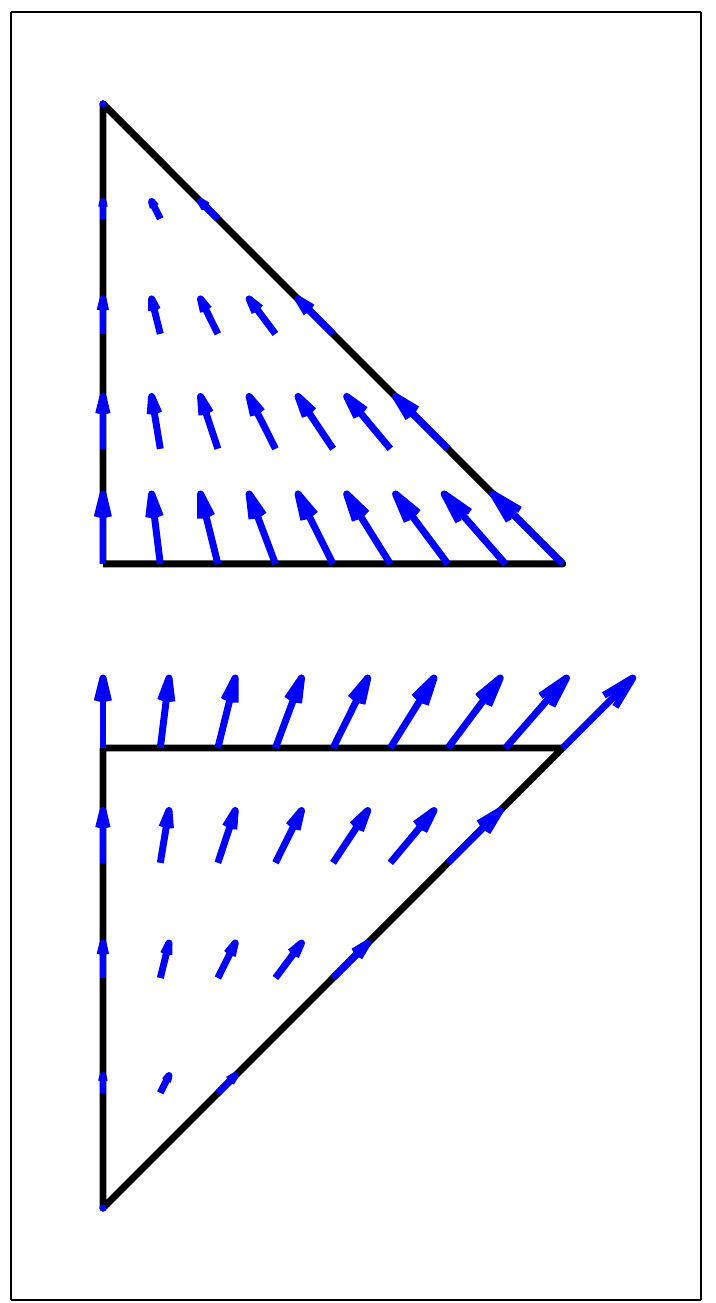}}
        \label{fig:interpol_vel}
    }
    \hspace{0.05\textwidth}
    \subfloat[]
    {
        \resizebox{0.17\textwidth}{!}{\includegraphics{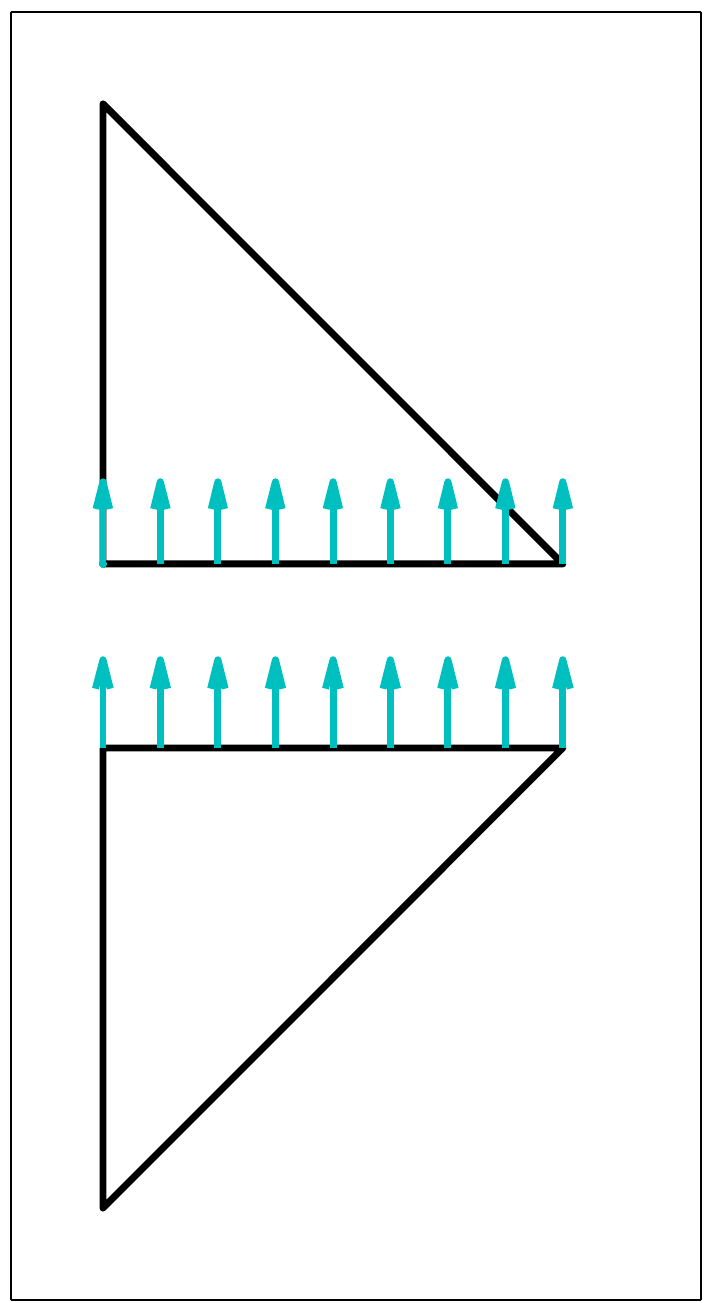}}
        \label{fig:interpol_vel_normal}
    }
    \caption{(a) Two enumerated adjacent elements $T_m$ and $T_{m+1}$ that share the same edge $j$. In this figure, $\hat n_{mj}$ and $\hat n_{m+1 j}$ are the outward unit normal vector to edge $j$ of elements $m$ and $m+1$, respectively, while $\hat e_j$ is
the unit normal vector to edge $j$ chosen with a global fixed orientation. This is such that $\hat e_j$ is equal to the outward unit normal vector of element with smaller index, i.e., $T_m$ (and hence points into $T_{m+1}$). Also based on this convention, the element with smaller index $T_m$ has $s_{mj} = 1$ whereas the one with larger index $T_{m+1}$ has $s_{m+1 j} = -1$. (b) velocity interpolation function for the common edge of two adjacent elements, $\tilde W_j(\bld x)$. (c) normal velocity along the common edge of two adjacent elements. Note the continuity of this component of velocity across element edge which results in element-wise mass conservation.}
    \label{fig:interpol_velocity}
\end{figure}

\textbf{Fluid velocity discretization}: 
\begin{equation} \label{eq:interpFunc2}
	\hat w (\bld x,t)= \sum_{j=1}^{\Nfc} \tilde W_j(\bld x) \tilde q_j(t)= W(\bld x) \bld q(t) , \qquad \tilde W_j(\bld x)=	 \begin{cases}
      			s_{mj}\, \frac{l_j}{2 \mathcal A_m}\, (\bld x-\bld x_j)& \text{if }  \bld x \in \mathcal {T}_m \\
      			0 & \text{otherwise}
    			\end{cases}
\end{equation}
Based on $\rt$ element, Darcy velocity DOFs $\tilde q_j$ (Figure \ref{fig:mixedElement}) are values of the normal velocity at the midpoints of edges of element $m$ (in fact, normal velocity is constant along each edge), defined as $\tilde q_j=\hat w_j \cdot \hat n_{mj}$ where $\hat w_j$ is fluid velocity vector defined on edge $j$. Also $\hat n_{mj}$ denotes the outward unit normal vector to edge $j$ of element $m$, while $\hat e_j$ is the unit normal vector to edge $j$ chosen with a global fixed orientation. This global convention is such that $\hat e_j$ points from element $\mathcal {T}_m$ with smaller index into its adjacent element with larger index (see Figure \ref{fig:neighberTrinagles}). Therefor, based on this convention, edges on the boundaries always have an outward global unit normal. 
In equation \eqref{eq:interpFunc} $\bld x_j$ is the position vector of the vertex node opposite the $j^{\text{th}}$ edge of element $m$ and $s_{mj}=\hat n_{mj} \cdot \hat e_j$ is a sign convention equal to $+1$ if $\hat e_j$ points outward and otherwise is $-1$. The effect of this sign convention can be seen in Figure\ref{fig:interpol_vel} where the velocity interpolation function, $\tilde W_j(\bld x)$, for the common edge ($j$) of the two adjacent elements is illustrated. Due to the convention, the element with smaller index $\mathcal {T}_m$ has $s_{mj}=1$ and the one with larger index $\mathcal {T}_{m+1}$ has $s_{m+1j}=-1$. Moreover, shown in Figure \ref{fig:interpol_vel_normal} is the normal component of velocity along the common edge of two adjacent elements illustrating the continuity of such normal component. Also $\mathcal A_m$ and $l_j$ refer to the area of element $m$ and the length of edge $j$, respectively.

\textbf{Pressure discretization}: Pressure DOF $\tilde p_m$ is set to be constant elemental pressure, thus its interpolation function $\tilde H_m(x)$ is a polynomial of degree 0 that takes a constant value of 1 at element $m$ and 0 at all other elements,
\begin{equation} \label{eq:interpFunc}
	\hat p (\bld x,t)= \sum_{m=1}^{\Nel} \tilde H_m(\bld x) \tilde{p}_m(t)= H(\bld x) \bld{p}(t) , \qquad \tilde H_m(\bld x)=	 \begin{cases}
      			1 & \qquad \, \text{if }  \bld x \in \mathcal {T}_m\\
      			0 & \qquad \,  \text{otherwise}
    			\end{cases} 
\end{equation}
Fluid Darcy velocity and pressure DOFs are collected in vector $\bld q(t)$ and $\bld p(t)$, respectively and are shown in Figure \ref{fig:mixedElement}.
\nomenclature{$\bld u$}{vector of skeleton displacement DOFs}
\nomenclature{$\bld p$}{vector of pressure DOFs}
\nomenclature{$\bld q$}{vector of Darcy velocity DOFs}
\nomenclature{$\linprt$}{linear Galerkin finite element coupled with $\rt$ finite element}
\nomenclature{$\quadprt$}{quadratic Galerkin finite element coupled with $\rt$ finite element}
\nomenclature{$\rt$}{Lowest order Raviart-Thomas element}
\nomenclature{$\hat w_j$}{velocity vector defined on edge $j$}
\nomenclature{$\tilde W_j(\bld x)$}{velocity interpolation function for edge $j$}
\nomenclature{$\hat n_{mj}$}{outward unit normal vector to edge $j$ of element $m$}
\nomenclature{$\hat e_j$}{unit normal vector to edge $j$ chosen with a global fixed orientation}
\nomenclature{$s_{mj}$}{sign convention equals to $\hat n_{mj} \cdot \hat e_j$}
\nomenclature{$\mathcal A_m$}{area of element $m$}
\nomenclature{$l_j$}{length of edge $j$}
\nomenclature{$N(\bld x)$, $W(\bld x)$, $H(\bld x)$}{skeleton displacement, fluid velocity, and pressure interpolation functions}

Given the described finite element approximation, we formulate a discrete version of the weak problem \eqref{eq:weakform1} by substituting approximation functions of $\hat u,\hat w ,\hat p$ in it.
\footnotesize
\begin{equation} \label{eq:weakform}
	\begin{aligned}
		&\rho\ \int_\Omega \bld{ \delta u}^\top N^\top N \ddot {\bld u} \, d\Omega\,+ \,\rhof \int_\Omega \bld{ \delta u}^\top N^\top W \dot {\bld q}   \, d\Omega\,+\,  \int_\Omega \bld{ \delta u}^\top  \mathcal B^\top \mathcal C \mathcal B \, \bld u \, d\Omega
		\,- \,\int_\Omega \bld{ \delta u}^\top  \mathcal B^\top \check{\delta} \, H \, \bld p \, d\Omega  \, =\,    \int_{\Gamma_N} \bld{ \delta u}^\top N^\top \, \tprsc \, d\Gamma\ \\
		&\rho_f  \int_\Omega \bld{ \delta q}^\top  W^\top  N \bld {\ddot u} \, d\Omega \, +\, \frac{\rho_f}{\nf} \int_\Omega \bld{ \delta q}^\top  W^\top W \bld {\dot q} \, d\Omega \,+\,  \frac{\rhof \, g}{\kh}
		\int_\Omega \bld{ \delta q}^\top  W^\top W \bld q \, d\Omega \,- \, \int_\Omega \bld{ \delta q}^\top  (\divg W)^\top H \bld  p \, d\Omega \, =\, - \int_{\Gamma_p}  \bld{ \delta q}^\top  (W^\top\cdot\hat n) p_\prsc  \, d\Gamma\\
		&\int_\Omega \bld{ \delta p}^\top H^\top \check{\delta}^\top B \, \bld {\dot u}  \, d\Omega \,+\,   \int_\Omega \bld{ \delta p}^\top H^\top (\divg W) \bld q \, d\Omega \, =\,0
	\end{aligned}
\end{equation}
\normalsize
The solution of the discretized problem must satisfies the integral equations for all test functions $(\delta{ \hat u},\,  \delta {\hat w}, \, \delta {\hat p}) \in \mathcal{\hat U}_0 \times \mathcal{\hat W}_0 \times \mathcal{\hat P}_0$. In equation \eqref{eq:weakform} $\check{\delta}$ denotes the Kronecker delta in the vector format, $\mathcal B(x)$ is the linearized deformation-displacement matrix, and $\mathcal C$ is the elastic constitutive matrix under the drained condition. Note that we dropped the dependence of functions to $x$ and $t$ due to clarity. The test functions and the shape functions coincide in this formulations. 
Moreover, equation \eqref{eq:weakform} represent a consistent variational formulation in which the same shape functions are applied for higher time derivatives. An alternative formulation is discussed in section \ref{sec:massRepresentation}. The prescribed skeleton displacement and fluid flux are imposed in a strong way, i.e., as essential boundary conditions, while the prescribed traction ($\tprsc$) and pressure ($ p_\prsc$) are entered the variational form as natural boundary conditions. Also, the constitutive equation \eqref{elasticity} is used in strong form as substitute in \subeq{eq:continuum_system}{1}.
The compact form of \eqref{eq:weakform} is the following semi-discrete expression
\begin{equation} \label{eq:spaceDiscritzed} 
	\begin{alignedat}{5} 
		& M \bld{\ddot u} &\,+&  \Mf \, \bld{\dot q}  & \,+\,&  K \, u &\,- \,&Q  \, \bld p   \, &=&\,  \mathbb{P}  \ \\
		& { \Mf}^\top  \bld{\ddot u}&\,+\,& A \bld{\dot q} &\,+\,&  \frac{n \, g}{\kh}\,A \bld q&\,-\, &B \, \bld p  \, &=&\, \mathbb{F}\\
		&Q^\top \bld {\dot u}  &\,+\,&  B^\top \bld q  \, &&&&&=& \,0
	\end{alignedat}
\end{equation}
where the following terms are defined
\begin{equation} \label{eq:Coef_matrix}
\begin{alignedat}{3}
& M=\rho \int_\Omega N^\top N  \, d\Omega &\qquad& A=\frac{\rhof}{n} \int_\Omega W^\top W  \, d\Omega\\
& \Mf=\rhof \int_\Omega  N^\top W  \, d\Omega&\qquad& K=\int_\Omega \mathcal B^\top \mathcal C \mathcal B  \, d\Omega\\
& Q=\int_\Omega \mathcal B^\top \check\delta \, H  \, d\Omega &\qquad& B=\int_\Omega (\divg W)^\top H  \, d\Omega\\
& \mathbb{P}= \int_{\Gamma_N}  N^\top \, \tprsc \, d\Gamma& \qquad&\mathbb{F}= - \int_{\Gamma_p}   (W^\top\cdot\hat n) p_\prsc  \, d\Gamma
\end{alignedat}
\end{equation}
We reiterate that such topology-based finite element discretization is physically compatible, and as a result all the above constant matrices (at the discrete level) imitates their corresponding operators at the continuous level in \eqref{eq:continuum_system}; for instance $B$ is the discrete version of divergence operator whose $m^{\text{th}}$ column represents the $m^{\text{th}}$ element where all entries are zeros except the rows representing the free edges belongs to that element. In those entries depending on the flux convention, $\pm 1$ is placed. Further elaboration on building these matrices is provided in Section \ref{Implementation}. Also using the described finite element, specification of pressure and flux boundary conditions are straightforward.
\section{Time Discretization} \label{TimeDiscritization}
The resulting system of equations \eqref{eq:spaceDiscritzed} not only consist of ordinary differential equations (ODE), but also algebraic equations. The numerical time integration of a system of Differential Algebraic Equations (DAE) is more complex than the solution of an ODE. Therefore, we utilize the Newmark constant average time integration scheme with a special attention to select a consistent initial condition (Section \ref{ConsistentInitial}). Based on the Newmark method, governing equations of \eqref{eq:spaceDiscritzed} is considered at $t_{n+1}$
\begin{equation} \label{eq:TimeDiscritzed1} 
	\begin{aligned}
& M \bld{\ddot u_{\text{n+1}}} &\,+&  \Mf \, \bld{\dot q_{\text{n+1}}}  & \,+\,&  K \, \bld{u_{\text{n+1}}} &\,- \,&Q  \, \bld{p_{\text{n+1}}}   \, &=&\,  \mathbb{P}_{n+1}  \\
		& { \Mf}^\top  \bld{\ddot u_{\text{n+1}}}&\,+\,& A \, \bld{\dot q_{\text{n+1}}} &\,+\,&  \frac{n \, g}{\kh}\,A \, \bld{q_{\text{n+1}}}&\,-\, &B \, \bld{p_{\text{n+1}}} \, &=&\, \mathbb{F}_{n+1}\\
		&Q^\top \bld{\dot u_{\text{n+1}}}  &\,+\,&  B^\top \bld{q_{\text{n+1}}}  \, &&&&&=& \,0
	\end{aligned}
\end{equation}
in which the following assumptions are made to approximate the values at each time step
\begin{equation} \label{eq:newmark1} 
	\begin{alignedat}{3}
&\bld{u_{\text{n+1}}}=\bld{\tilde{u}_{\text{n}}}+\frac{\Delta t}{2}\bld{\dot u_{\text{n+1}}} \qquad &
&\bld{\ddot u_{\text{n+1}}}=\bld{\tilde{\ddot u}_{\text{n}}}+\frac{2}{\Delta t}\bld{\dot u_{\text{n+1}}} \qquad&
&\bld{\dot q_{\text{n+1}}}=\bld{\tilde{\dot q}_{\text{n}}}+\frac{2}{\Delta t}\bld{q_{\text{n+1}}}
	\end{alignedat}
\end{equation}
where $\Delta t=t_{n+1}-t_n $ and 
\begin{equation*} 
	\begin{alignedat}{3}
&\bld{\tilde{u}_{\text{n}}}=\bld{u_{\text{n}}}+\frac{\Delta t}{2}\bld{\dot u_{\text{n}}} \qquad 
&\bld{\tilde{\ddot u}_{\text{n}}}=-(\bld{\ddot u_{\text{n}}}+\frac{2}{\Delta t}\bld{\dot u_{\text{n}}})\qquad 
&\bld{\tilde{\dot q}_{\text{n}}}=-(\bld{\dot q_{\text{n}}}+\frac{2}{\Delta t}\bld{q_{\text{n}}})
	\end{alignedat}
\end{equation*}
using approximations of \eqref{eq:newmark1} in the discretized equations of \eqref{eq:TimeDiscritzed1} results in the final time discretized version
\begin{equation} \label{eq:TimeDiscritzed2} 
	\begin{alignedat}{7}
& \bar{M} \, \bld{ \dot u}_{\text{n+1}}&\,+&  {M}_f  \, \bld q_{n+1}  & \,- &\, \bar{Q} \,   \bld{p}_{\text{n+1}}  \, &=\, & \bar{\mathbb{P}}_{\text{n+1}}   \\
		& {{M}_f}^\top \, \bld{\dot u}_{\text{n+1}} &\,+&\bar{A}  \, \bld q_{n+1}  &\,- &\, \bar{B} \,  \bld{p}_{\text{n+1}}  \, &=&\,  \bar{\mathbb{F}}_{\text{n+1}} \\
		&\-\bar{Q}^\top \, \bld{ \dot u}_{\text{n+1}}  &\,+&  \bar{B}^\top \, \bld q_{n+1}&& & =&\,0
\end{alignedat}
\end{equation}
where the following definitions are made
\begin{equation} \label{eq:TimeDiscritzed3} 
	\begin{aligned}
		& \bar{M}=M+\frac{\Delta t^2}{4} K \, ;\qquad  \bar{A}=(1+\frac{n \, g}{\kh}\frac{\Delta t}{2}) A \\
		& \qquad \bar{Q}=Q \, \frac{\Delta t}{2} \, ;\qquad \bar{B}=   B \, \frac{\Delta t}{2} \\
		&\bar{\mathbb{P}}_{\text{n+1}}=\frac{\Delta t}{2} \, (\mathbb{P}_{n+1}-M \bld{\tilde{\ddot u}_{\text{n}}}-  \Mf \, \bld{\tilde{\dot q}_{\text{n}}} -  K \bld{\tilde{u}_{\text{n}}} )\\
		&\bar{\mathbb{F}}_{\text{n+1}}=\frac{\Delta t}{2} \, (\mathbb{F}_{n+1}-{ \Mf}^\top \bld{\tilde{\ddot u}_{\text{n}}}-  A \, \bld{\tilde{\dot q}_{\text{n}}})
	\end{aligned}
\end{equation}
Rearranging and reverting  Equation \eqref{eq:TimeDiscritzed2} to matrix form yields  
\begin{equation}    \label{Matrix form} 
   \begin{bmatrix}   \bar{M} &{M}_f  & -\bar{Q} \\ 
			 {{M}_f}^\top &  \bar{A} & -\bar{B} \\ 
			-\bar{Q}^\top &  -\bar{B}^\top & 0 \end{bmatrix}
			 \begin{bmatrix} \bld {\dot u}_{\text{n+1}}\\ \bld q_\text{n+1} \\  \bld{p}_{\text{n+1}} \end{bmatrix}
			= \begin{bmatrix} \bar{\mathbb{P}}_{\text{n+1}}\\ \bar{\mathbb{F}}_{\text{n+1}} \\ 0 \end{bmatrix}
\end{equation}
$0$ represent a zero matrix of appropriate dimension. Above system of equations is symmetric with diagonal blocks $\bar{M}$ and $\bar{A}$ both symmetric positive definite.
\subsection{Consistent initial condition} \label{ConsistentInitial} 
Following the three-field spatial approximation of the problem, the incompressibility condition, which is embedded in the mass conservation equation \subeq{eq:spaceDiscritzed} {3}, forms an algebraic constraint. The presence of this constraint, makes \eqref{eq:spaceDiscritzed} represent a DAE system in time rather than an ODE in time. This fact is also emphasized by Diebels and Ehlers \cite{Diebels_Ehlers_1996}. We note the absence of the algebraic variable (pressure) from the constraint \subeq{eq:spaceDiscritzed} {3}. We also note that not only the same matrix block $\begin{bmatrix} -\bar{Q}\\ -\bar{B}\end{bmatrix}$ appears in \subeq{Matrix form}{1,2} and \subeq{Matrix form}{3}, but also the Schur complement of the upper left block matrix
$$ - \begin{bmatrix} -\bar{Q}^\top &  -\bar{B}^\top  \end{bmatrix}^\top \begin{bmatrix} \bar{M} &{M}_f \\ {{M}_f}^\top &  \bar{A}\end{bmatrix}^{-1} \begin{bmatrix} -\bar{Q}\\ -\bar{B}\end{bmatrix}$$
is a nonsingular matrix with a bounded inverse. As a result, \eqref{Matrix form} yields to a Hessenberg index-2 system of DAEs as indicated by Ascher and Petzold \cite{Ascher:1998:CMO:551054} and has the form of an augmented Lagrangian equation. Consistency of initial conditions is one of the key aspects in DAEs\footnote{An alternative approach to deal with incompressible poroelasticity system of DAEs is to utilize the energy preserving method proposed by Bauchau \cite{Bauchau_2003}. His approach satisfies the algebraic constraint at the end point (if written in terms of displacement), while imposes the algebraic variable at the midpoint; hence, it does not require an initial condition on the algebraic variable. This method reduces to the constant average Newmark method over unconstrained problems and it preserve the energy of the system by eliminating of the work done by algebraic variable. The result of our Numerical Examples in Section \ref{NumericalExamples} using the constant average Newmark method and Bauchau's method perfectly match.}; the initial values for variables in a DAE must satisfy not only the original equations in the system but also their differentials with respect to time \cite{pantelides1988consistent,Brenan_etal_1989}. Accordingly, to select consistent initial conditions for acceleration terms and pressure, we must satisfy the original equations in the system as well as the differential of the mass conservation equation with respect to time, at time zero which reads
\begin{equation} \label{eq:constInital} 
	\begin{alignedat}{5} 
		& M \bld{\ddot u_0} &\,+&  \Mf \bld{\dot q_0}   &\,- \,&Q  \, \bld {p_0}   \, &=&\,  \mathbb{P}_0  \ \\
		&  { \Mf}^\top  \bld{\ddot u_0}&\,+\,&A \bld{\dot q_0} &\,- \,&B \, \bld {p_0}  \, &=&\, \mathbb{F}_0\\
		&Q^\top \bld {\ddot u_0}  &\,+\,&  B^\top \bld{\dot q_0}  \, &&&=& \,0
	\end{alignedat}
\end{equation}
Solving this equation, $\bld{\ddot u_0}, \bld{\dot q_0}, \bld {p_0}$ is obtained. Note that in deriving \ref{eq:constInital}, fluid and solid phases are assumed to be initially at rest $(\bld{\dot u_0}=0, \bld{q_0}=0)$ without any deformation $(\bld{u_0}=0)$, which is an intuitive set of initial conditions that also satisfy the mass balance equation at time zero. We observe that an inconsistent set of initial conditions results in instability of pressure and acceleration time histories illustrated by nonphysical oscillations as numerical artifacts.
\subsection{Equation of energy} \label{EquationEnergy}
In this section, we derive the work done by various forces acting on the system during a time step. We show that the work done by pressure (Lagrange multiplier) vanishes, as expected for incompressible media, so the total mechanical energy of the system is preserved. 

The total incremental energy of the system is obtained by multiplying the discretized equations of motion \subeq{eq:spaceDiscritzed}{1,2} at the midpoint $n+\frac{1}{2}$ (i.e., average of the equations at time $n$ and $n+1$) by their respective displacement increments, and then adding them together. The resulting increment of the total energy is
\begin{equation} \label{eq:energy2} 
	\begin{aligned}
		&\delta E^{\text{K}_\text{s}}+\delta E^{\text{K}_\text{f}}+\delta E^{\text{D}}+\delta E^{\text{S}}+\delta E^{\text{C}}=\delta E^{\text{In}}\\
	\end{aligned}
\end{equation}
where $\delta E^{(\cdot)}=E^{(\cdot)}_\text{n+1}-E^{(\cdot)}_\text{n}$. Also $\delta E^{\text{K}_\text{s}}$, $\delta E^{\text{K}_\text{f}}$, $\delta E^{\text{D}}$, $\delta E^{\text{S}}$, $\delta E^{\text{C}}$, and $\delta E^{\text{In}}$ denote the increments of kinetic energy of skeleton, kinetic energy of fluid, damping energy, strain energy of skeleton, work done by constraint force, and input energy respectively. Each of these terms are defined as
\begin{equation} \label{eq:energy3} 
	\begin{aligned}
		&\delta E^{\text{K}_\text{s}} := \frac{1}{2} (\bld{ \dot u}^\top_{\text{n+1}} M  \bld{ \dot u_{\text{n+1}}}-\bld{ \dot u}^\top_{\text{n}} M \bld{ \dot u_{\text{n}}}) \quad \delta E^{\text{K}_\text{f}} := (\bld{ \dot u}^\top_{\text{n+1}} \Mf \,  \bld{ q_{\text{n+1}}}-\bld{ \dot u}^\top_{\text{n}} \Mf \, \bld{ q_{\text{n}}})\,+\, \frac{1}{2} (\bld{ q}^\top_{\text{n+1}} A \,  \bld{ q_{\text{n+1}}}-\bld{ q}^\top_{\text{n}} A \, \bld{ q_{\text{n}}}) \\
		&\delta E^{\text{D}} := \frac{\nf \, g\, \Delta t}{4 \kh} (\bld{q}_\text{n+1}+\bld{ q}_\text{n})^\top\,A\, (\bld{q}_\text{n+1}+\bld{ q}_\text{n}) \quad \delta E^{\text{S}} := (\bld{ u}^\top_{\text{n+1}} K  \bld{ u_{\text{n+1}}}-\bld{ u}^\top_{\text{n}} K \bld{ u_{\text{n}}}) \\
		&\delta E^{\text{In}} := \frac{\Delta t}{4} \left( (\bld{\dot u}_\text{n+1}+\bld{\dot u}_\text{n})^\top  (\mathbb{P}_{n+1}+\mathbb{P}_{n})+(\bld{q}_\text{n+1}+\bld{ q}_\text{n})^\top (\mathbb{F}_{n+1}+\mathbb{F}_{n}) \right)\\
&\delta E^{\text{C}} := -\frac{\Delta t}{2} \left( (\bld{\dot u}^\top_\text{n+1} Q +\bld{ q}^\top_\text{n+1}B)\, +(\bld{\dot u}^\top_\text{n} Q+\bld{q}^\top_\text{n} B) \right)  \bld{ p}_{\text{n}+\frac{1}{2}}
	\end{aligned}
\end{equation}
Where we employ the relationships of \eqref{eq:newmark1} for the displacement increments of solid and fluid. The statement of energy balance \eqref{eq:energy2} says that the increment of the input energy must be equal to the increment of the kinetic energy, strain energy and damping forces energy. We note that the increment of the work done by constraint force is zero, $\delta E^{\text{C}}=0$, since the terms inside the parentheses in \subeq{eq:energy3}{4} are the discretized constraint equations that are set to be zero. Moreover, since we satisfy the constraint at initial time, we get $E^{\text{C}}=0$ as well. 

Using \eqref{eq:energy3}, the energy components at each time step is recovered as
\begin{equation} \label{eq:energy4} 
	\begin{aligned}
		&E^{\text{K}_\text{s}}_{\text{n+1}} = \frac{1}{2} (\bld{ \dot u}^\top_{\text{n+1}} M  \bld{ \dot u_{\text{n+1}}}) \quad E^{\text{K}_\text{f}}_{\text{n+1}} = (\bld{ \dot u}^\top_{\text{n+1}} \Mf \,  \bld{ q_{\text{n+1}}})\,+\, \frac{1}{2} (\bld{ q}^\top_{\text{n+1}} A \,  \bld{ q_{\text{n+1}}}) \\
		&E^{\text{D}}_{\text{n+1}} =E^{\text{D}}_{\text{n}}+\delta E^{\text{D}} \quad E^{\text{S}}_{\text{n+1}} = \bld{ u}^\top_{\text{n+1}} K  \bld{ u_{\text{n+1}}} \quad E^{\text{In}}_{\text{n+1}} =E^{\text{In}}_{\text{n}}+\delta E^{\text{In}}
	\end{aligned}
\end{equation}
The damping and the input energy are time history dependent.
\section{Implementation} \label{Implementation} 
The algorithm for implementation of the proposed numerical method for dynamic analysis of porous media is summarized in Procedure \ref{alg:DynaPoro}.
\begin{algorithm}
    \begin{algorithmic}[1]
	\State Form the coefficient matrices, $M$, $\Mf$, $A$, $K$ , $Q$, $B$ using Procedure \eqref{alg:coef_Matrix}
	\State Select a time increment $\Delta t$ (for an accurate result use Equation \eqref{eq:TimeStepFormula})
	\State Calculate $ \bar{M}$, $ \bar{A}$, $\bar Q$, $\bar B$ using Equation \eqref{eq:TimeDiscritzed3} 
	\State Initialize $\bld{u}_0$, $\bld{\dot u}_0$,  $\bld q_0$ to satisfy \subeq{eq:spaceDiscritzed}{3}
\State Initialize a consistent $\bld{p}_0$ to satisfy \subeq{eq:spaceDiscritzed}{1,2} and \subeq{eq:constInital}{3} \Comment{an example is given in Section \ref{ConsistentInitial}}
	 \For{$\text{n} \gets 0, (\text{numTimInc}-1)$}  
	\State Given $(\tprsc)_{\text{n}}$, $(\tprsc)_{\text{n}+1}$, $(p_\prsc)_{\text{n}}$ and $(p_\prsc)_{\text{n}+1}$, calculate $\mathbb{P}_{\text{n}}$, $\mathbb{P}_{\text{n}+1}$, $\mathbb{F}_{\text{n}}$ and $\mathbb{F}_{\text{n}+1}$ using \eqref{eq:Coef_matrix}.
	\State Update $\bar{\mathbb{P}}_{\text{n+1}}$ and $\bar{\mathbb{F}}_{\text{n+1}}$ using Equation \eqref{eq:TimeDiscritzed3} 
	\State Solve \eqref{Matrix form} for $\bld{\dot u}_{\text{n}+1}$, $\bld q_{\text{n}+1}$, and $\bld p_{\text{n}+1}$
	\State Calculate $\bld u_{\text{n}+1}$ (if required evaluate accelerations as well) using \eqref{eq:newmark1}
	\State For the energy balance analysis:

       Calculate $E^{\text{K}_\text{s}}_{\text{n+1}}$, $E^{\text{K}_\text{f}}_{\text{n+1}}$ and 		$E^{\text{S}}_{\text{n+1}}$ using \eqref{eq:energy4} 

     Calculate $\delta E^{\text{D}}_{\text{n+1}}$ and $\delta E^{\text{In}}_{\text{n+1}}$ using \eqref{eq:energy3}

	Calculate, $E^{\text{D}}_{\text{n+1}}=E^{\text{D}}_{\text{n}}+\delta E^{\text{D}}_{\text{n+1}}$ and 		$E^{\text{In}}_{\text{n+1}}=E^{\text{In}}_{\text{n}}+\delta E^{\text{In}}_{\text{n+1}}$
        \EndFor
    \end{algorithmic}
    \caption{\textsc{Algorithm to solve dynamic poroelasticity}}
    \label{alg:DynaPoro}
\end{algorithm}
Procedure \ref{alg:coef_Matrix} gives an outline of the construction of coefficients defined in \eqref{eq:Coef_matrix} for a two dimensional problem. Section \ref{FEDescription} provides more details on the notations used here. $\Nfn$ and $\Nffc$ denote the number of free nodes and the number of free edges. In building these coefficient matrices, in addition to the connectivity (element to node) and nodal coordinates matrices used in the typical finite element problems, some additional data structures such as element to edge and node to edge are required. An example of the RT element implementation is elaborately discussed by Bahriawatiand and Carstensen \cite{Bahriawati_2005}. To avoid clattering in Procedure \ref{alg:coef_Matrix}, we introduce
\begin{equation}\label{eq:coefMatrix_sum}
\begin{aligned}
	&N_\text{m}=\begin{bmatrix} \tilde N_1(\bld x) &0&\cdots &\tilde N_i(\bld x)&0 \\ 0&\tilde N_1(\bld x)&\cdots&0&\tilde N_i(\bld x) \end{bmatrix}_{i=3 \text{ or } 6 \text{ for linear and quadratic triangle, respectively}}\\
	&W_\text{m}=\frac{1}{2 \mathcal A_m} 	\begin{bmatrix}   s_{m1} \, {l_1} \, (\bld x-\bld x_1) & s_{m2} \, {l_2} \, (\bld x-\bld x_2) & s_{m3} \, {l_3} \, (\bld x-\bld x_3)\end{bmatrix}\\
	&M_\text{m}=\int_{A_m}  {N_m}^\top {N_m} \, dA \quad , \quad A_m= \int_{A_m}  {W_m}^\top {W_m} \, dA
\quad , \quad ({\Mf})_m= \int_{A_m}  {N_m}^\top {W_m} \, dA \\
\end{aligned}
\end{equation}
These integrations over the element domains can be numerically evaluated by Gauss integration method.The following comments and Gauss quadrature rules over triangles \cite{cowper1973gaussian} are considered for numerical examples of section \ref{NumericalExamples}
\begin{itemize}
\item The order of the Gauss integration must be selected to accurately calculate the integrations. For instance, at least a six-point formula (three-point formula), which gives accurate results for up to fourth (second) degree of precision, is required to approximate a mass matrix $M_\text{m}$ with quadratic (linear) interpolation functions. 
\item A reduction in the order of Gauss numerical integration from the order that is required to evaluate these matrices exactly, leads to inaccurate results with long lasting and high amplitude nonphysical oscillations.
\item To decrease the computational effort in case of a quadratic interpolation function, we use a sub-parametric element rather than an iso-parametric one, since all the elements in the numerical examples are constructed as straight edges, and the geometry can be interpreted to a lower degree than the nodal variable's interpolation.
\end{itemize}
Moreover, in forming matrices of Procedure \ref{alg:coef_Matrix}, we utilize the following properties of the $\rt$ mixed element
\begin{equation}    
\begin{aligned}
	& \tilde W_j \cdot \hat e_j=\begin{cases}
      			1   & \text{along edge } j \\
      			0 & \text{otherwise}
    			\end{cases} \quad , \quad \divg \tilde W_j =\begin{cases}
      			 s_{mj} \, \frac{l_j}{\mathcal A_m}    & \text{on } \mathcal {T}_m  \\
      			0 & \text{otherwise}
    			\end{cases}
\end{aligned}
\end{equation}
The right hand side equality is used in construction of matrix $B$. The left hand side equality indicates that the velocity interpolation functions at edge $j$, have the property of producing a unit flux through the $j^{\text{th}}$ edge and 0 flux through all other edges. Furthermore, according to our notation, in a typical triangle, the enumeration on the three vertex nodes is counterclockwise, while the edge enumeration is set such that the edge opposite to vertex node $i$ is the $i^{\text{th}}$ edge of the triangle.  
\begin{algorithm}
    \begin{algorithmic}[1]
        \State Given definitions in \eqref{eq:Coef_matrix}
        \State $A \gets 0_{\Nffc}$ ; $M \text{ and } K \gets 0_{\Nfn}$ ; $\Mf \gets 0_{\Nfn \times \Nffc}$; $Q \gets 0_{\Nfn \times \Nel}$
; $B \gets 0_{\Nffc \times \Nel}$
        \For{$m \gets 1, \Nel$}
		\State Assemble $A_m$, $M_m$ and ${\Mf}_m$ from \eqref{eq:coefMatrix_sum} into $A $, $M$ and $\Mf$ respectively
		\State	Assemble $Q_m=\int_{A_m} \mathcal B^\top \begin{bmatrix} 1&1&0 \end{bmatrix}  dA$  into $Q$
		\State	Assemble  $B_m=\frac{1}{A_m} \int_{A_m} \begin{bmatrix}s_{m1}\, l_1 & s_{m2}\, l_2 &s_{m3}\, l_3 \end{bmatrix}^\top dA  $ into $B$
		\State	Assemble  $K_m=\int_{A_m} { \mathcal B_m}^\top  \mathcal C { \mathcal B_m}  dA$ into $K$
        \EndFor
    \end{algorithmic}
    \caption{\textsc{Formation of coefficient matrices for two dimensional problems}}
    \label{alg:coef_Matrix}
\end{algorithm}
\nomenclature{$\Nffc$}{number of free edges}
\nomenclature{$\Nfn$}{number of free nodes}
\section{Numerical Examples} \label{NumericalExamples}
To demonstrate the capabilities of our formulation to simulate the incompressible dynamical behavior of saturated porous media, we perform two benchmark numerical examples and verify them with analytical solutions (Example 1) and boundary element solutions (Example 2).
As was noted in section \ref{LimitingCases}, it has been found that when the hydraulic conductivity $\kh$ is very small, locking (in terms of failing the LBB condition) occurs \cite{Liu_2004,Phillips_Wheeler_2009}. To explore this effect for our formulation, we consider a third set of examples to evaluate the stability of the method when locking is often a concern. Overall, we found that the presented numerical scheme is stable and sufficiently accurate for a wide range of material parameters including large and small values of hydraulic conductivity. 
\subsection{Example 1: half space -- soil column} \label{Example1_halfSpace}
\begin{figure}
    \centering  
        \resizebox{0.7\textwidth}{!}{\includegraphics {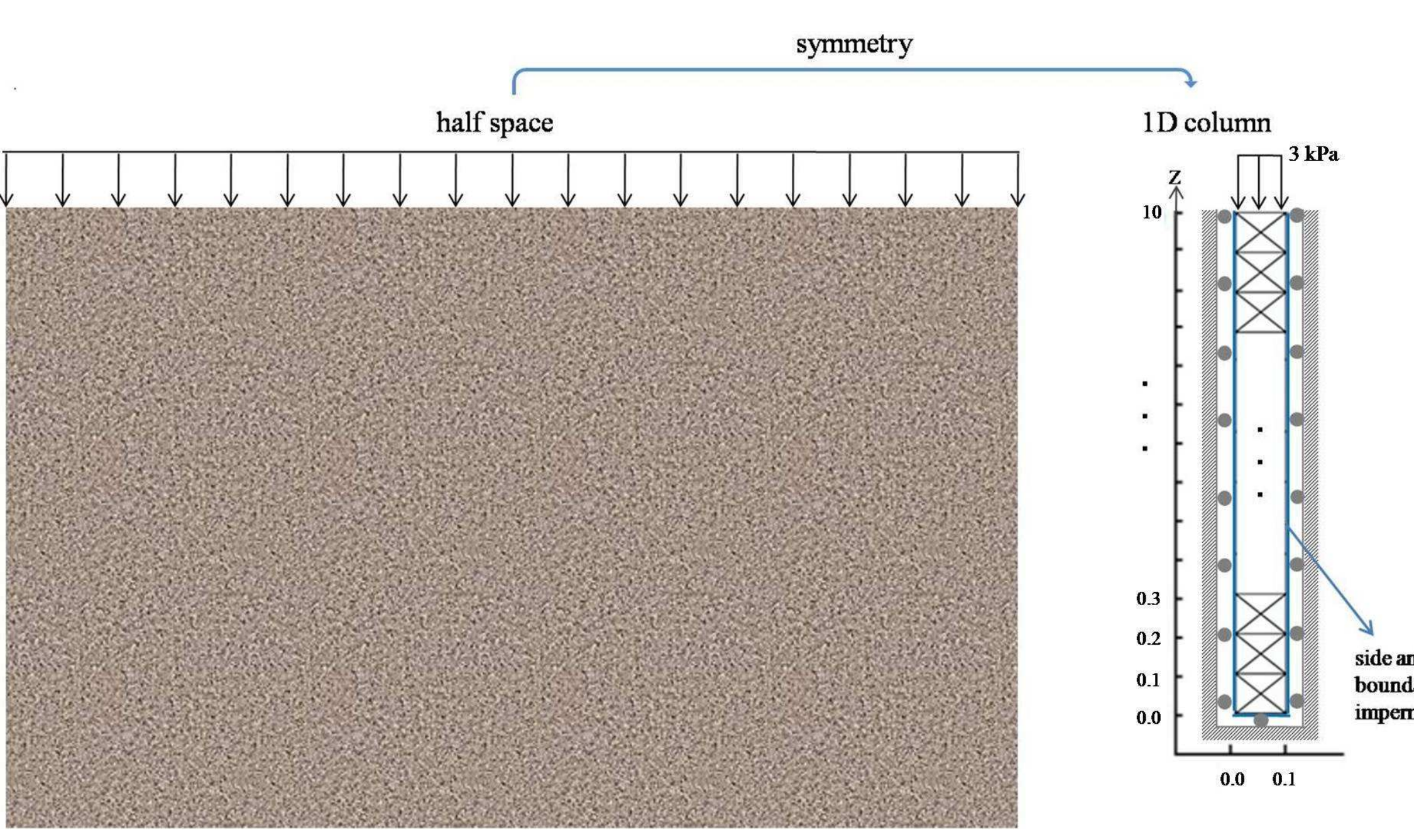}}
    \caption[Example 1 --- model]{Example 1 ---- a column of soil under step load: geometry and finite element mesh. We have isolated a 10 m length by 0.1 m width for analysis. The figure is not in scale. Three dots show the continuation.}
    \label{fig:example1_1}
\end{figure}

The first benchmark example is the application of a uniform load on the surface of a half space\cite{Prevost_1985,Chen_Dargush_1995,Li_etal_2004,Diebels_Ehlers_1996,Zienkiewicz_Shiomi_1984,Arduino_Macari_II_2001}. Due to symmetry, only a representative column of small width, 0.1 m,
is studied and the problem is of one dimension. The model is taken as an incompressible saturated soil with plane strain behavior. The material data are shown in Table \ref{tab:Material_properties} and are taken from the same benchmark example as in \cite{Diebels_Ehlers_1996,deBoer_etal_1993}. The geometry and finite elements are illustrated in Figure \ref{fig:example1_1}, where the side walls and the bottom are impermeable and the displacements normal to their surfaces are constrained. The upper boundary is perfectly drained and subjected to a step load of 3 $\frac{\text{kN}}{\text{m}^2}$ amplitude, which results in a total load of $0.3 \text{kN}$ on top boundary. 400 finite elements with a crisscross mesh pattern are used to model the soil column. We also use a time step of $\Delta t=1 \times 10^{-4} $ (see section \ref{sec:time step}) that is adequately small for the purpose of wave propagation analysis. We have provided the results utilizing the $\linprt$ element and a lumped mass matrix unless otherwise mentioned.
\begin{table}
            \centering
            \begin{tabular}{ccccccc}
            	\toprule   
            	 & $E$ $\frac{\text{kN}}{\text{m}^2}$ & $\nu$ & $\rhos$ $\frac{\text{kg}}{\text{m}^3}$ & $\rhof$ $\frac{\text{kg}}{\text{m}^3}$ & $ n $ & $\kh$  $\frac{\text{m}}{\text{s}}$ \\
             \midrule 
       		 Example 1& 14.516e+3   & 0.3 & 2000  &1000 
		&0.33& 1e-2  \\
		Example 2& 14.516e+3   & 0.3 & 2700  &1000 
		&0.42& 1e-1 and 1e-4 \\
		Example 3& 10   & 0.4 & 2667  &1000 
		&0.4& 1e-7 \\
            	\bottomrule
            \end{tabular}
            \caption{Material properties for numercial examples}
            \label{tab:Material_properties}
\end{table}

\noindent With this example, we pursue many goals as listed below:
\begin{enumerate}
\item Verification with analytical solutions
	\begin{itemize}
	\item Short-term behavior: we compare the numerical results with the analytical solution of de Boer et al. \cite{deBoer_etal_1993} for wave propagation in an infinitely long column of incompressible saturated soil. The influence of the total depth of the column is negligible \cite{Diebels_Ehlers_1996}; therefore, for our numerical analysis, we are able to isolate a 10-m length and still compare the short term behavior of our system with de Boer's results for an infinite long column.
	\item Long-term behavior: we compare long term behavior with Terzaghi's equation for consolidation. 
	\end{itemize}
\item Wave propagation analysis: according to the well known derivation of Biot \cite{Biot_lowfrequency_1956}, there are three types of waves that pass through a porous medium. Two dilatational waves (fast and slow) and one shear wave \footnote{Accounting for apparent mass, there are two shear waves, one in fluid and one in solid skeleton\cite{Biot_lowfrequency_1956}.}. It is proved analytically that as a consequence of the incompressibility condition, the fast dilatational wave propagates with an infinite speed and practically, only the slow dilatational wave exists \cite{deBoer_etal_1993,Schanz_Pryl_2004}. Also, note that for this one-dimensional problem, there is no shear waves.
\item Pressure build-up and diffusion: we are interested in studying the initial and long term behavior of the pressure field. 
\item Effect of instantaneous loading: the impact and shock nature of the step load creates a challenging condition for the incompressible porous medium. The effect of this condition on the response history of the medium is studied, and some remedies are suggested to avoid the associated numerical artifacts.
\item Effect of mass representation: it is suggested in the literature, that such numerical artifacts can be alleviated by appropriate representation of mass matrix. We explore different mass lumping methods and compare their results with those of the consistent mass matrix method.
\item Energy balance analysis: in order to highlight the working of the presented numerical algorithm, we monitor the balance of the energy in the porous medium. In particular, by means of the energy balance equation, we can evaluate our time integration scheme and size of the time step.
\end{enumerate}
\subsubsection{Verification with analytical solutions}
\begin{figure}
    \centering
 \resizebox{0.5\textwidth}{!}{\includegraphics{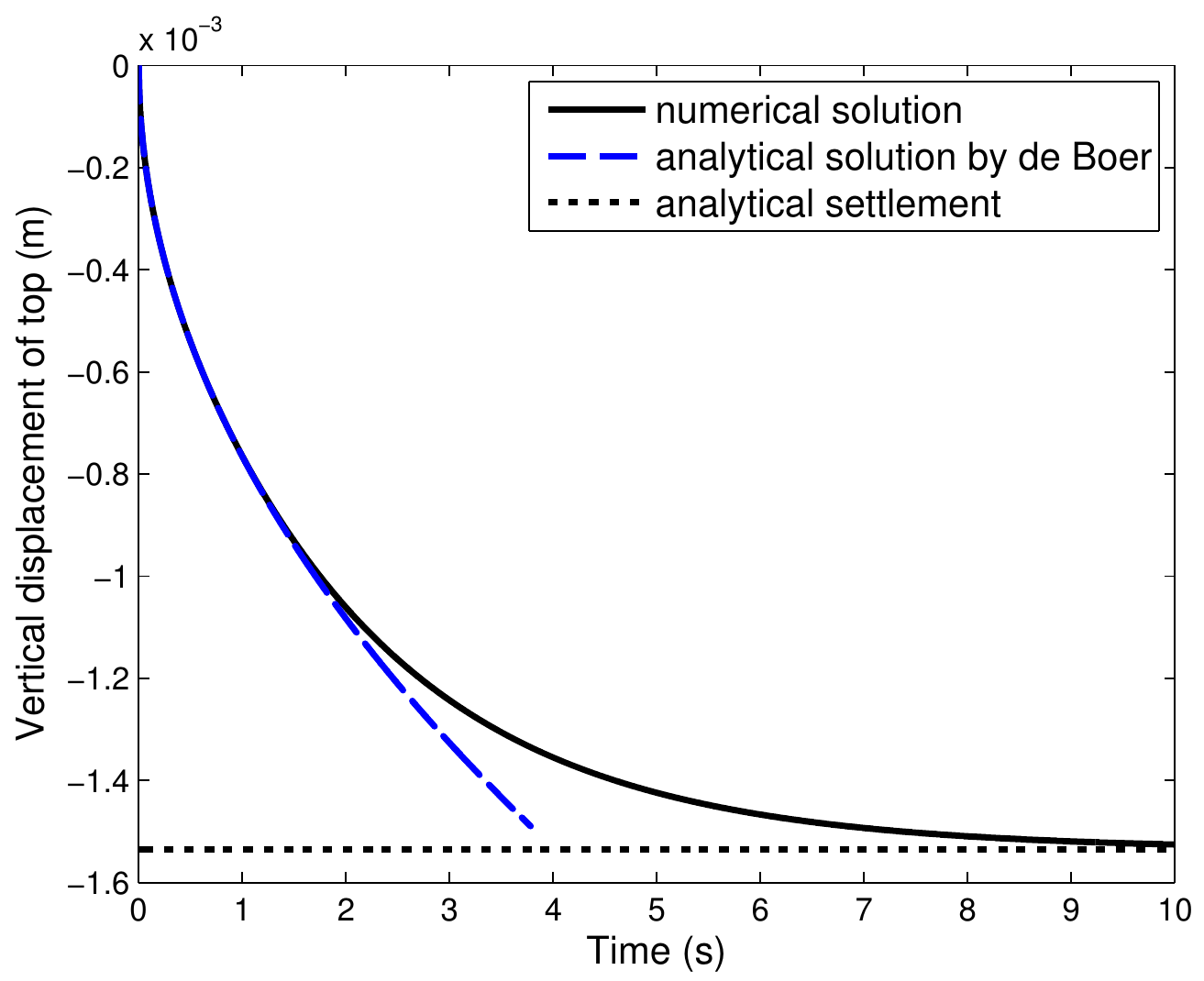}}     
        \caption[Example 1 --- displacement time history]{Example 1 --- displacement time history. Good agreement between numerical and the analytical computations is observed for both the transient behavior given by de Boer et al. \cite{deBoer_etal_1993} and the long term behavior given by Terzaghi \ref{terzaghi}}
    \label{fig:example1_disp}
\end{figure}

\mbox{}\\
Shown in Figure \ref{fig:example1_disp} is the displacement response at the top boundary. Our numerical results match well with de Boer's analytical solution in the transient region \cite{deBoer_etal_1993}. We point out that de Boer's solution is for an infinitely long column where there is no reflection of wave from the bottom boundary; however, we isolate a finite length for computations, and therefore our numerical plots separates from de Boer's solution after a while. Moreover, to verify the long term behavior, we use Terzaghi's equation for steady state consolidation. According to that equation, the total settlement at the top boundary \cite[equation 2-3]{Pant_2007} can be simplified for such one dimensional problems as:
\begin{equation} \label{terzaghi}
\delta s = m_\text{v}\, f \, L
\end{equation}
where $f$ is the total load, $H$ is the length of the soil column, and $m_\text{v}$ is the the slope from a one dimensional test,
\begin{equation}
m_\text{v}= \frac{d\varepsilon}{d\sigma}=\lambda+2G
\end{equation}
where a plain strain assumption is made. 
\subsubsection{Wave propagation analysis}\label{sec:time step}
\begin{figure}
    \centering
 \resizebox{0.5\textwidth}{!}{\includegraphics{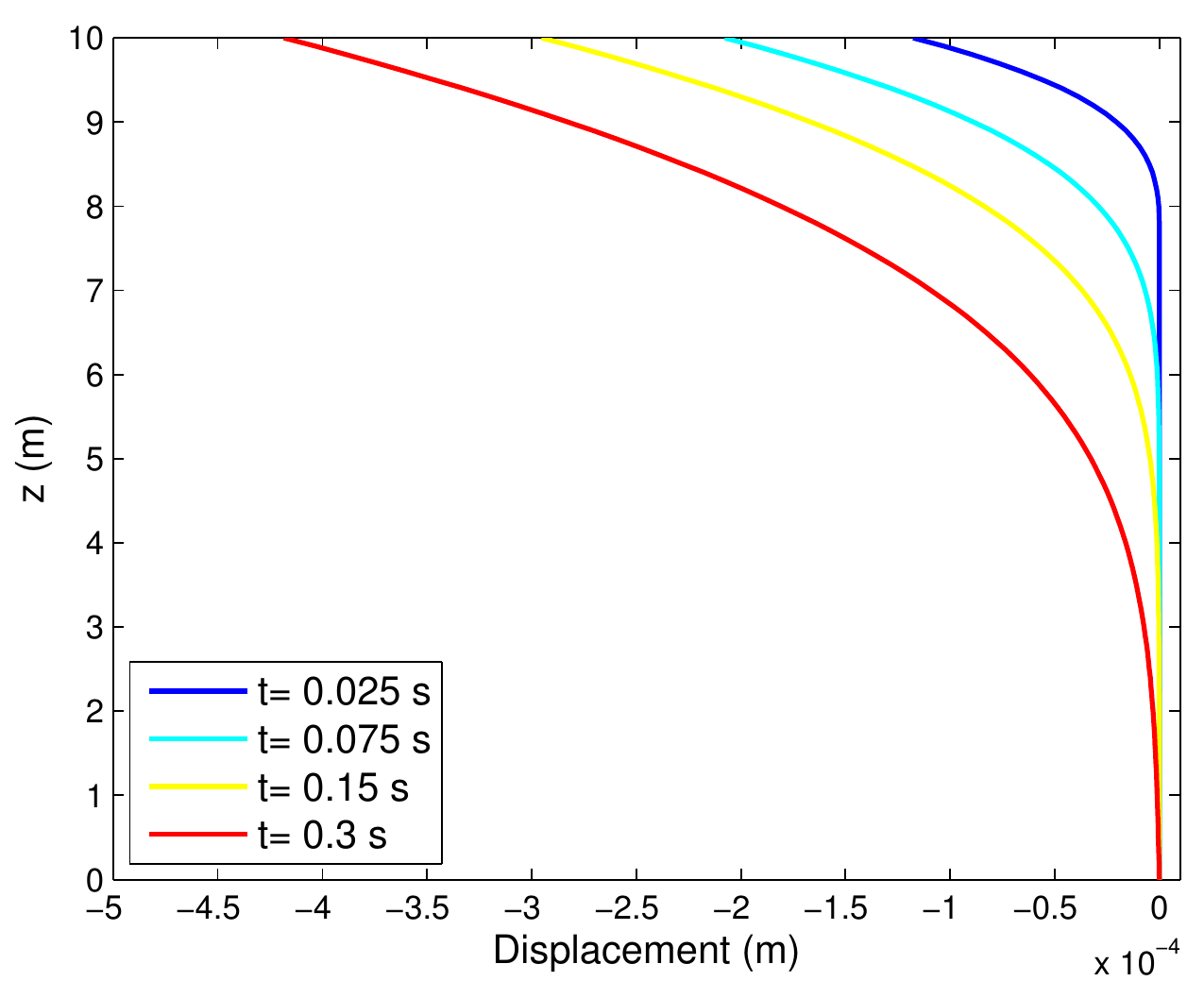}}  
        \caption[Example 1 --- snapshots of the wave propagation]{Example 1--- snapshots of wave propagation through the medium. The propagation speed reduces over the time as roughly calculate in \eqref{eq:waveSpeed_2}.}
    \label{fig:wave}
\end{figure}

\mbox{}\\ 
\textit{wave speed and choice of time step} \\ 
The significance of the choice of an adequate time step is clear in a wave propagation problem. In order to represent the wave travel accurately, it is recommended that the time step is determined by the approximate Courant-Friedrichs-Lewy (CFL) condition (equation 9.102,\cite{bathe2006finite})
\begin{equation}\label{eq:TimeStepFormula}
\Delta t \leq \frac{\Delta x_{\text{min}}} {C_0}
\end{equation}
This is the smallest time that is required for a dilatational wave to pass across any of the elements.
Given the smallest element size, $\Delta x_{\text{min}}$, we can compute $\Delta t$ knowing the wave speed $C_0$. An estimate of the slow dilatational wave speed is given by de Boer et al.\cite{deBoer_etal_1993} as,
\begin{equation} \label{eq:waveSpeed}
	C_0=\sqrt{\frac{m_v}{\rho-\rhof(2-\frac{1}{n})}}
\end{equation}
Using the parameters of Example 1 in the above relationship, we get $C_0=85.1 \ms$. In the next section, we compute wave velocity using the plots of Figure \ref{fig:wave} and demonstrate correlation with this estimate.

Setting the minimum mesh size, $\Delta x= \frac{0.1}{\sqrt{2}}$ m, to be the short-edge length for triangles, equation \eqref{eq:TimeStepFormula} results in a time step of $\Delta t \leq 8.3 \times 10^{-4}$, so we apply a time step of $\Delta t=1 \times 10^{-4}$.  

\noindent \textit{Wave velocity analysis}
\\ We are interested in observing the passage of the slow dilatational wave which is clear in Figure \ref{fig:wave} recorded at some snapshots of time. The time span was selected such that no reflections of the wave affect the results. In this plot, we see the propagation of the wave induced by the load at $z=10$ m and its propagation through the medium more and more as time progresses. Using this figure, the approximate wave velocity can be calculated as: $C_0=\frac{\delta z}{\delta t}$. For instance, taking the first three plots we get:
\begin{equation} \label{eq:waveSpeed_2}
\begin{aligned}
	& C_0=\frac{10-7.8}{0.0025-0}=88 \ms\\
	& C_0=\frac{10-4.6}{0.075-0}=72 \ms\\
	& C_0=\frac{10-1.5}{0.150-0}=56 \ms
\end{aligned}
\end{equation}
Recall that an approximate value of 85 $\ms$ was calculated earlier using de Boer's equation \eqref{eq:waveSpeed}; that is close to the early time value of our numerical speed in Equation \eqref{eq:waveSpeed_2}. However, our numerical results show a decreasing trend of wave speed as the wave travels through the media which is predicted by Schanz and Pryl \cite{Schanz_Pryl_2004} as well. One possible reason for this difference is the finite length of the column in our example compared to infinite length of de Boer analytical formulation.

Time histories of skeleton velocity at different locations $(z)$, shown in Figure \ref{fig:example1_velZoom}, also reveals the passage of the wave through the medium.  We observe that moment when the wave reaches to a certain depth, the velocity profile of that depth separates from the rest of the profiles with zero velocity.
\subsubsection{Pressure build-up and diffusion}
\begin{figure}
    \centering
    \subfloat[]
    {
        \resizebox{0.4\textwidth}{!}{\includegraphics{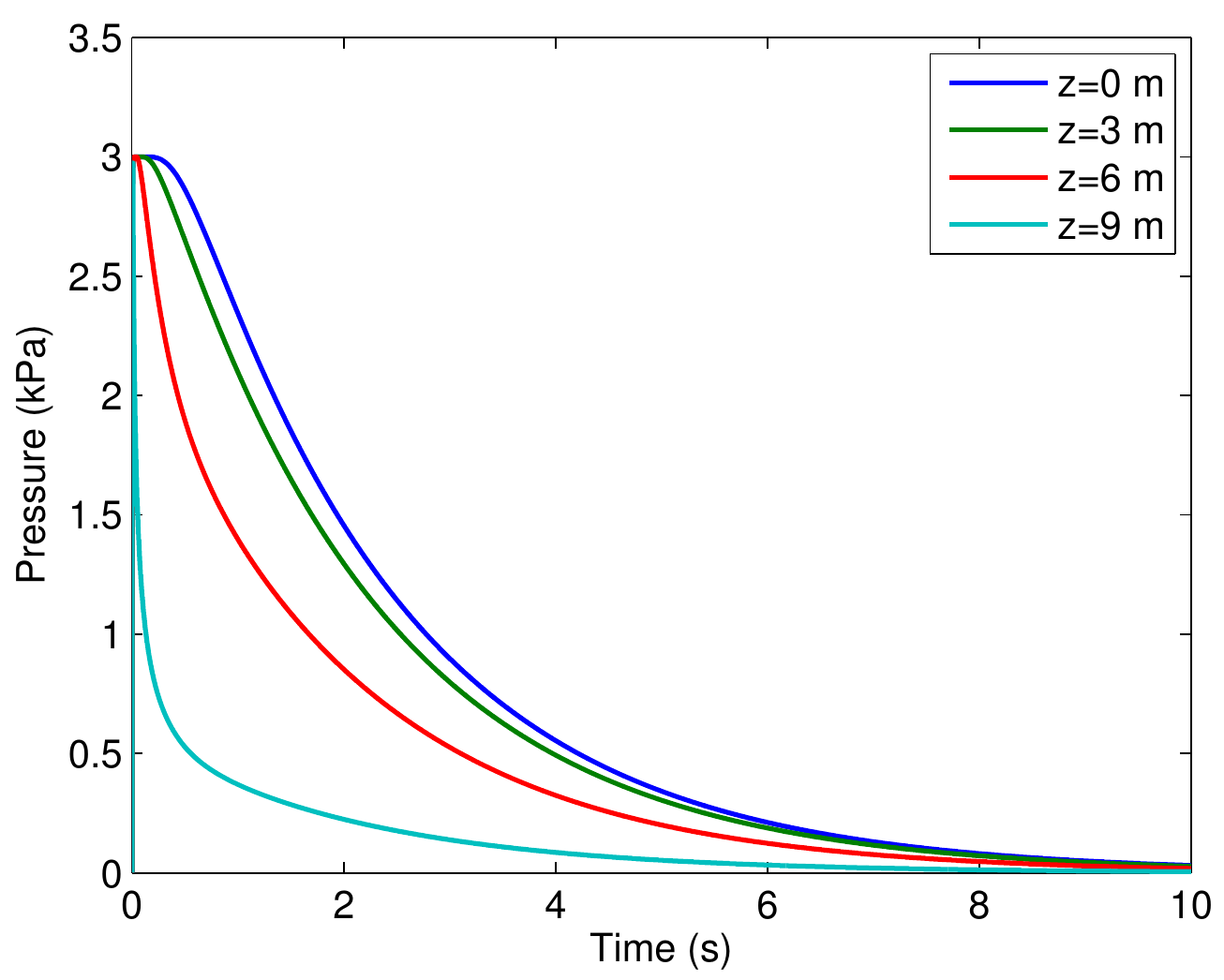}}
        \label{fig:example1_prs}
    }
    \hspace{0.05\textwidth}
    \subfloat[]
    {
        \resizebox{0.4\textwidth}{!}{\includegraphics{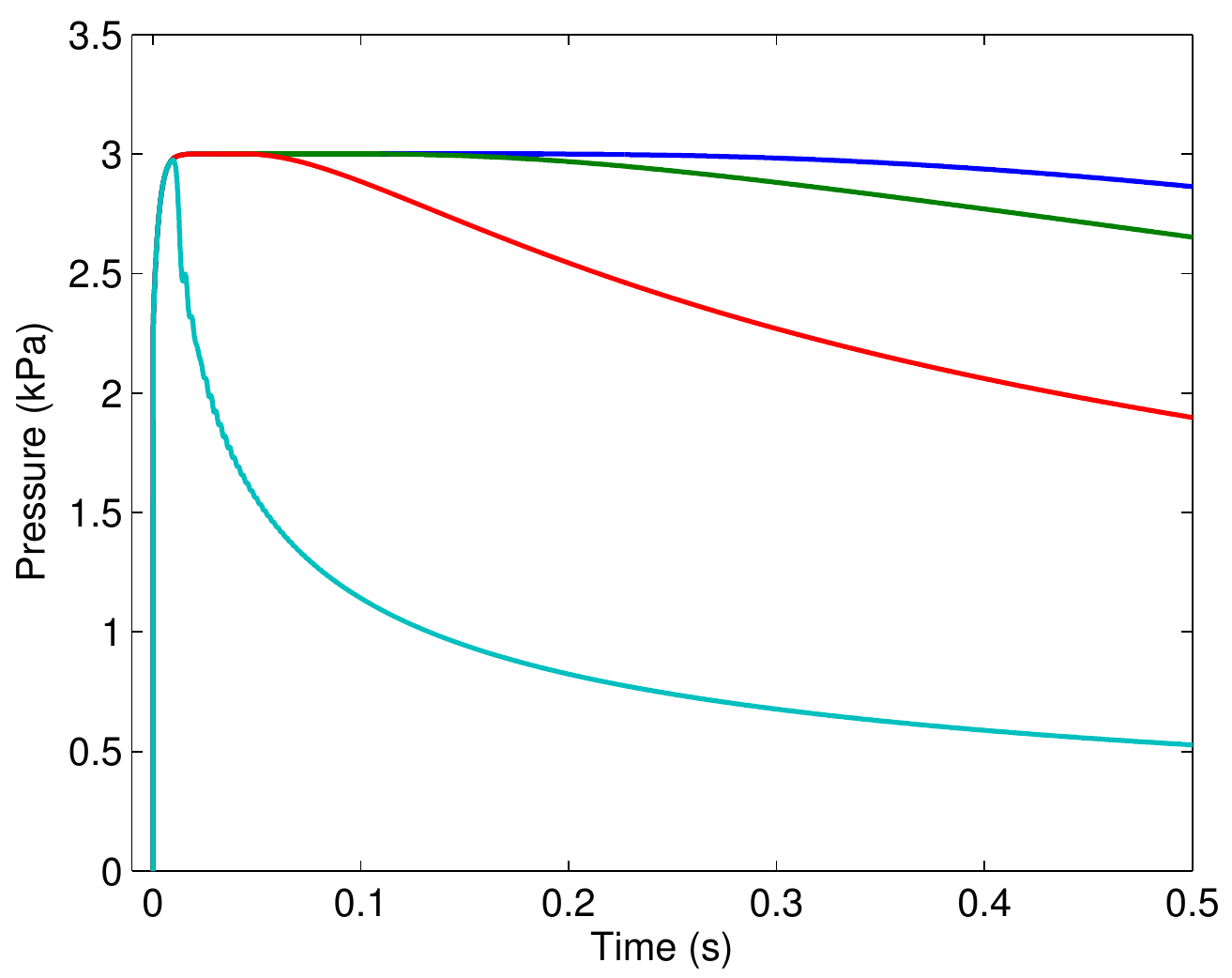}}
        \label{fig:example1_prs_zoom}
    }
    \caption[Example 1--- pressure time history]{Example 1--- (a) Pressure time histories. Diffusion of pressure is apparent over the long time. (b)  Zoom in of first moments. Pressure jumps upon loading at time zero. Also note the diffusion begins right after the passage of the stress wave.}
    \label{fig:example1_prs}
\end{figure}
\mbox{}\\Pore fluid pressure time histories computed at different depths is displayed in Figure \ref{fig:example1_prs} for short and long term. Upon the instantaneous application of the load, a pressure build-up occurs at all points of the medium, indicating the propagation of the pressure at infinite speed \cite{gresho1987pressure}(Figure \ref{fig:example1_prs_zoom}). After the stress wave hits a depth, it induces a change in pressure which results in diffusion of the pressure back toward its initial state.
\subsubsection{Effect of instantaneous loading}
\mbox{}\\
\begin{figure}
    \centering
    \subfloat[]
    {
        \resizebox{0.45\textwidth}{!}{\includegraphics{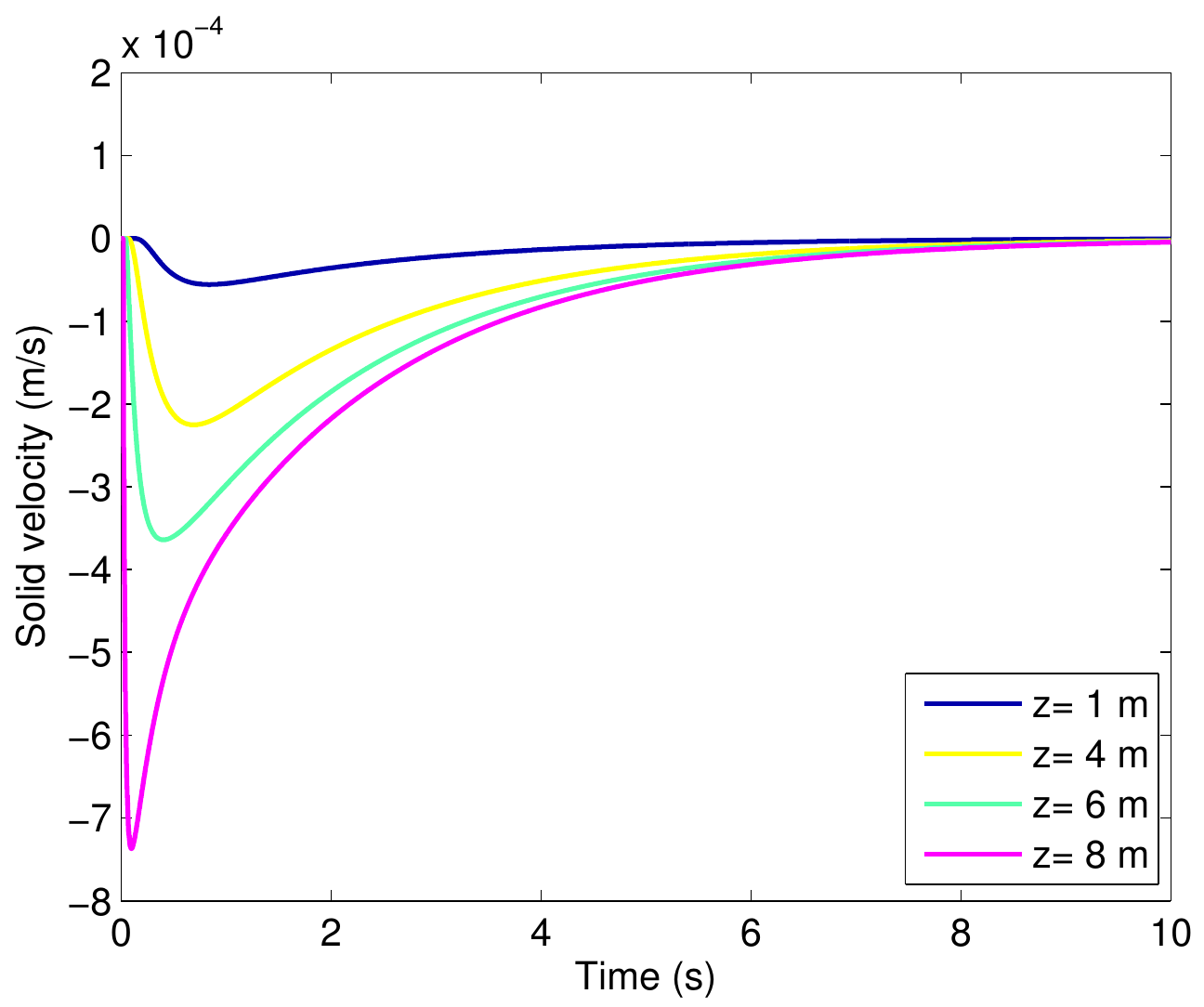}}
        \label{fig:example1_vel}
    }
    \subfloat[]
    {
        \resizebox{0.45\textwidth}{!}{\includegraphics{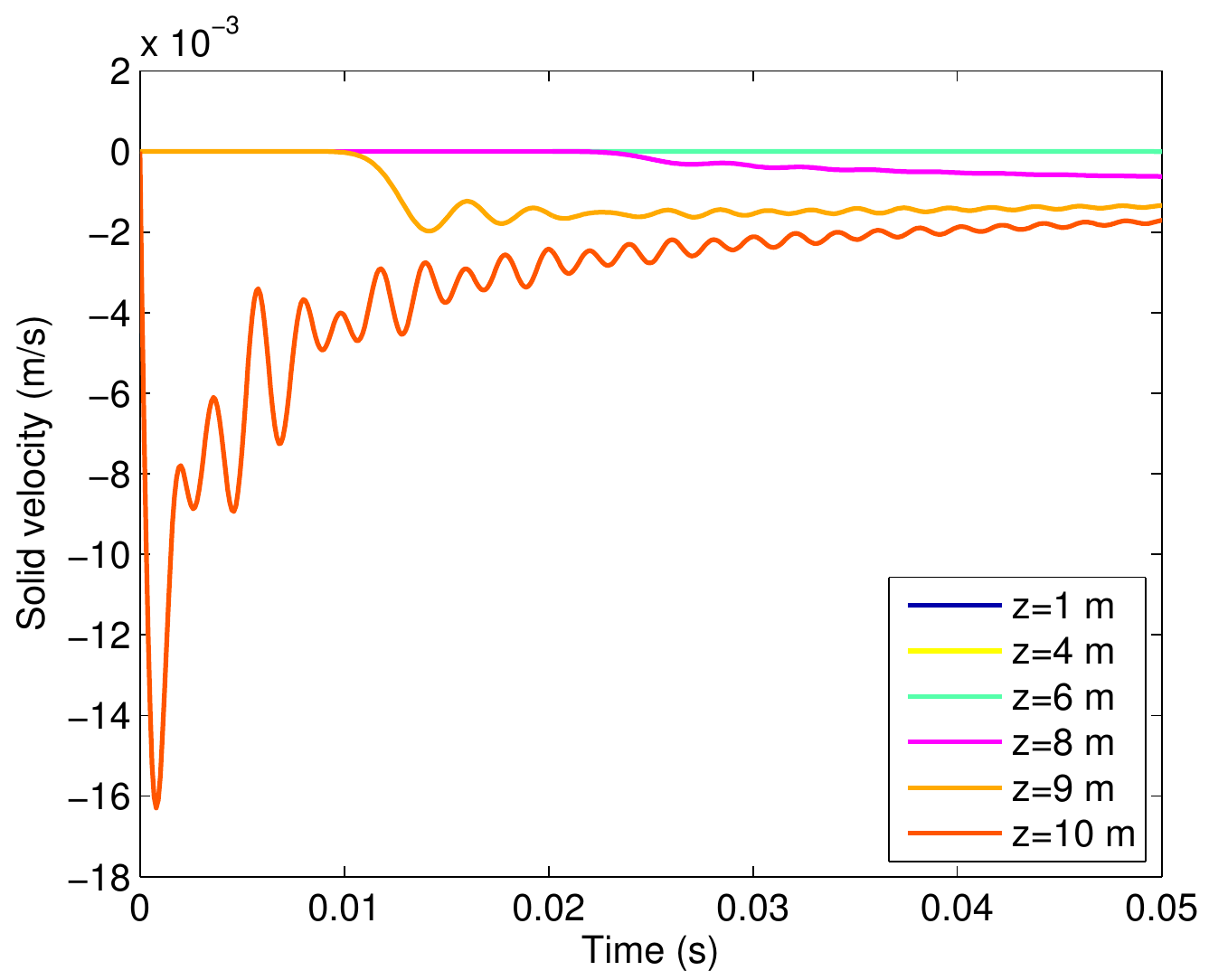}}
        \label{fig:example1_velZoom}
    }
    \caption[Example 1--- skeleton velocity time history]{Example 1--- (a) Solid skeleton velocity time histories for some representative depths. The trend, however, is the same at all locations with their velocities ultimately reach to zero. (b) Zoom in of first moments at some critical depths. The high frequency oscillations, byproducts of step loading, are more pronounced close to the top boundary, but they subside as time marches forward.}
    \label{fig:example1_velocity}
\end{figure}
The time history of the skeleton velocity is shown in Figure \ref{fig:example1_velocity}. Over time, as energy dissipates and diffusion occurs, the velocity tends to zero (Figure \ref{fig:example1_vel}). The trend in this plot is the same at all depths, so we only display the response at a few representative locations. 

\begin{figure}
    \centering
 \resizebox{0.5\textwidth}{!}{\includegraphics{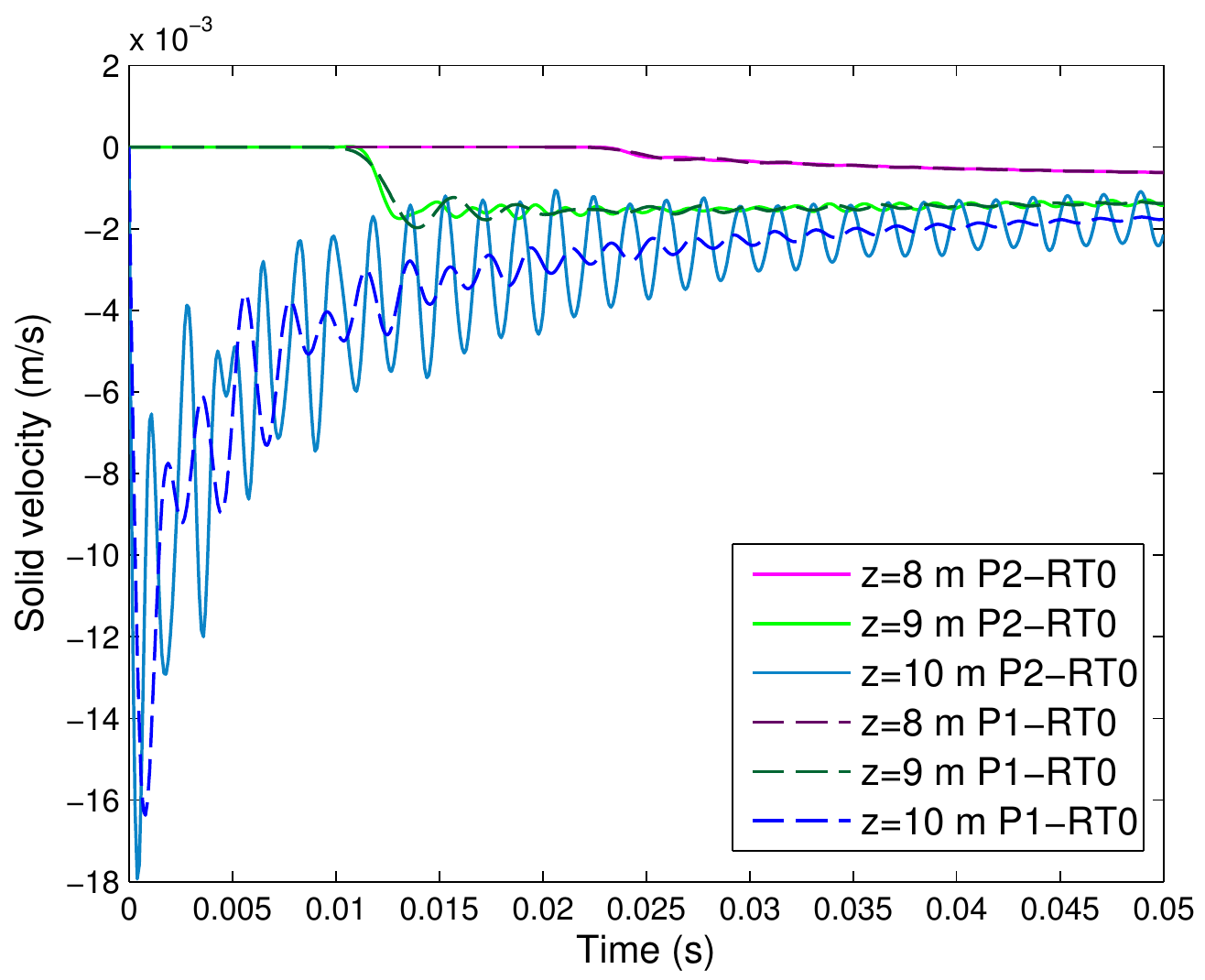}}  
        \caption[Example 1 --- a comparison between linear and quadratic elements]{Example 1 --- a comparison of the solid skeleton velocity time histories. In the proximity of the loading boundary, the $\quadprt$ element exhibits higher amplitudes and frequency oscillations compared to the linear $\linprt$ element. After a while, both elements converge to identical results.}
    \label{fig:LST_CST_1d}
\end{figure}

The step loading due to its abruptness creates a jump in pressure and inertia forces. This in turn excites some higher modes of vibration that survive only for a short time interval due to presence of physical damping in the form of fluid diffusion. These early high frequency oscillations are pronounced in the vicinity of the top boundary with rapid pressure change as shown in Figure \ref{fig:example1_velZoom}. We recognize that these high frequency oscillations are a numerical artifact as their frequencies become unbounded upon mesh refinement. We note that decreasing the time step does not remove these spurious oscillations either. We also found that these nonphysical oscillations are not affected significantly (or removed) by finite element type, as presented in Figure \ref{fig:LST_CST_1d}. In these plots, we compare velocity response of the medium using $\linprt$ and $\quadprt$ elements with the same number of elements. It is observed that the response using $\quadprt$ element has higher amplitude and frequency oscillations in comparison with $\linprt$ element solution. However, as we move away form the loading boundary, the responses of both finite elements are indistinguishable. 

To avoid these numerical oscillations, we explore different remedies -- (1)  Replacing the step load with a ramp load (even with a relatively short rise time) smooths out the initial oscillations (the results of this study are not shown here); (2) An alternative approach is to apply time integration methods containing algorithmic dissipation such as the Generalized $\alpha$ method \cite{chung1993time}. We note that the results shown here are computed by the constant average Newmark scheme where no algorithmic dissipation is used and consequently, the total energy of the system \eqref{eq:energy2} is preserved. This fact is illustrated in Figure \ref{fig:example1_engComp} as well; (3) While such modifications of time integration algorithm can be easily implemented, we instead search for remedies at the element formulation level. In \cite{Gosz_Michael_FEBook}, it is suggested that such oscillations may be associated with representation of mass matrix. This is discussed in the next section \ref{sec:massRepresentation}. We also point out that these early numerical oscillations do not appear in displacement time histories since the integration process automatically filters out the high frequency components.
\subsubsection{Effects of mass representation} \label{sec:massRepresentation}
\begin{figure}
    \centering
 \subfloat[]
    {
        \resizebox{0.45\textwidth}{!}{\includegraphics{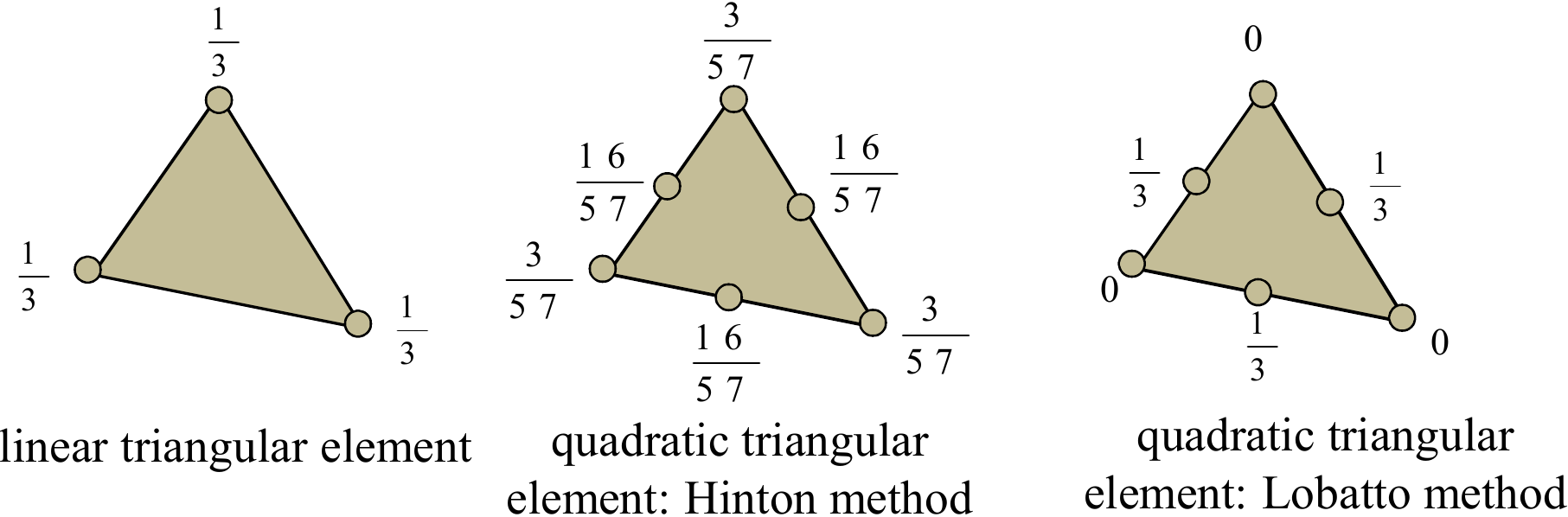}} 
        \label{fig:lumping_weights}
    }
    \vspace{-0.5\baselineskip}
    \subfloat[]
    {
        \resizebox{0.45\textwidth}{!}{\includegraphics{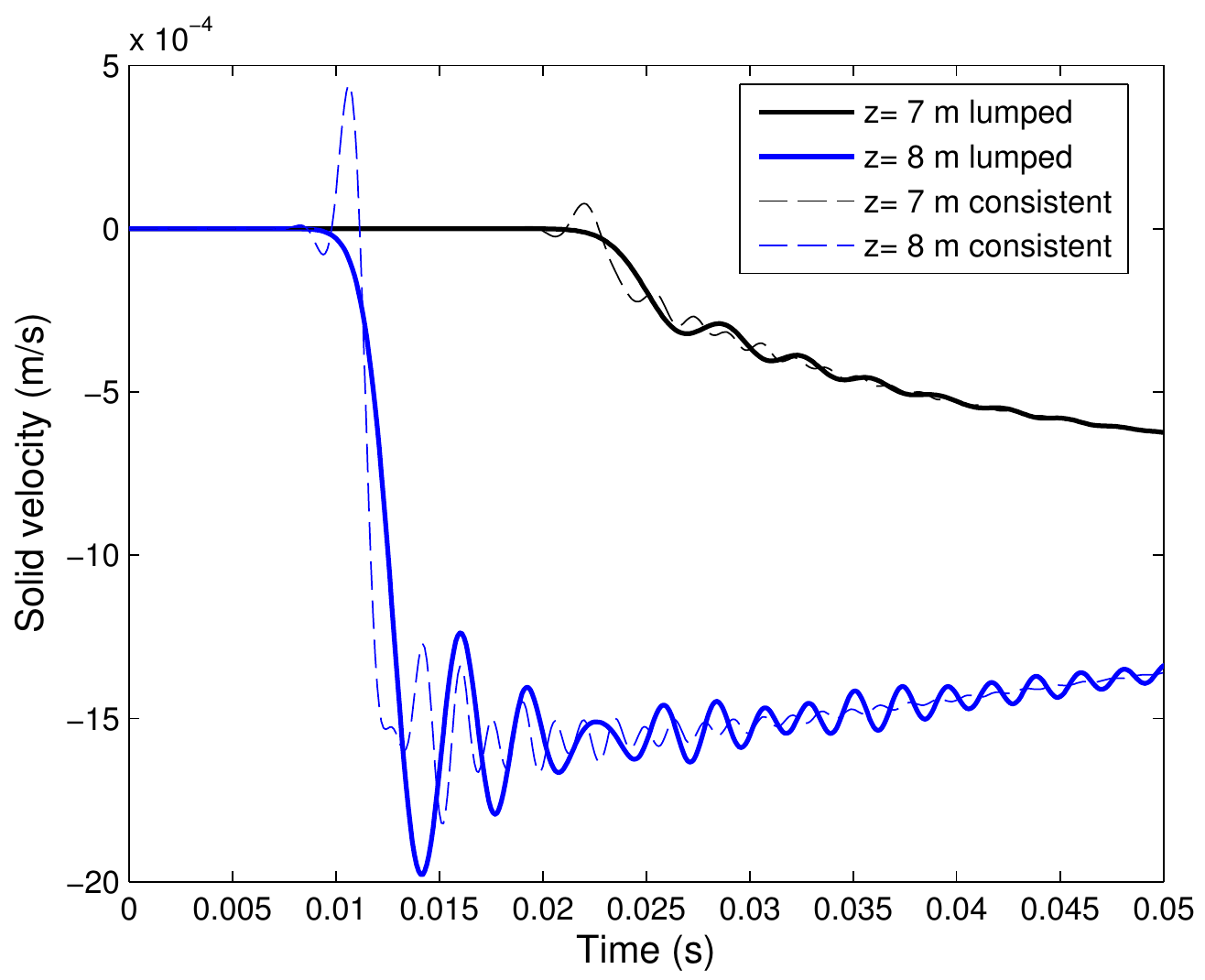}}
        \label{fig:example1_lumpingAnal}
    }  
    \subfloat[]
    {
        \resizebox{0.46\textwidth}{!}{\includegraphics{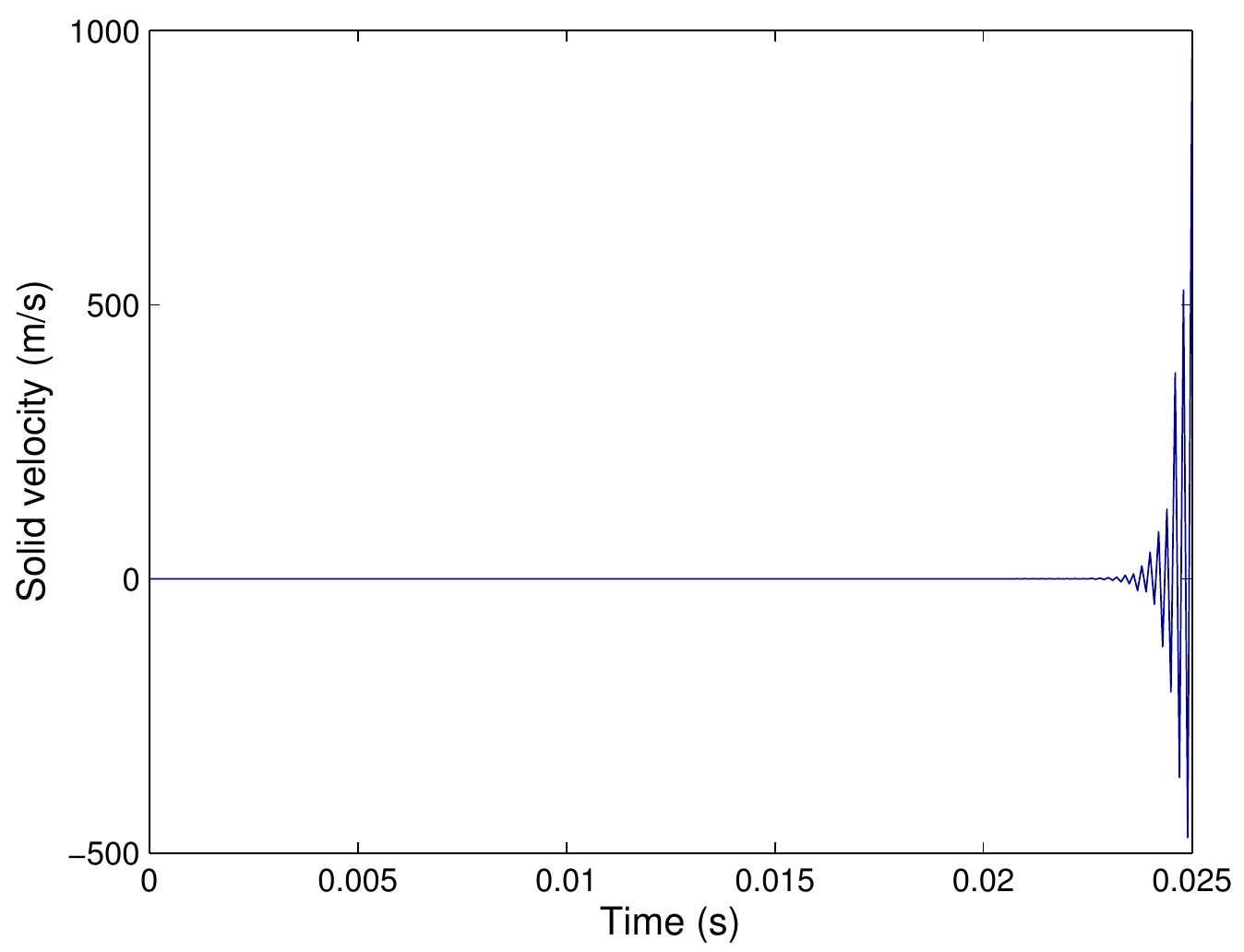}}
        \label{fig:vel_unstable}
    }  
    \caption[Example 1 --- mass lumping]{Example 1 --- (a) Weighting factor of each nodal DOF in the element lumped mass matrix. Note that over the linear triangle element, $\linprt$, both lumping procedures result in identical nodal weights. (b) Comparison of lumped and consistent mass matrix responses for the $\linprt$ element. The consistent mass matrix invokes a greater amount of high frequency oscillations, especially upon the arrival of the wave front. (c) Unstable response of the quadratic triangle element, $\quadprt$, when a Lobatto quadrature lumping method is used.}
    \label{fig:example1_lump_cstVSlst}
\end{figure}

\mbox{}\\In constructing the matrix $M$ defined in \eqref{eq:Coef_matrix}, we use the same shape functions for velocities and displacements, resulting in a consistent variational formulation as presented in \eqref{eq:weakform}. Accordingly, $M$ represents a consistent mass matrix, which is generally non-diagonal.
An alternative to the consistent mass matrix is to employ a diagonal mass matrix based on direct lumping. In this method, the total mass of element is directly apportioned to nodal DOFs, ignoring any cross coupling. There are several ways of constructing a diagonally lumped mass matrix. Here, we test two of them, namely, mass lumping by nodal Lobatto quadrature \cite{hughes2012finite} and the lumping method proposed by Hinton \cite{hinton1976note}.
We note that in this example the lumping methods are applied on the mass matrices associated with the block related to the skeleton displacement, i.e, $M$ and $\Mf$ \footnote{The lumping process can be viewed as substituting a different interpolation function for velocity $N_L$ (as opposed to using the same interpolation function as used for displacement) to approximate the mass matrix as, $M_L=\int_V {N_L}^\top N_L  \, d\Omega$. Accordingly, to build the coupled mass matrix $\Mf$ we must use the same interpolation function $N_L$ in order to tune it with the lumped mass $M_L$.}, while we keep using the same mass matrix for the fluid velocity block, i.e., $A$; since there is no notion of fluid displacement in our formulation and instead the primary variable filed is fluid velocity. The resulting weighting factor of each DOF in the element mass matrix is illustrated in Figure \ref{fig:lumping_weights}. Note that over the linear triangle which is used in the $\linprt$ element, both procedures give rise to the identical nodal weights.\\
Shown in Figure \ref{fig:example1_lumpingAnal}, we provide the numerical results employing both the lumped and consistent mass matrices for the linear triangle element, $\linprt$. The lumped mass matrix results are displayed in solid line while the consistent mass matrix results are in dashed line. We see that the consistent mass matrix creates oscillations with higher frequencies compared to the lumped mass matrix response. In other words, it is more noisy compared to lumped mass response especially upon the arrival of the wave front. Nevertheless, in both cases, noise is damped out after a short while.

With the quadratic triangle element, $\quadprt$, we see the same trend as for linear triangle element when the Hinton et al. method of lumping is applied. However, the use of the Lobatto quadrature lumping method gives rise to an unbounded and unstable solution as plotted in Figure \ref{fig:vel_unstable}. This numerical instability occurs regardless of the mesh size, the time step size, and the hydraulic conductivity. In fact, we observe that the Lobatto quadrature lumping method over the quadratic triangular element produces unstable numerical solution even in a typical dynamic elasticity problem. We also emphasize that when using a quadratic triangle element, a consistent load vector as defined in \eqref{eq:Coef_matrix} should be exerted regardless of mass matrix begin lumped or consistent. Otherwise, very inaccurate results are obtained as reported in \cite{bathe2006finite} as well.
\subsubsection{Energy plots}
\begin{figure}
    \centering
    \subfloat[]
    {
        \resizebox{0.42\textwidth}{!}{\includegraphics{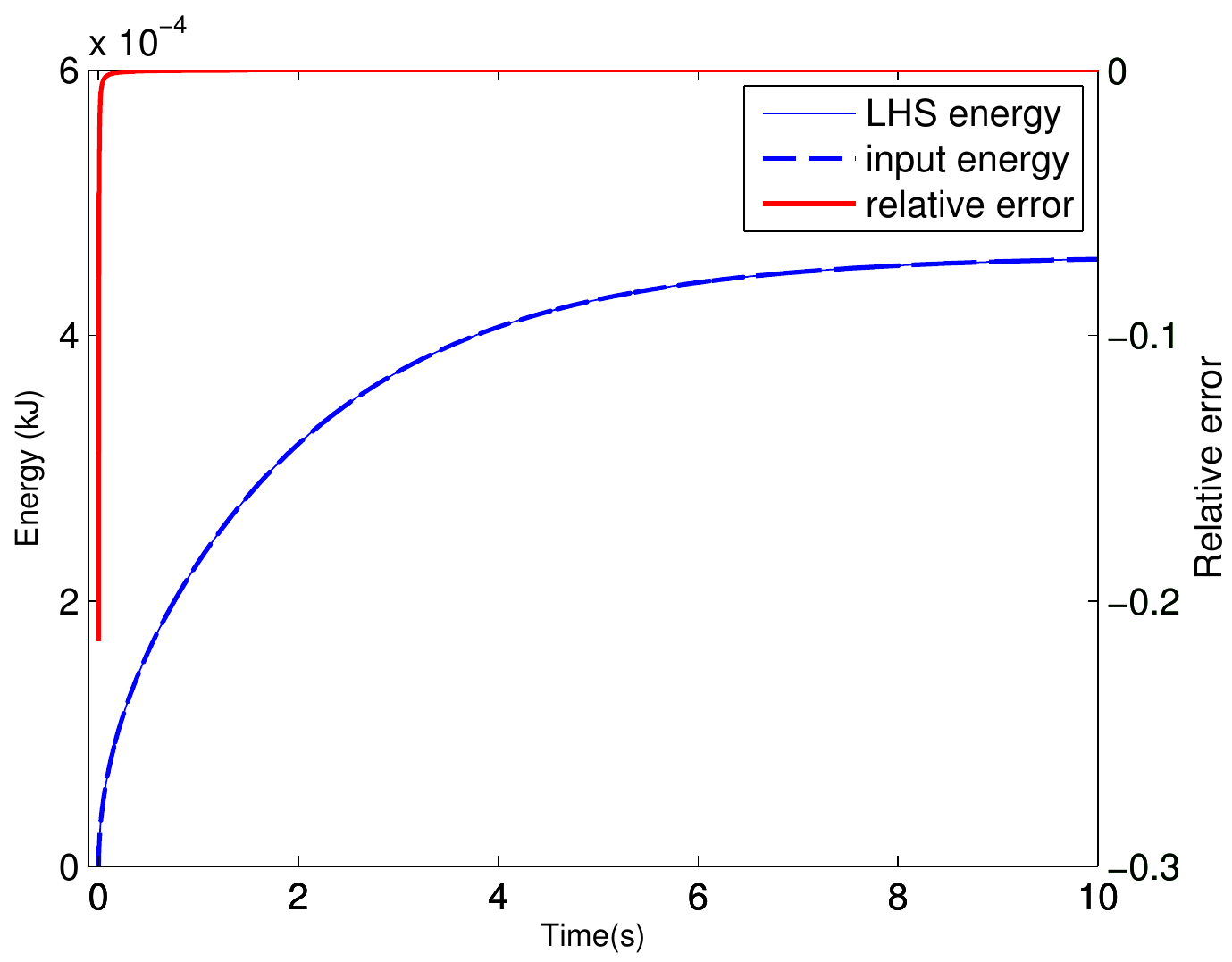}}
        \label{fig:example1_engyErr}
    }
    \hspace{0.05\textwidth}
    \subfloat[]
    {
        \resizebox{0.4\textwidth}{!}{\includegraphics{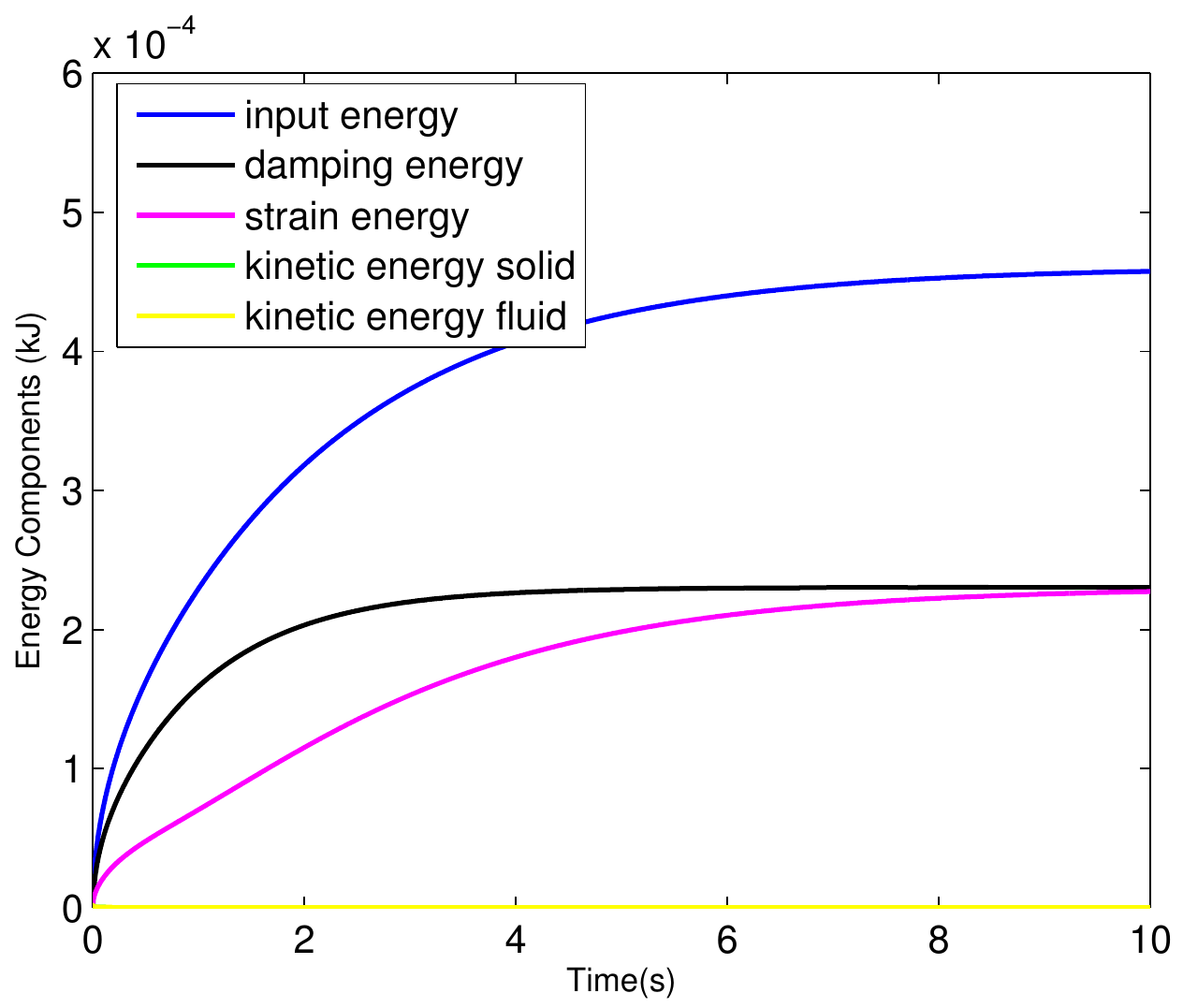}}
        \label{fig:example1_engComp}
    }
    \caption[Example 1 --- balance of energy]{Example 1 --- (a) Energy balance error. (b) Components of energy as defined in \eqref{eq:energy4}}
    \label{fig:example1_2}
\end{figure}

\mbox{}\\We have developed the energy balance equation for the entire porous medium in section \ref{EquationEnergy}. As a means of highlighting the performance and accuracy of the
presented time integration approach, the relative error in the energy balance of the solution is computed as shown in Figure \ref{fig:example1_engyErr}. Left hand side (LHS) energy is defined as the summation of all the energy terms except the input energy as expressed in \eqref{eq:energy4}. The relative error is defined as the absolute value of the difference of the LHS energy and input energy divided by the input energy. The error is very small and the time step scheme is accurate. Also all the components of energy are illustrated in Figure \ref{fig:example1_engComp}.
\subsubsection{Summary of findings}
\begin{itemize}
\item As shown in Figure \ref{fig:example1_disp}, a good agreement between the numerical and the analytical results is observed for both the transient behavior given by de Boer et al. \cite{deBoer_etal_1993} for wave propagation analysis, and the long term behavior given by Terzaghi's consolidation equation.
\item We observe a decreasing trend of wave speed as it travels through the medium.
\item Pressure diffusion begins right after the passage of the stress wave as shown in Figure \ref{fig:example1_prs}.
\item Due to the instantaneous nature of the step loading, high frequency numerical modes are excited early in time. These are shown in the velocity time history of Figure \ref{fig:example1_velZoom}. We suggest using a small ramp or a time integration method with algorithmic dissipation to alleviate this problem. Also a mass lumping method is shown to be effective in reduction of these types of numerical artifacts. We note that these high frequencies are not apparent in displacement time history.
\item The use of the Lobatto quadrature mass lumping method gives rise to an unbounded and unstable solution as shown in Figure \ref{fig:vel_unstable}.
\item The small relative error in the energy balance of numerical solution expresses the satisfactory performance and accuracy of the presented time integration method.
\end{itemize}
\subsection{Example 2: soil block with a partially loaded surface} \label{Example2_soilBlock}
\begin{figure}
    \centering
	 \resizebox{0.6\textwidth}{!}{\includegraphics{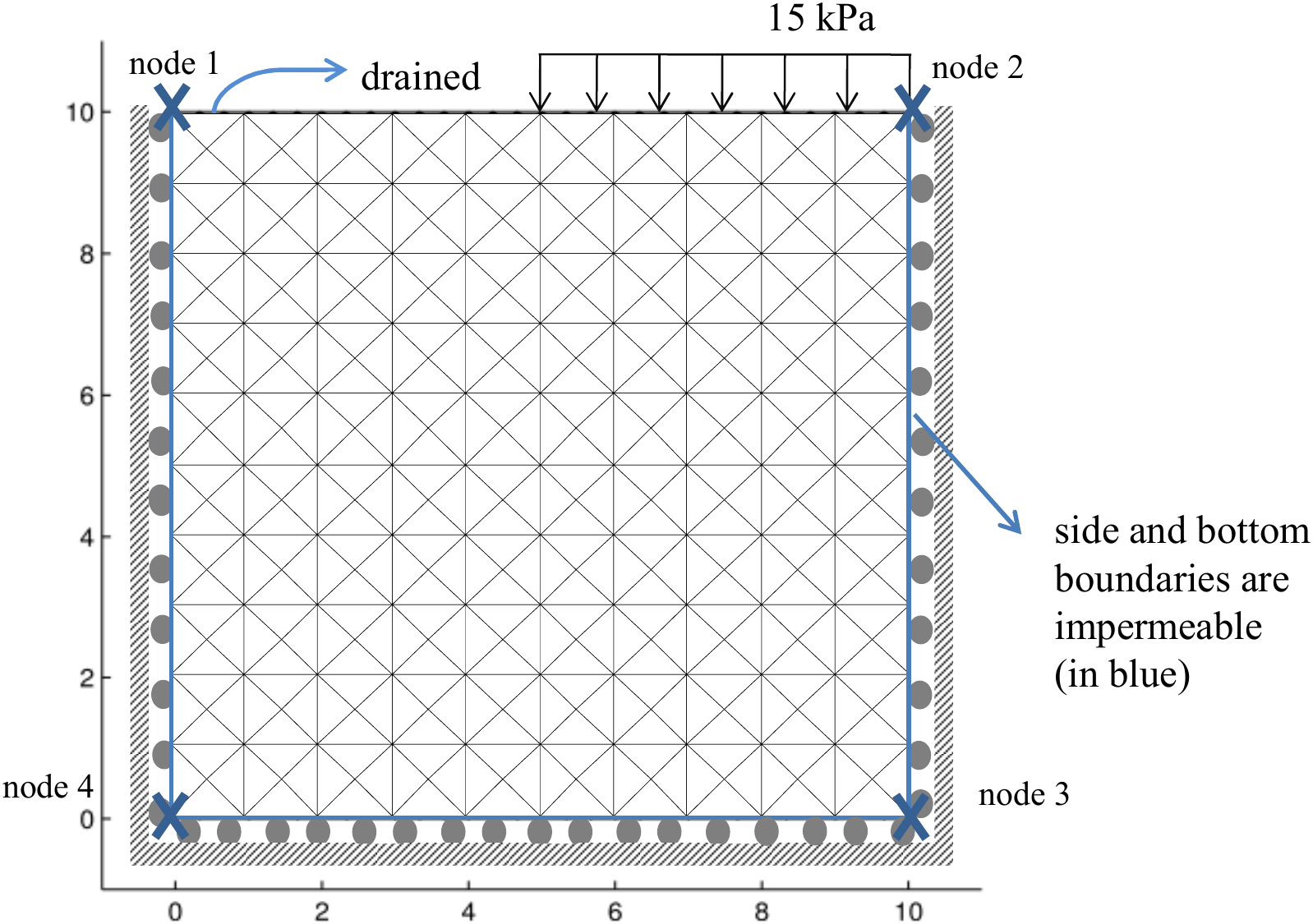}} 
	 \caption{Example 2 --- geometry and finite element mesh}
   	 \label{fig:example2_model}
\end{figure}

A key benchmark for verification of numerical scheme is the analysis of a saturated soil block with a partially loaded surface as shown in Figure \ref{fig:example2_model}. Due to symmetry only half of the problem is modeled.
All the boundaries except the unloaded area of the upper boundary are impermeable. The displacements normal to the surfaces of the side walls and the bottom are constrained. We adopt the material properties that are used by Li et al. \cite{Li_etal_2004} and shown in Table \ref{tab:Material_properties} which are also similar to those used by Diebels and Ehlers \cite{ Diebels_Ehlers_1996}. Two different values of hydraulic conductivity is considered to cover the typical and the small range of this parameter. The model deforms in plane strain under a vertical step load of 15 $\frac{\text{kN}}{\text{m}^2}$ amplitude. The time increment is fixed at the value of $5 \times 10^{-3}$. Other time steps including $\Delta t=1 \times 10^{-3}, 1 \times 10^{-4}$ are also tested but they appear to have no significant influence on the response. 

\noindent With this example as well, we pursue many goals as listed below:
\begin{enumerate}
	\item Verification with boundary element solutions
         \begin{itemize}
			\item Significance of finite element parameters: We show that the dynamic behavior of porous media with small value of hydraulic conductivity is very sensitive to the choice of finite element, mesh pattern, and mesh size. 
         \end{itemize}
	\item Damping and frequency content: the dynamic response of the medium with different material parameters is studied in terms of damping and frequency content. 
	\item Fluid and pressure diffusion: to check the smoothness of the variable field results, the passage of flow through the porous medium and the pressure contours are illustrated at different moments of time.
\end{enumerate}
\subsubsection{Verification with boundary element solutions}
\mbox{}\\Shown in Figure \ref{fig:example2_dispComparison} are the vertical displacement time histories at node 1 and 2 (the left and right upper corner nodes shown in Figure \ref{fig:example2_model}). Our numerical results with a crisscross pattern of meshing (shown in black color) match exactly with solutions obtained using an implementation of the boundary element formulation of Dargush and coworkers \cite{Chen_Dargush_1995,Apostolakis_Dargush_2012}. The boundary element solutions are not shown in the figure because they overlap our solutions. We remark that several finite element parameters such as mesh pattern, mesh size, and type of element can significantly affect numerical results in terms of dissipation and frequency content. This topic is discussed in the next section \ref{Poro_EffectMesh}.
\subsubsection{Significance of finite element parameters} \label{Poro_EffectMesh}
\begin{figure}
    \centering
    \subfloat[]
    {
        \resizebox{0.5\textwidth}{!}{\includegraphics{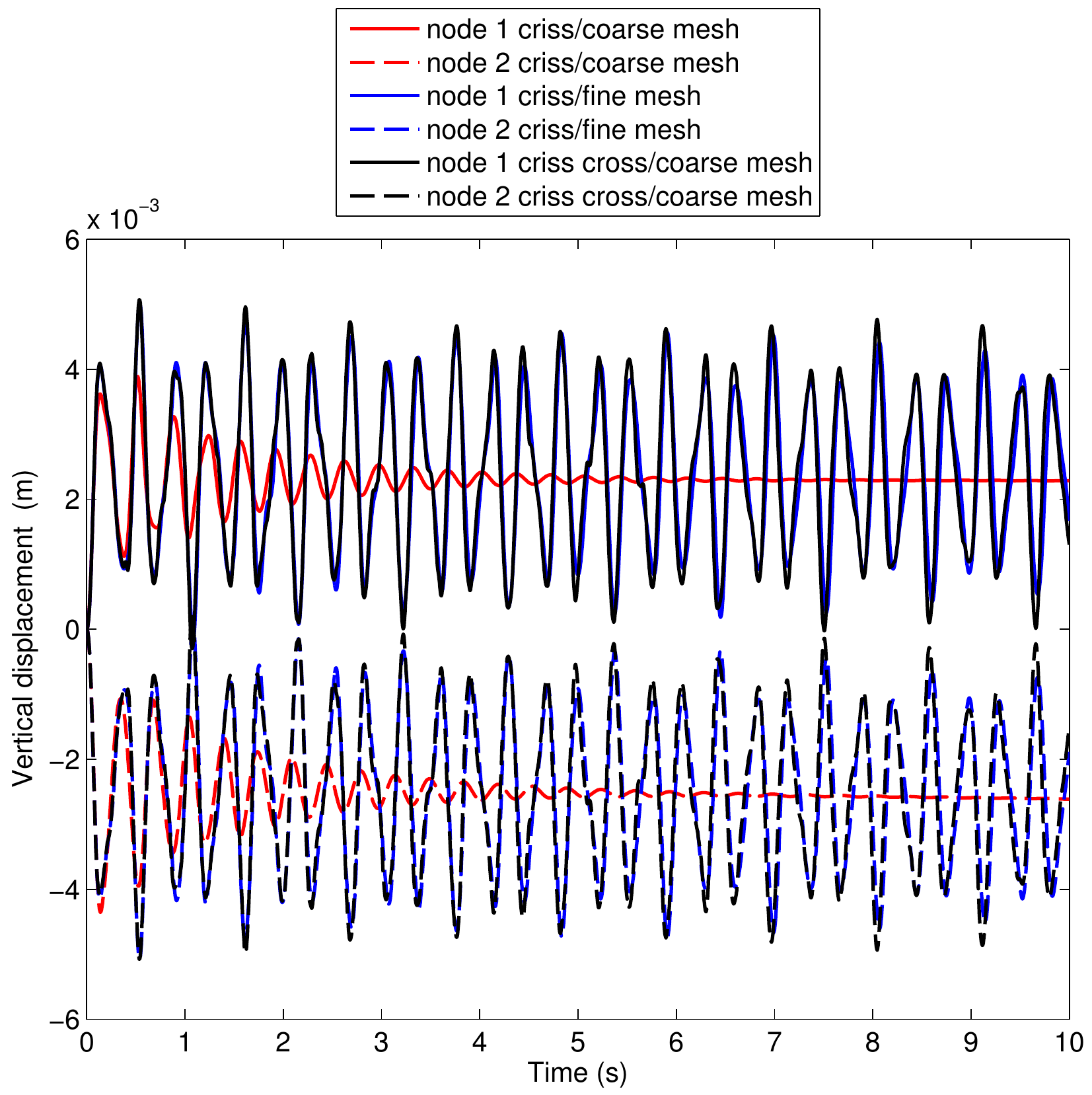}}
        \label{fig:vel_lowPermeab_2d}
    }
    \subfloat[]
    {
        \resizebox{0.5\textwidth}{!}{\includegraphics{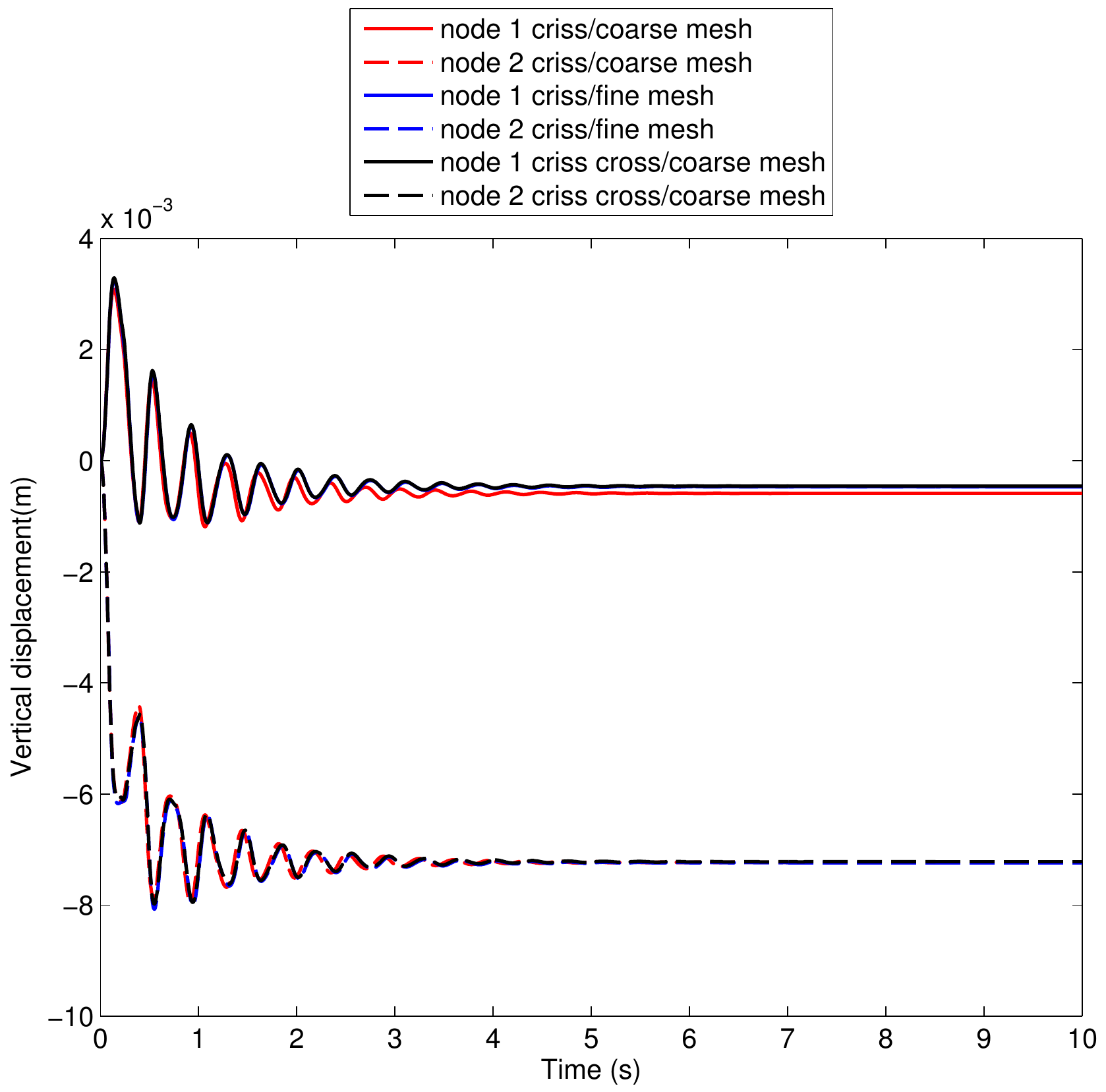}}
        \label{fig:vel_highPermeab_2d}
    }
    \caption[Example 2 --- displacement time history] {Example 2 --- displacement time histories for two different hydraulic conductivities. (a) Low hydraulic conductivity $\kh=10^{-4} \ms$: numerical results are sensitive to finite element parameters such as mesh pattern and size. For example criss pattern with coarse mesh leads to erroneous results in terms of dissipation and frequency content (in red). Most of literature studies report such results. Mesh refinement fixes the problem (in blue). (b) Normal hydraulic conductivity $\kh=10^{-1} \ms$: numerical results are not very sensitive to finite element parameters; criss pattern with a coarse mesh did slightly skews the settlement and changes the frequency, whereas crisscross pattern works properly under a coarse mesh.}
    \label{fig:example2_dispComparison}
\end{figure}

\mbox{}\\In order to demonstrate the effect of different parameters on the dynamic response of porous media, we analyze several types of mesh patterns (see Figure \ref{Meshpatterns_problem}) and sizes over our two finite elements, i.e., $\linprt$ and $\quadprt$. The results are summarized in Table \ref{tab:summary_elementType} and are shown in Figure \ref{fig:example2_dispComparison} for the case of $\linprt$ element. 
An important observation is made in Figure \ref{fig:vel_lowPermeab_2d} for porous media with small hydraulic conductivity value -- when using a criss pattern for griding, a coarse mesh would cause a dramatic numerical error in terms of dissipation and frequency content (shown in red color). However, after sufficiently refining the mesh (shown in blue color), the numerical solution converges to the correct one (boundary element solution). We remark that most of the studies, such as \cite{Diebels_Ehlers_1996,Arduino_Macari_II_2001,Li_etal_2004} present results similar to those shown in red color possibly due to insufficient mesh refinement and using low accuracy methods. 

On the other hand, in case of a typical value of hydraulic conductivity $\kh=1 \times 10^{-1} \ms$, we observe that the dynamic behavior of the porous media is not significantly affected by aforementioned finite element parameters -- only combination of a coarse mesh with criss pattern could slightly contaminate the frequency content and the amount of settlement (see red color results in Figure \ref{fig:vel_highPermeab_2d}) which is expected due to an asymmetric pattern of meshing. In contrast to numerical results of the $\linprt$ element , the displacement time history results using the quadratic element $\quadprt$, are satisfactory for all types of mesh patterns and sizes regardless of the hydraulic conductivity value.
\begin{table}
            \centering
           \begin{tabular}{ccccp{5cm}}
            	\toprule
		element type & mesh pattern & mesh size & performance & Remarks \\   
	\midrule
		\multirow{4}{*}{$\linprt$} & \multirow{2}{*}{criss} & coarse& X & \footnotesize{$\kh=1 \times 10^{-4} \ms$: exhibit nonphysical dissipation and inaccurate frequency see Figure \ref{fig:vel_lowPermeab_2d}. $\kh=1 \times 10^{-1} \ms$: slightly contaminate the frequency content and the amount of settlement }\\
		&&fine& \checkmark& \\
		& \multirow{2}{*}{crisscross \& union jack} & coarse& \checkmark & \\
		&&fine& \checkmark& \\
  		  \midrule 
		\multirow{4}{*}{$\quadprt$} & \multirow{2}{*}{criss} & coarse& \checkmark & \\
		&&fine& \checkmark&{\footnotesize{In all four analysis, time history of velocity develop high frequency oscillations}} \\
		& \multirow{2}{*}{criss ross \& union jack} & coarse& \checkmark &{\footnotesize{while displacement time history is smooth}} \\
		&&fine& \checkmark& \\
            	\bottomrule
            \end{tabular}
            \caption{Example 2 --- performance of the method using different element types, mesh patterns and mesh sizes}
            \label{tab:summary_elementType}
\end{table}

; the black color results using the crisscross pattern perfectly match the boundary element solution that is not shown here

\subsubsection{Damping and frequency content -- effects of hydraulic conductivity}
\mbox{}\\Comparing Figures \ref{fig:vel_lowPermeab_2d} and \ref{fig:vel_highPermeab_2d}, it is seen that the frequency of the dynamic response of the porous medium remains almost the same for models with different values of hydraulic conductivity, suggesting that the wave speed is not affected by this parameter. In terms of damping, when hydraulic conductivity is larger, the response decays faster in time and the steady state consolidation solution is achieved over a short period of time (Figure \ref{fig:vel_highPermeab_2d}) due to the greater drainage and fluid diffusion. In contrast, in case of small hydraulic conductivity, oscillations at the corner nodes persist over a long period of time and the steady state consolidation solution is not achieved in that period. In this condition, a nearly undrained response (similar to a nonporous incompressible elastic body) is observed where both the corner nodes move in skew symmetric manner. Also such low amount of drainage does not allow for any damping as shown in Figure \ref{fig:vel_lowPermeab_2d}. This behavior can be explained by considering the damping energy relationship of Equation \subeq{eq:energy3}{2}. This energy $E^{\text{D}}$ is proportional not only to the damping coefficient but also to the fluid Darcy velocity. In case of small hydraulic conductivity $\kh$, although the damping coefficient increases, the fluid Darcy velocity drops to such an extent that results in the overall reduction of damping energy. Comparing the slope of the damping energy in Figures \ref{fig:example2__k1_engComp} and \ref{fig:example2__k1_engComp} reveals the same fact. 
\subsubsection{Fluid and pressure diffusion}
\begin{figure}
    \centering
    \subfloat[]
    {
        \resizebox{0.4\textwidth}{!}{\includegraphics{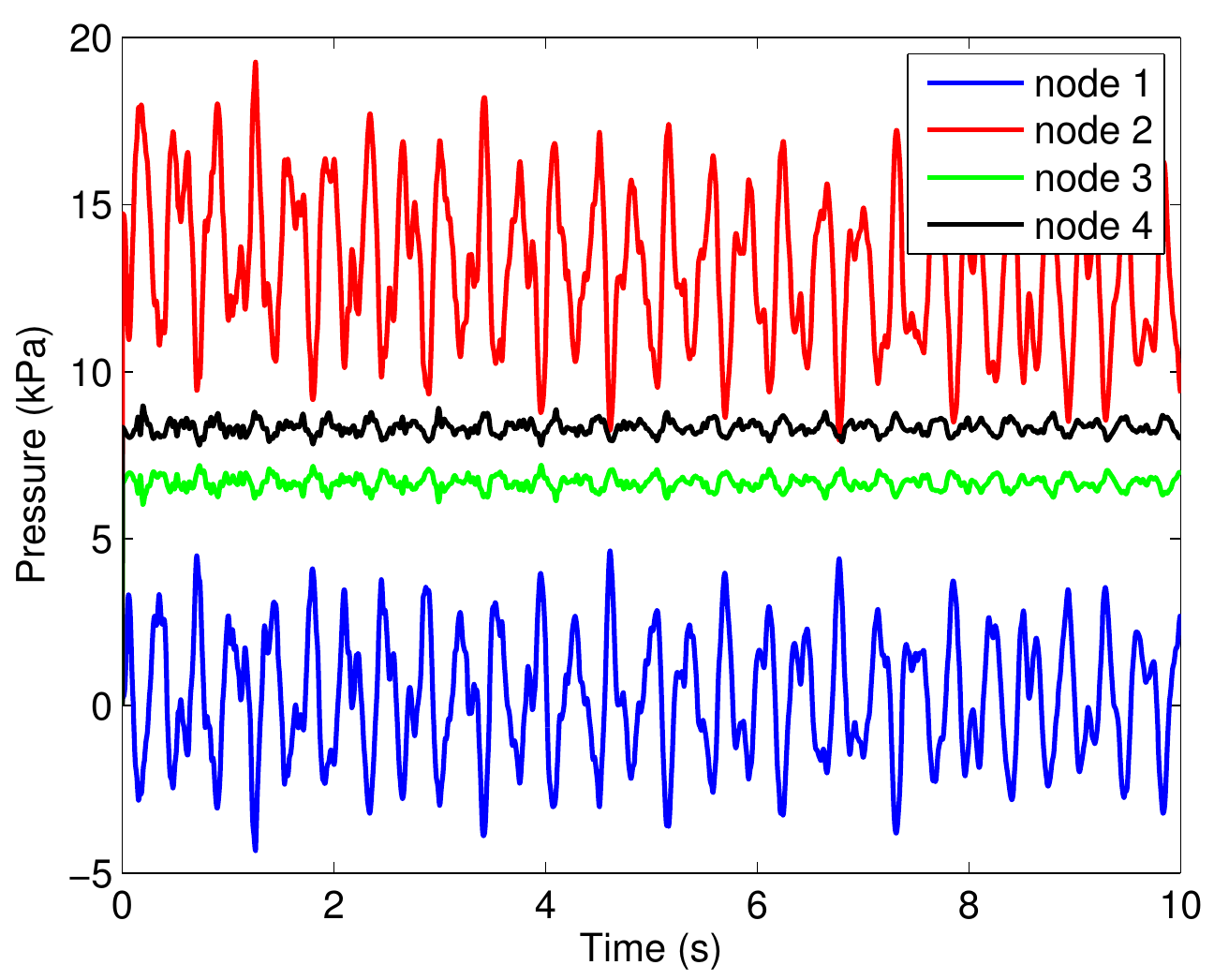}}
        \label{fig:example2_k4_presCritic}
    }
    \hspace{0.1\textwidth}
    \subfloat[]
    {
        \resizebox{0.4\textwidth}{!}{\includegraphics{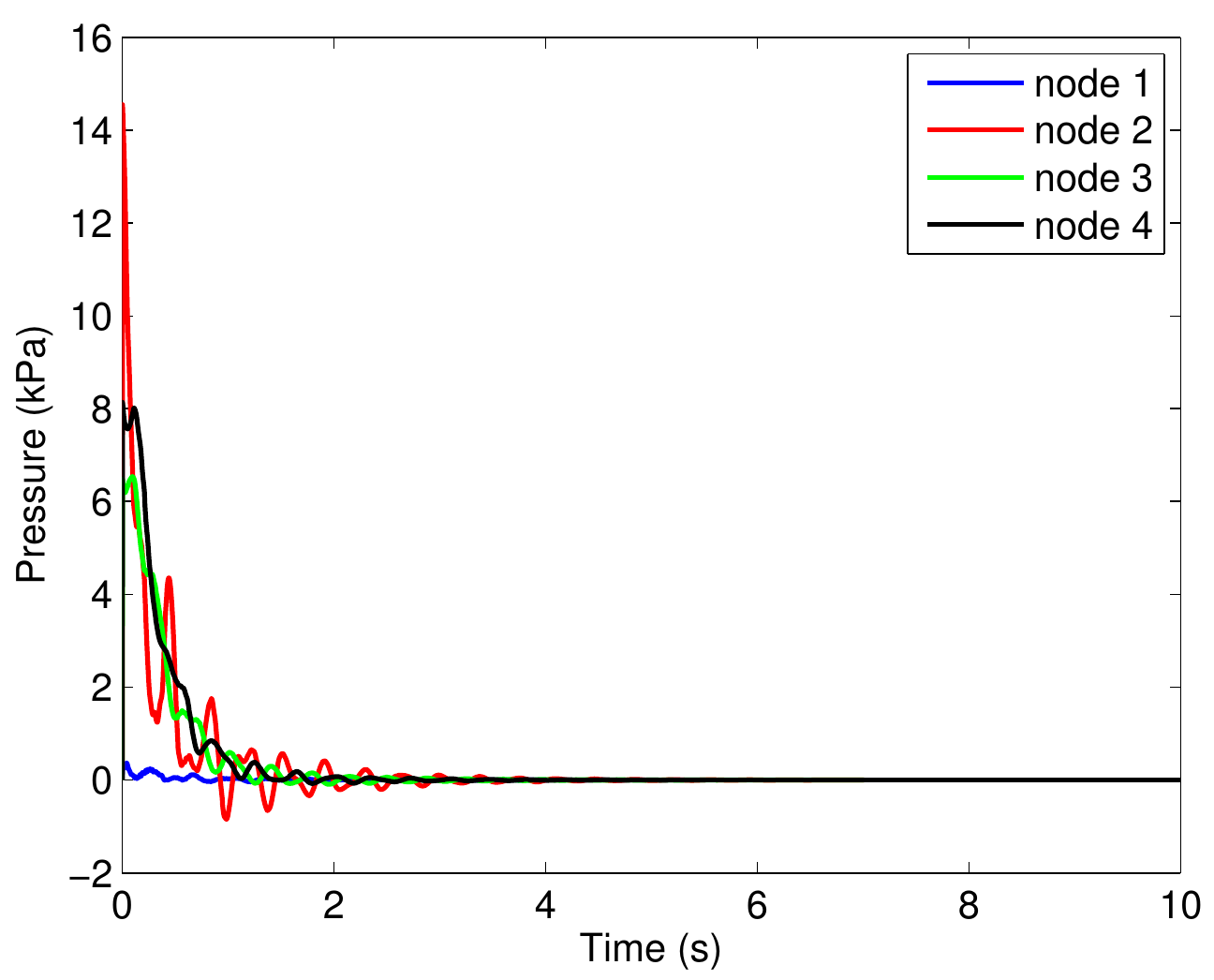}}
        \label{fig:example2_k1_presCritic}
    }
    \caption[Example 2 --- pressure time history]{Example 2 --- pressure time histories. (a) $\kh=1 \times 10^{-4} \ms$. (b) $\kh=1 \times 10^{-1} \ms$. }
    \label{fig:example2_presur}
\end{figure}

\begin{figure}
    \centering
    \subfloat[]
    {
      \resizebox{0.4\textwidth}{!}{\includegraphics{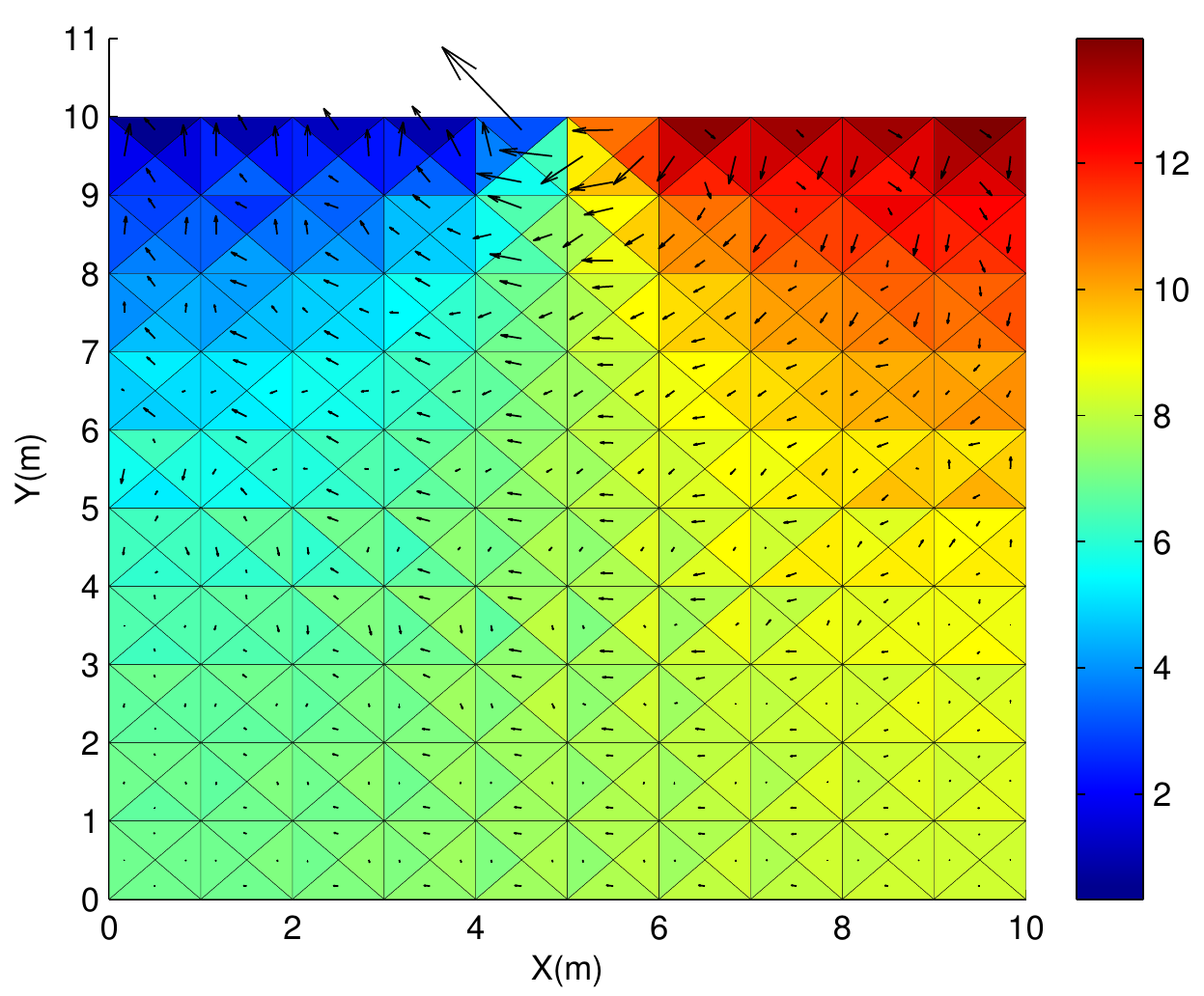}}
        \label{fig:example2_k4_presContour1}
   }
    \hspace{0.1\textwidth}
    \subfloat[]
    {
        \resizebox{0.4\textwidth}{!}{\includegraphics{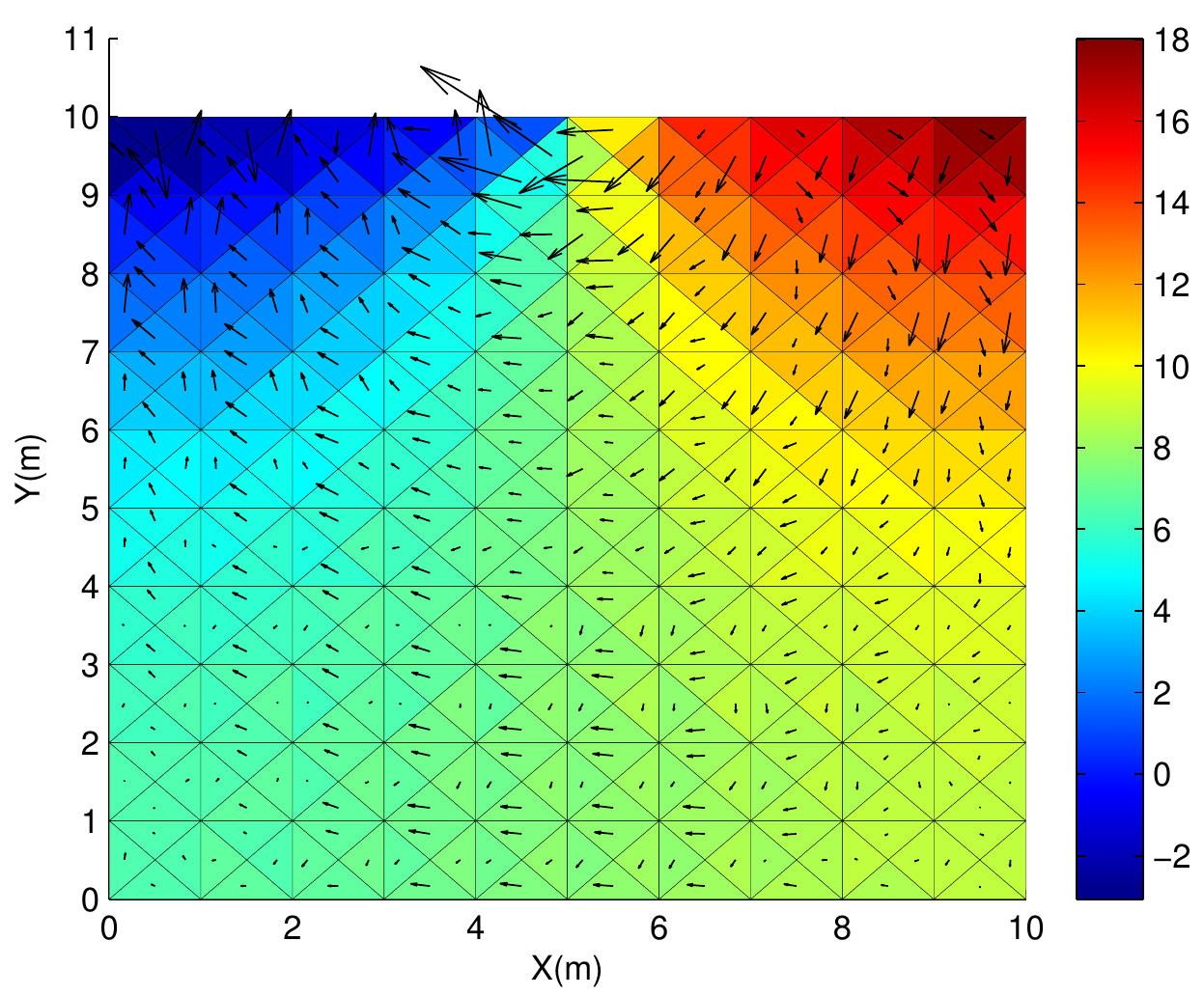}}
        \label{fig:example2_k4_presContour2}
    }
  \caption[Example 2 --- pressure contours and qualitative velocity vector field at $\kh=1 \times 10^{-4} \, \ms$] {Example 2 --- pressure contours and qualitative velocity vector field at $\kh=1 \times 10^{-4} \, \ms$. The results are smooth. Velocity field is larger close to the common loading and free boundary on top. (a) t=0.1 s. (b) t=0.9 s}
    \label{fig:example2_presurContourk4}
\end{figure}

\begin{figure}
    \centering
    \subfloat[]
    {
      \resizebox{0.4\textwidth}{!}{\includegraphics{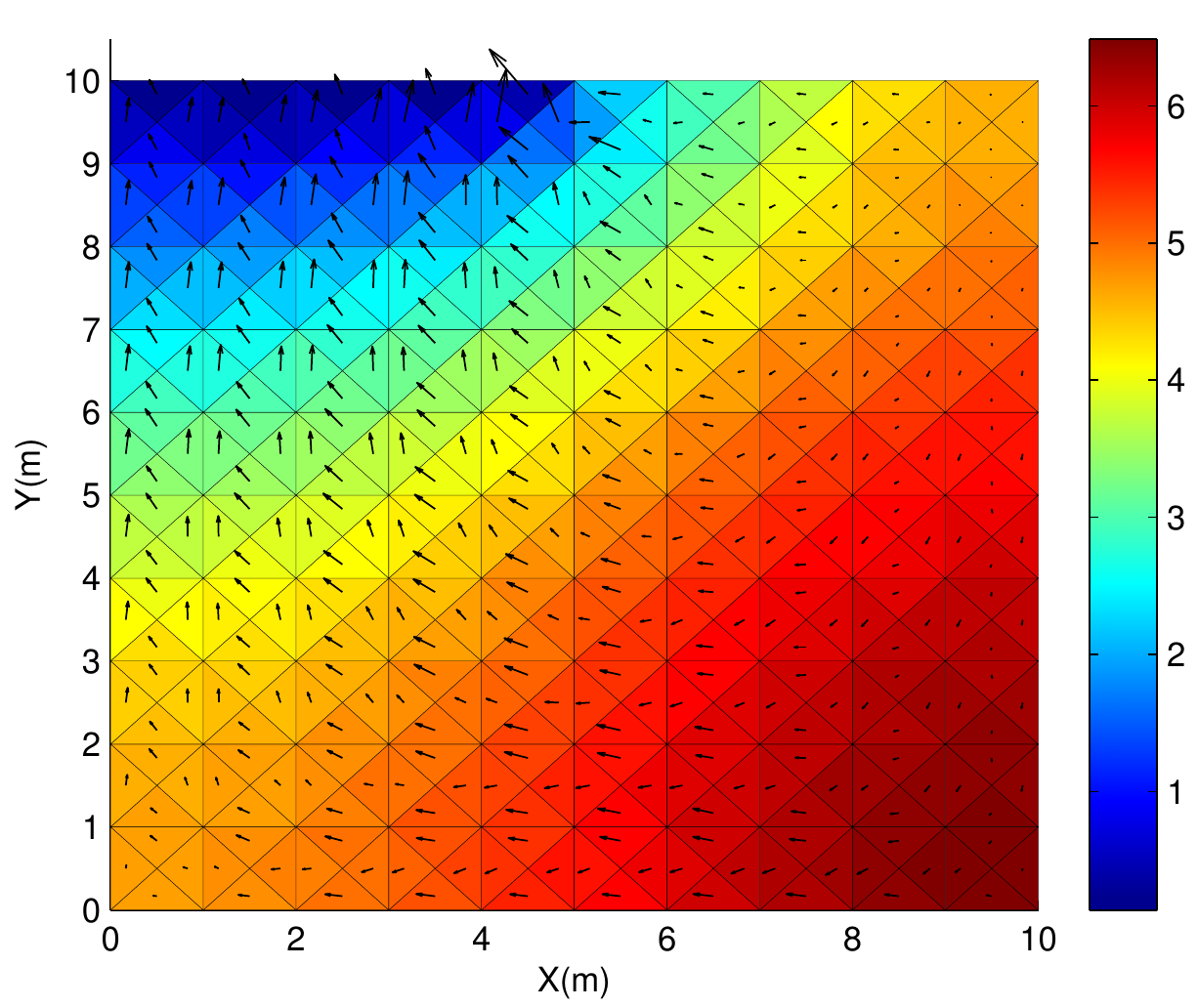}}
        \label{fig:example2_k1_presContour1}
   }
    \hspace{0.1\textwidth}
    \subfloat[]
    {
        \resizebox{0.4\textwidth}{!}{\includegraphics{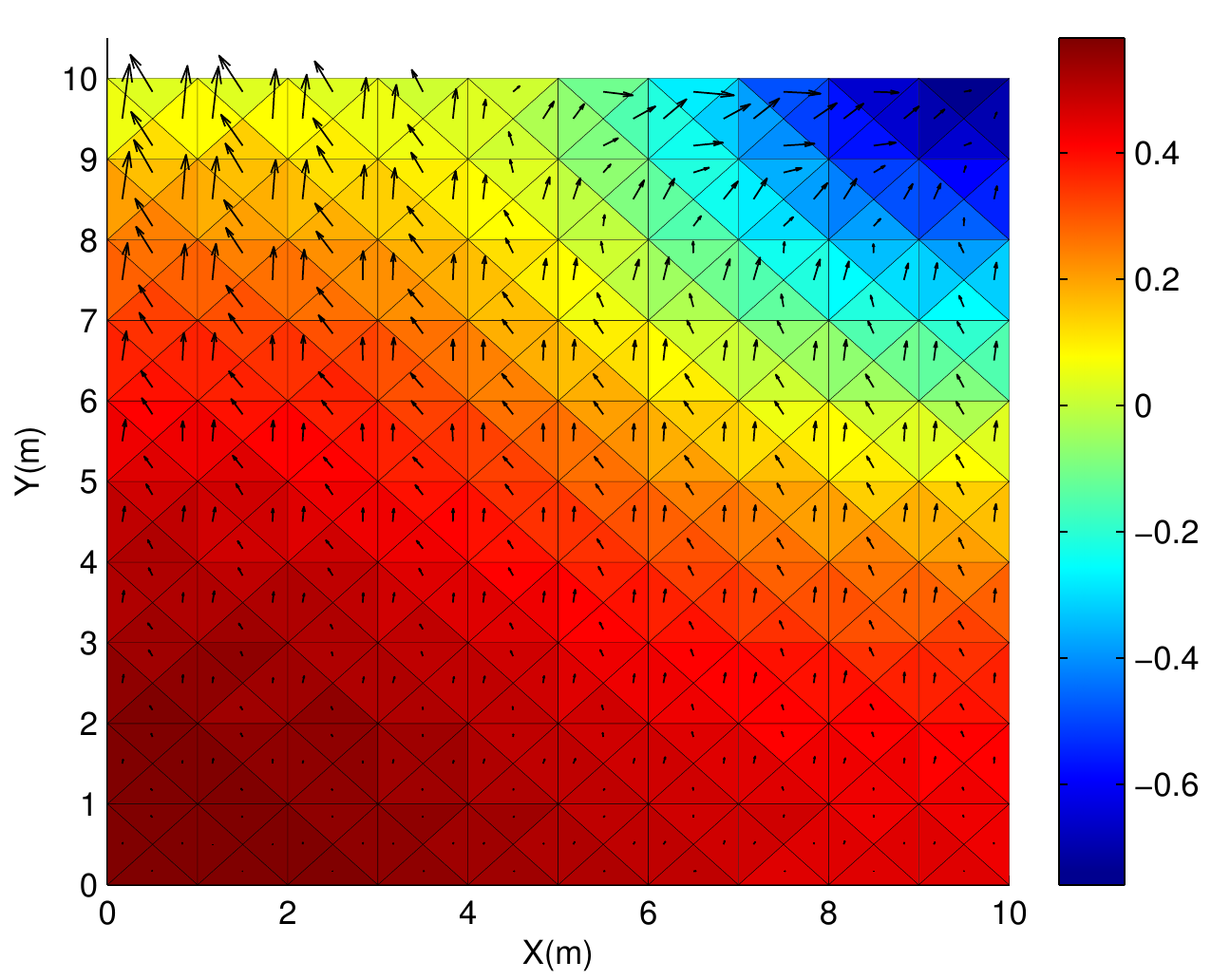}}
        \label{fig:example2_k1_presContour2}
    }
\\
  \subfloat[]
    {
      \resizebox{0.4\textwidth}{!}{\includegraphics{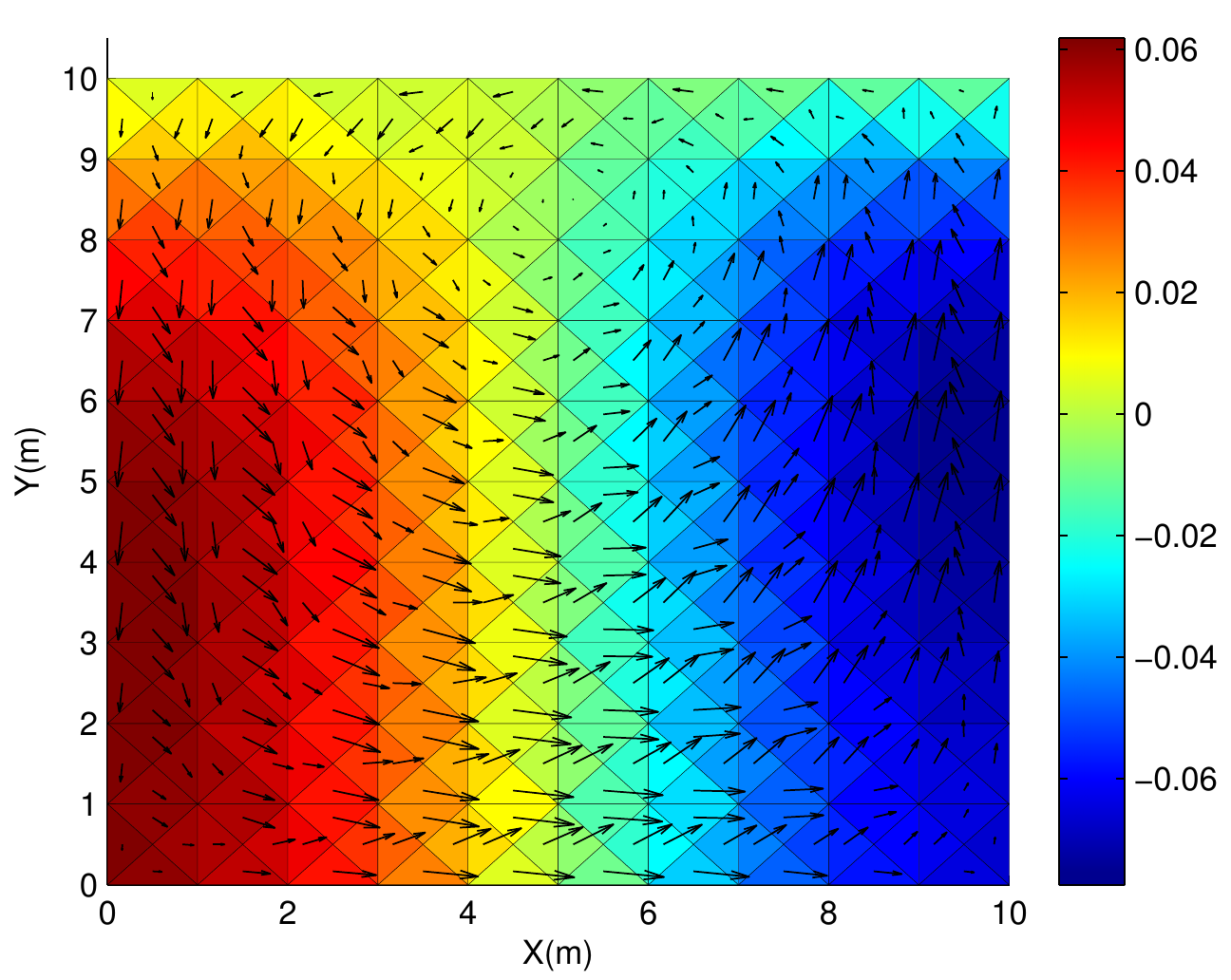}}
        \label{fig:example2_k1_presContour3}
   }
    \hspace{0.1\textwidth}
    \subfloat[]
    {
        \resizebox{0.4\textwidth}{!}{\includegraphics{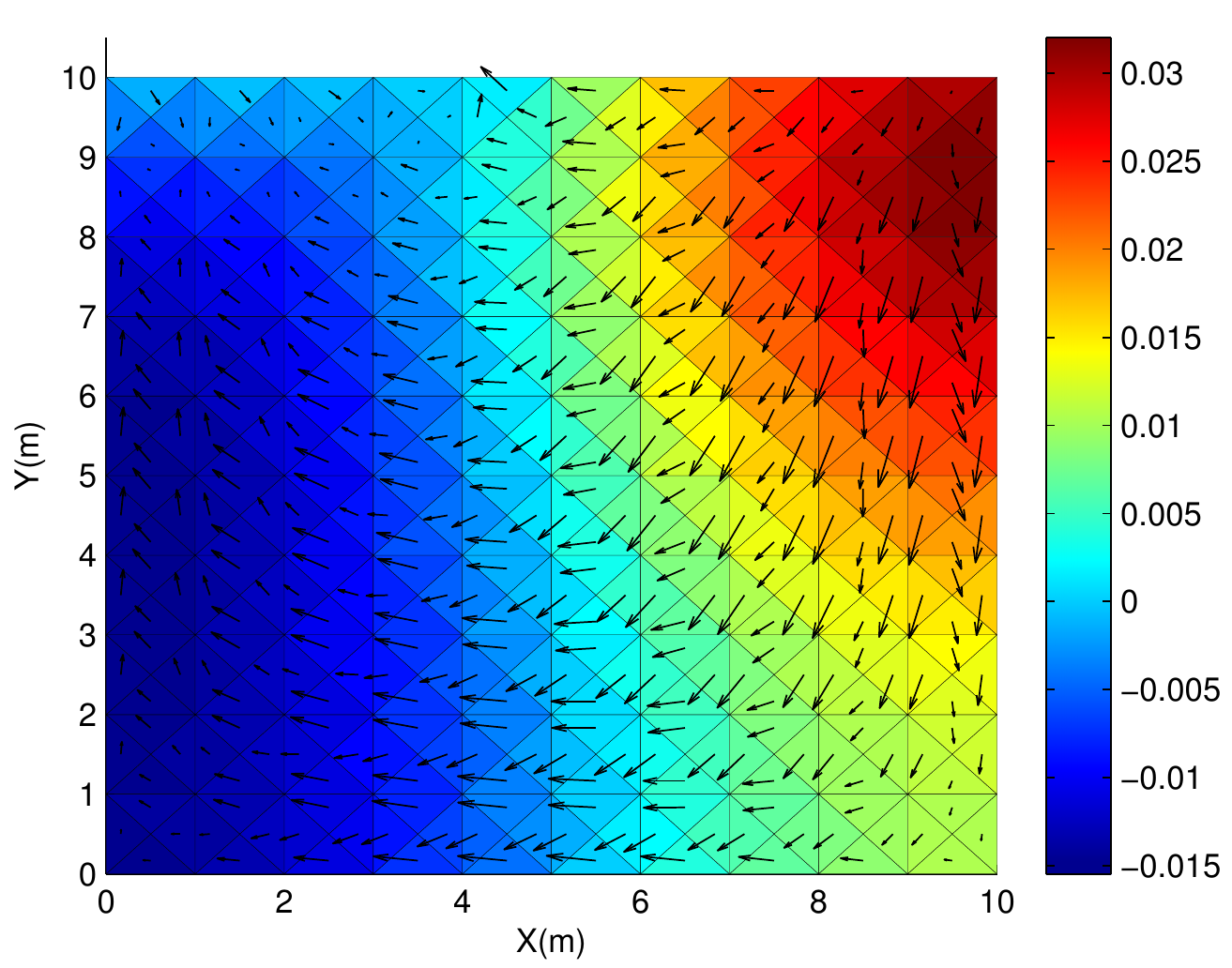}}
        \label{fig:example2_k1_presContour4}
    }
  \caption [Example 2 --- pressure contours and qualitative velocity vector field at $\kh=1 \times 10^{-1} \ms$]{Example 2 --- pressure contours and qualitative velocity vector field at $\kh=1 \times 10^{-1} \ms$. (a) t=0.2 s. (b) t=1 s. (c) t=2.2 s. (d) t= 3.8 s.}
    \label{fig:example2_presurContourk1}
\end{figure}

\mbox{}\\Time history of the pressure is plotted in Figure \ref{fig:example2_presur} for all the corner points. For a typical value of conductivity $\kh=1 \times 10^{-1} \ms$, Figure \ref{fig:example2_k1_presCritic} shows a general trend of pore pressure decay. The initial pressure build-up is apparent in this plot. On the other hand, the pressure time history for the small value of conductivity, Figure \ref{fig:example2_k4_presCritic}, shows a small amount of decay while oscillating with the same frequency as that of skeleton and fluid velocity.
Shown in Figures \ref{fig:example2_presurContourk4} and \ref{fig:example2_presurContourk1} are several snapshots at different moments of time of the pore pressure contours and qualitative fluid velocity vector fields. The velocities in these figures are computed at the centroid of each triangular element using the interpolation function in \eqref{eq:interpFunc}.

\begin{figure}
    \centering
    \subfloat[]
    {
        \resizebox{0.4\textwidth}{!}{\includegraphics{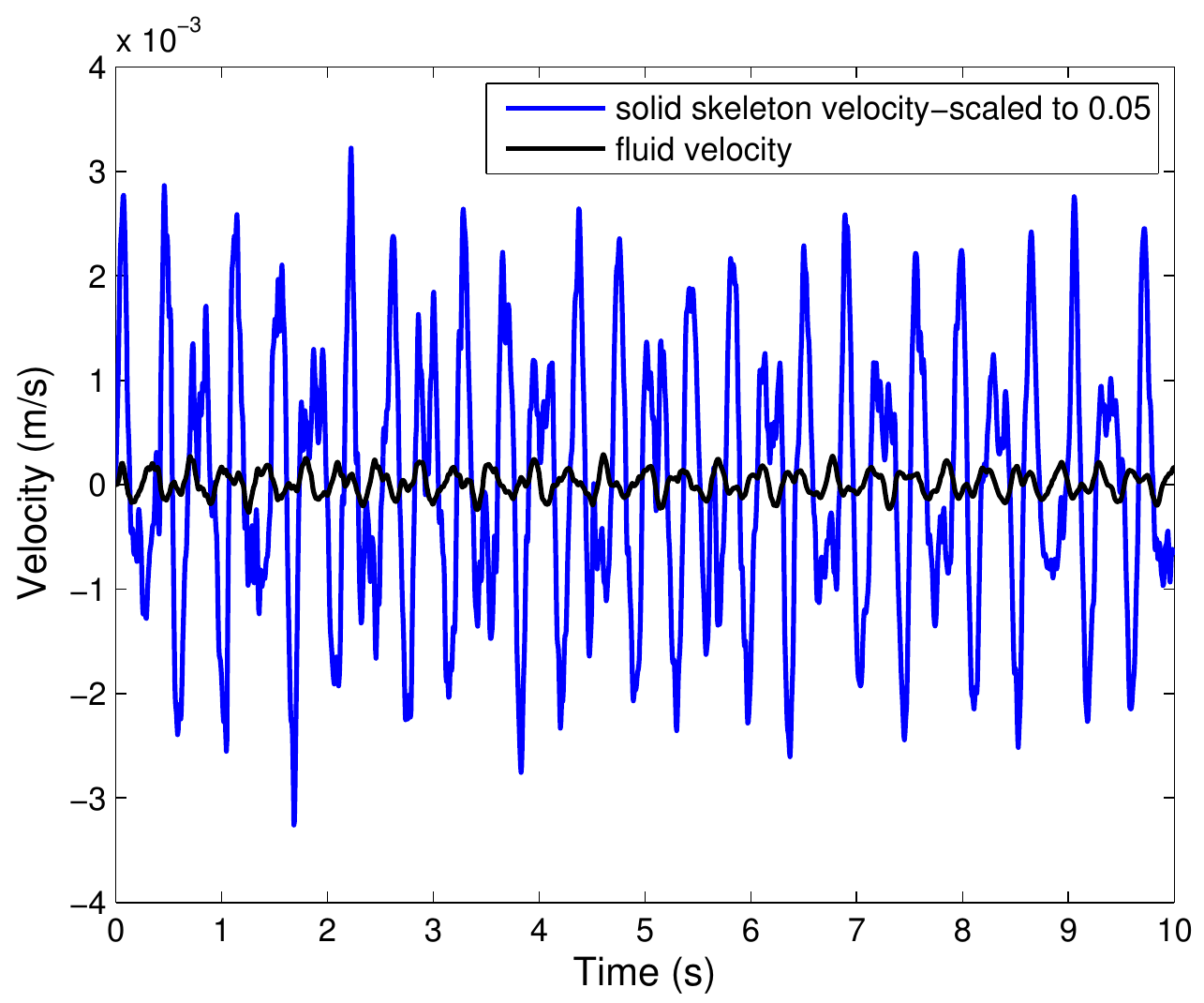}}
        \label{fig:velComparison_k4_2d}
    }
    \hspace{0.1\textwidth}
    \subfloat[]
    {
        \resizebox{0.4\textwidth}{!}{\includegraphics{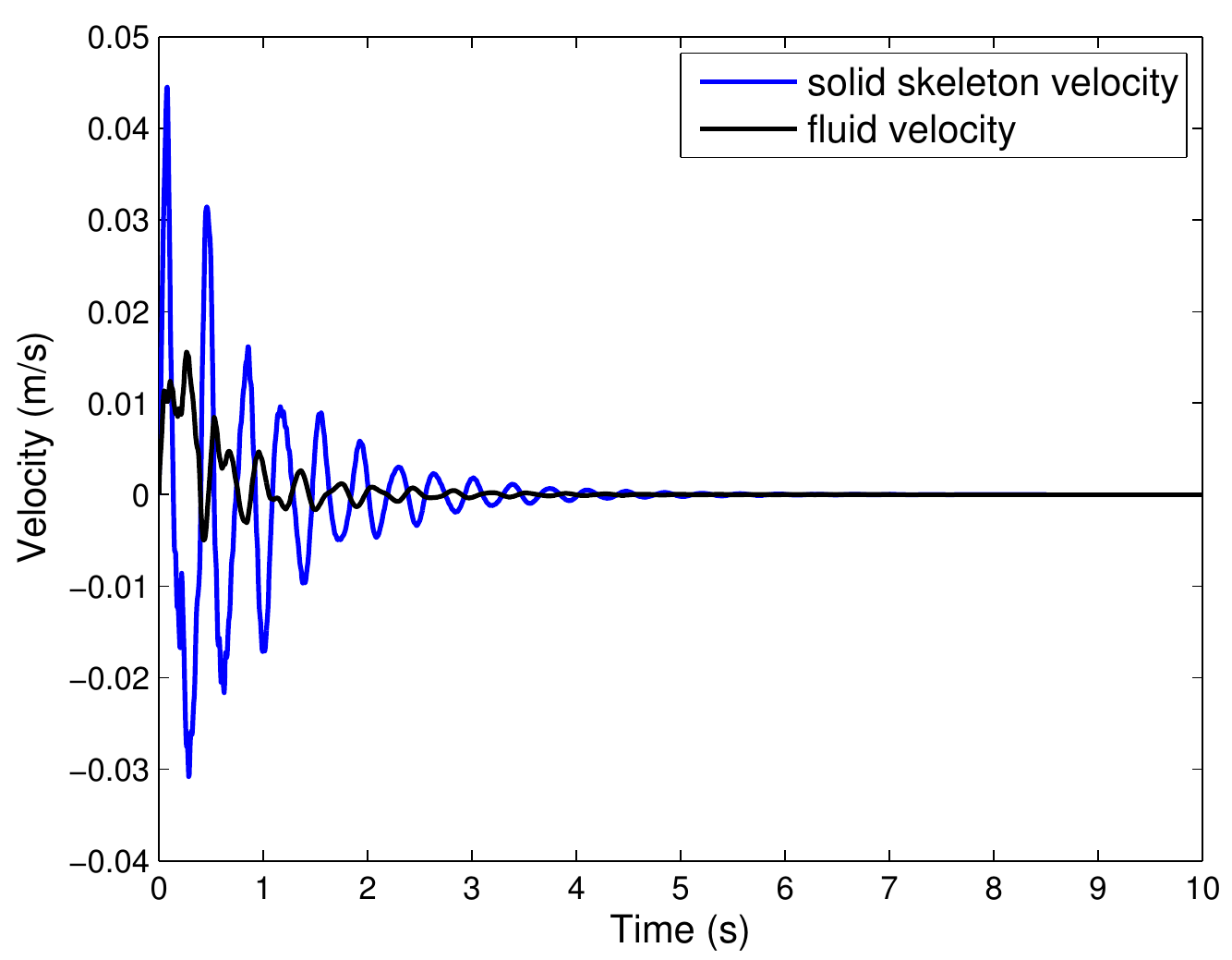}}
        \label{fig:velComparison_k1_2d}
    }
    \caption[identification of slow dilatational wave]{Example 2 --- skeleton and fluid Darcy velocity time histories, (a) $\kh=1 \times 10^{-4} \ms$, (b) $\kh=1 \times 10^{-1} \ms$. In an incompressible porous medium, the slow dilatational wave is characterized as moving of the fluid and solid in opposite phases.}
    \label{fig:velComparison_2d} 
\end{figure}

Figure \ref{fig:velComparison_2d} illustrates a comparison of fluid Darcy velocity and skeleton velocity. We observe that in an incompressible porous medium only the slow dilatational wave exists, which is characterized by motion of the fluid and skeleton with opposite phase angles \cite{Schanz_Pryl_2004}.

\begin{figure}
    \centering
    \subfloat[]
    {
        \resizebox{0.5\textwidth}{!}{\includegraphics{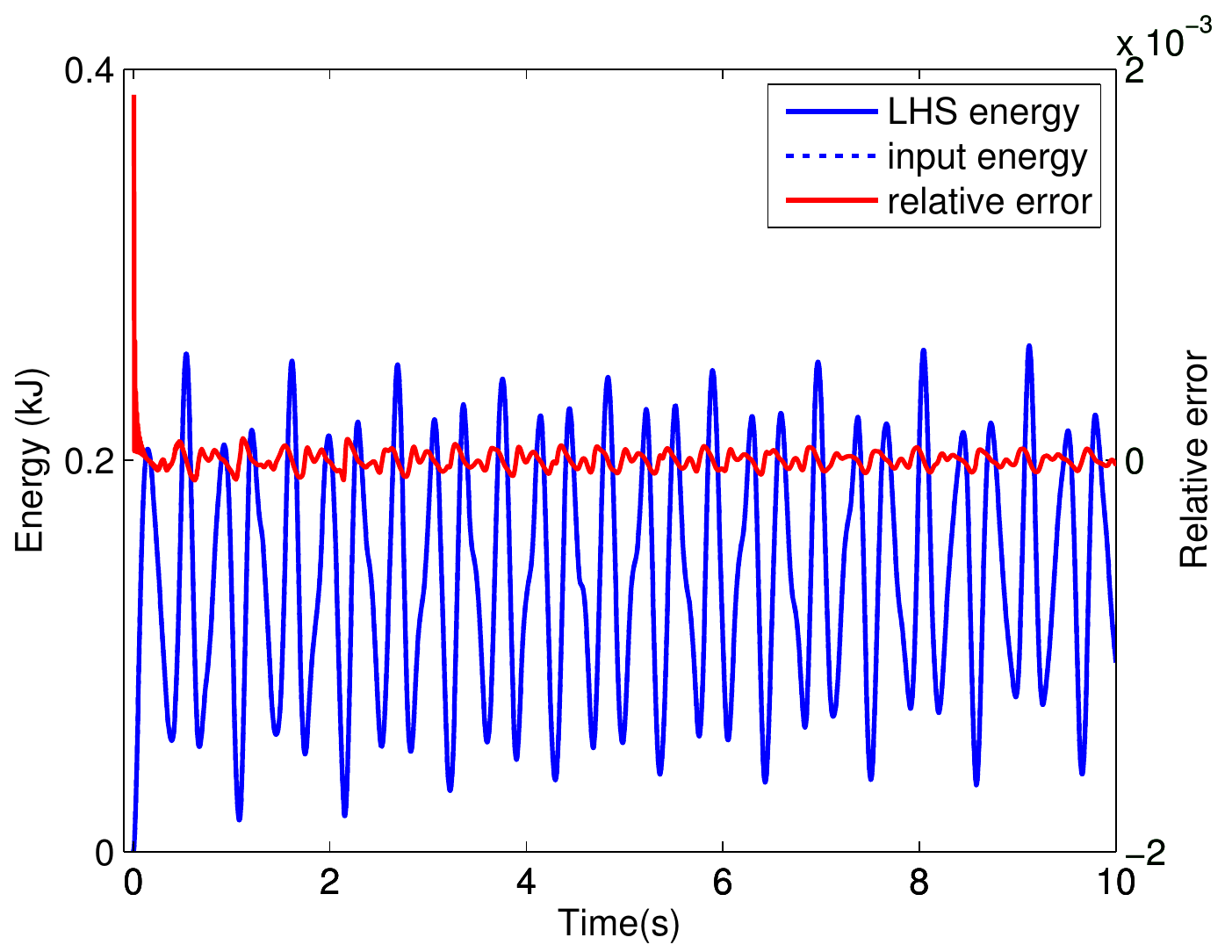}}
        \label{fig:example2_k1_engyErr}
    }
    \subfloat[]
    {
        \resizebox{0.46\textwidth}{!}{\includegraphics{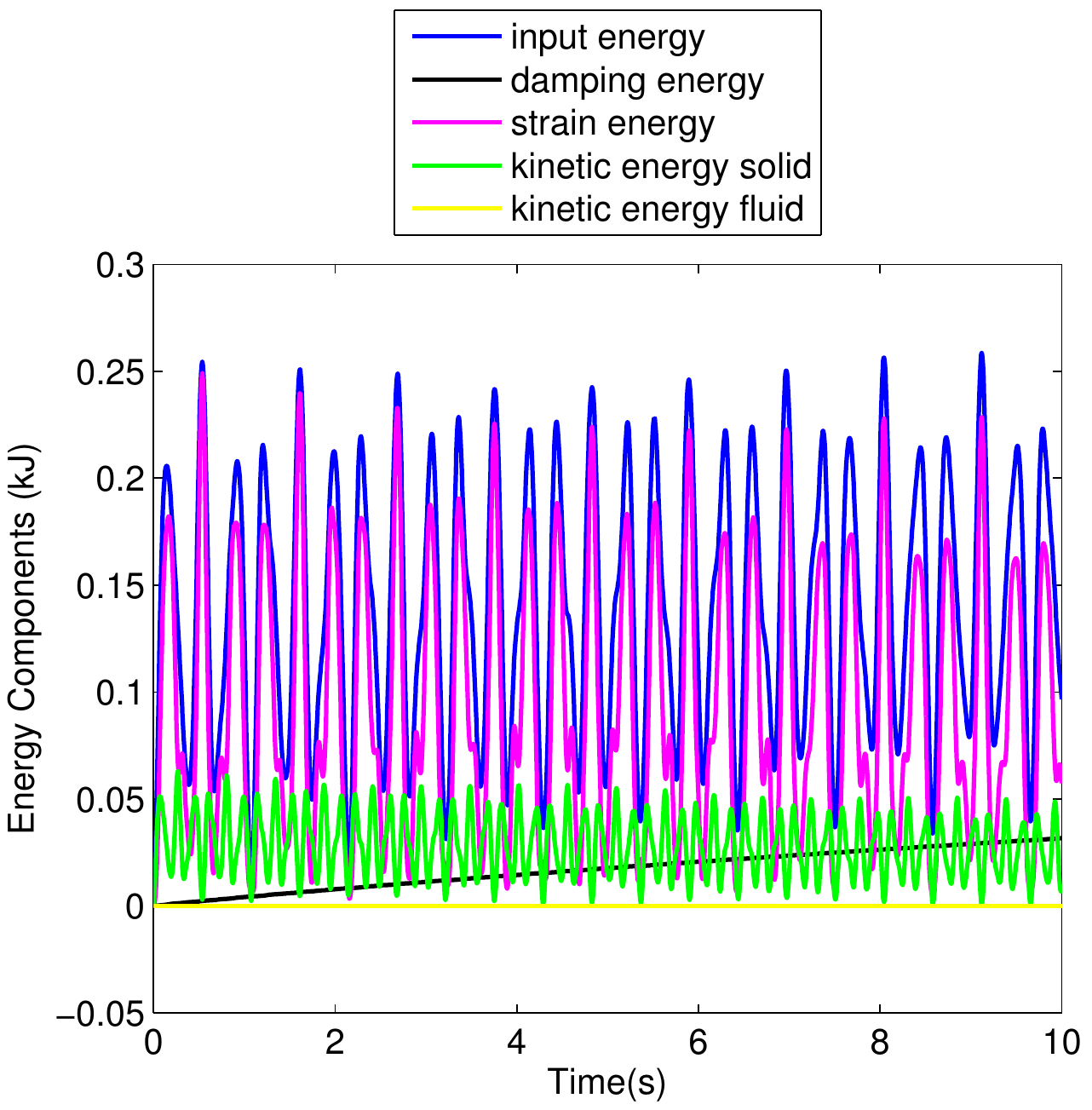}}
        \label{fig:example2__k1_engComp}
    }
\\
   \subfloat[]
    {
        \resizebox{0.5\textwidth}{!}{\includegraphics{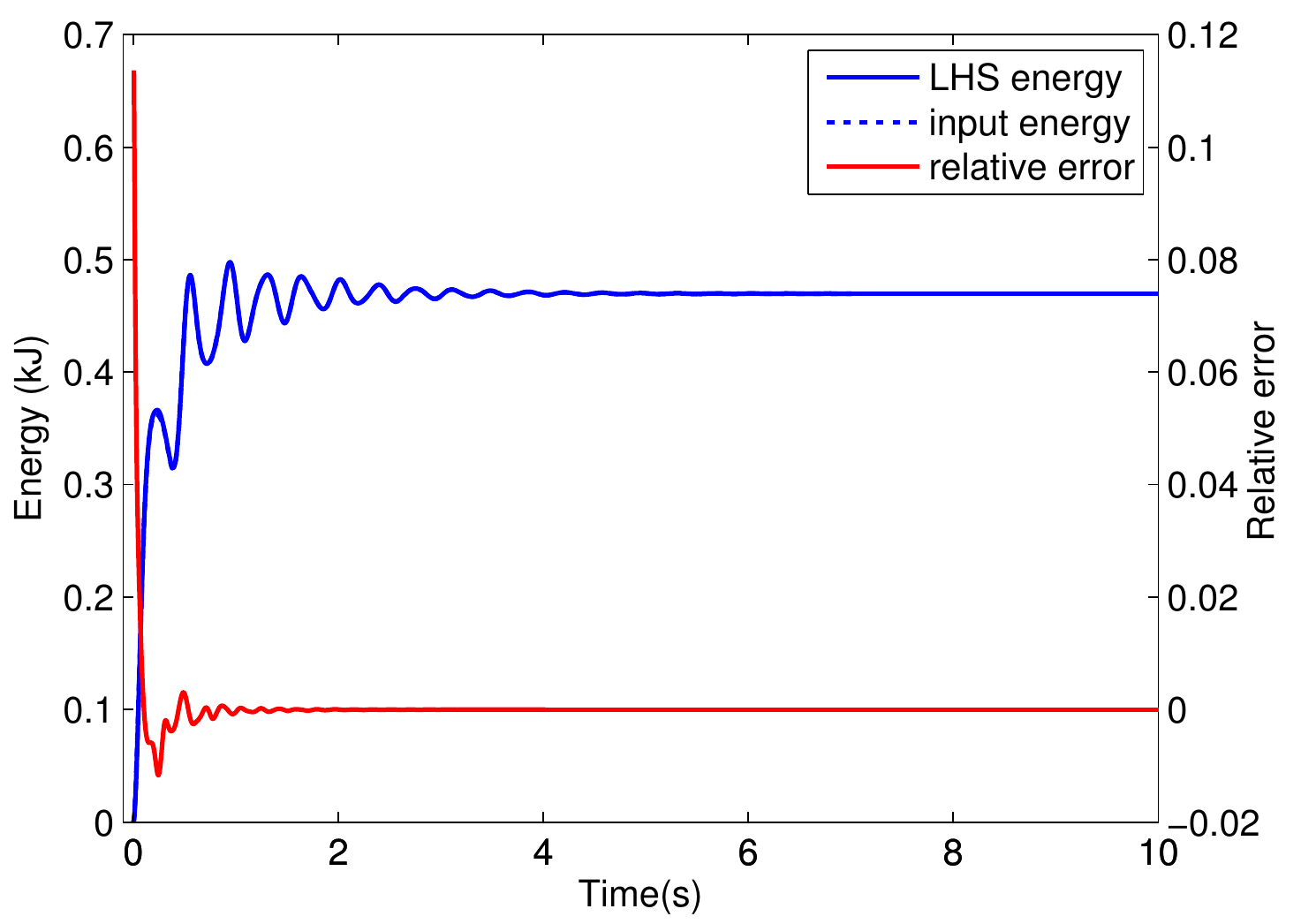}}
        \label{fig:example2_k4_engyErr}
    }
    \subfloat[]
    {
        \resizebox{0.46\textwidth}{!}{\includegraphics{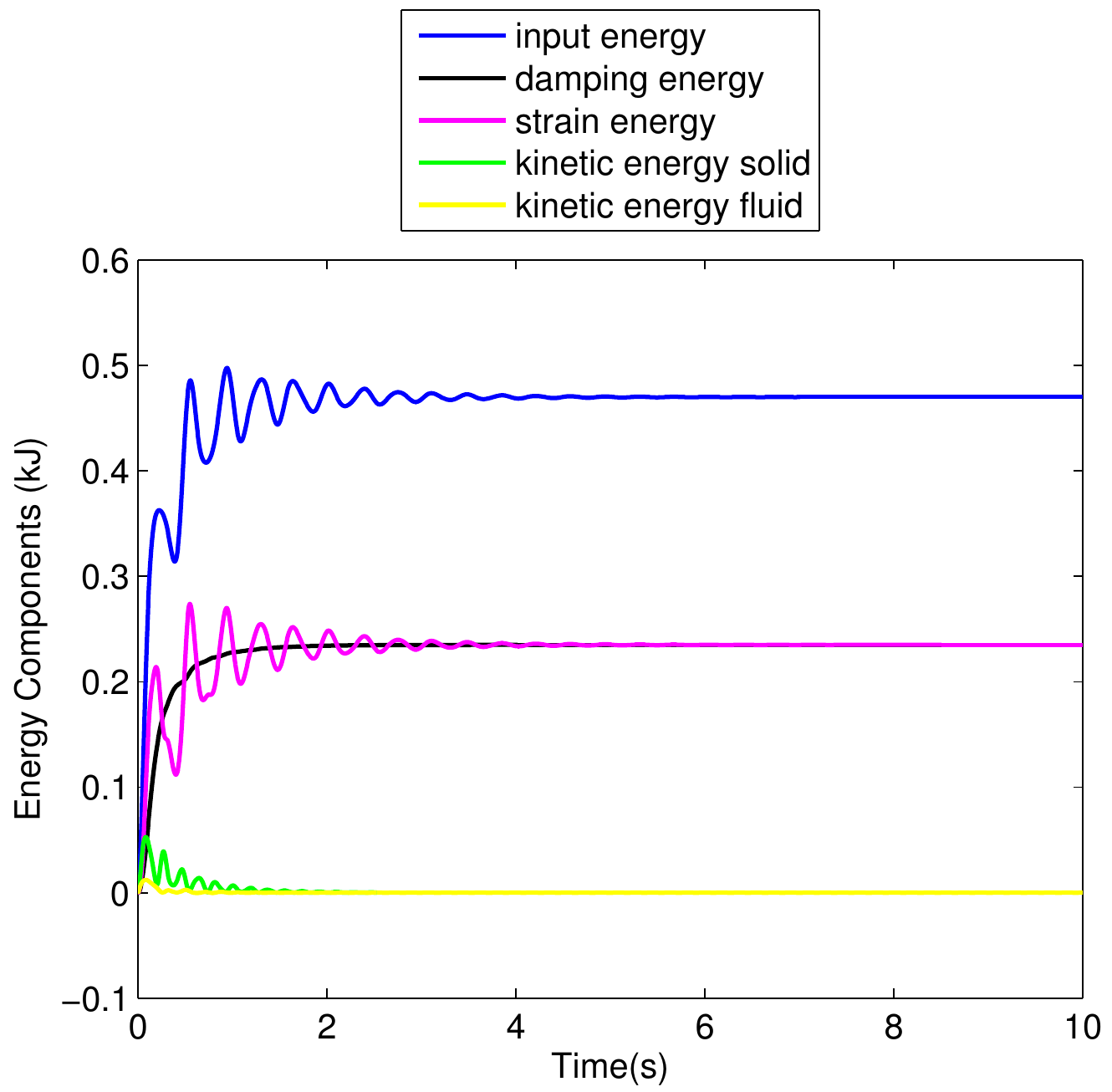}}
        \label{fig:example2_k4_engComp}
    }
    \caption[Example 2 --- energy plots]{Example 2 --- energy plots. (a,b) $\kh=1 \times 10^{-4} \ms$. (c,d) $\kh=1 \times 10^{-1} \ms$ }
    \label{fig:example2_engy}
\end{figure}

Furthermore, similar to Example 1, the satisfactory performance and accuracy of the time discretization technique is evaluated through the error in satisfying the energy balance of the system as shown in Figure \ref{fig:example2_engy}.
\subsubsection{Summary of findings}
\begin{itemize}
\item Numerical results are verified with boundary element solutions for a wide range of values of hydraulic conductivity (Figure \ref{fig:example2_dispComparison}).
\item The numerical simulation of dynamics of porous media with small value of hydraulic conductivity can be significantly biased with with an unsuitable choice of finite element method, size, or mesh pattern (Figure \ref{fig:vel_lowPermeab_2d}). Most studies in the literature, such as \cite{Diebels_Ehlers_1996,Arduino_Macari_II_2001,Li_etal_2004}, present erroneous results in terms of dissipation and frequency content possibly due to insufficient mesh refinement or low accuracy methods. These results in the literature suggest an \emph{increase} in damping with decreasing hydraulic conductivity, while a proper finite element discretization leads to the opposite result.
\item  The frequency of the response remains almost the same for different hydraulic conductivity values (Figure \ref{fig:example2_dispComparison}), suggesting that the wave speed is not affected by the value of the hydraulic conductivity. However, damping is shown to be dependent to the hydraulic conductivity value, i.e., the larger this value is, the response is damped out faster. In contrast, in case of small hydraulic conductivity (Figure \ref{fig:vel_lowPermeab_2d}) oscillations persist over a long period of time and a nearly undrained response (similar to a nonporous incompressible elastic body) is observed. 
\item  The slow dilatational wave in an incompressible porous medium is characterized by motion of the fluid and skeleton with opposite phase angles as shown in Figure \ref{fig:velComparison_2d}s
\end{itemize}
\subsection{Example 3: challenge of locking -- model with very low hydraulic conductivity} \label{sec:NumericalExmpl3}
As discussed in section \ref{LimitingCases}, there are two limiting cases at which the finite element must satisfy the stability criterion. The rigid skeleton limit requires that flow variables $(w,p)$ be compatible with each other in terms of LBB condition; Raviart-Thomas mixed finite element is theoretically proved to satisfy this condition. The other limiting problem is the very low hydraulic conductivity case ($\kh \rightarrow 0$) in which the system degenerates to an incompressible elastic system with $u$ and $p$ as the primary variables. For this limiting case, the finite element discretization of $(u,p)$ must comply with LBB criterion. In this example we aim at studying such a limiting case under a special configuration that typically induces spurious pressure mode and locking. In this work, the phenomena of locking refers to the situation when a finite element discretization is not able to pass the inf-sup condition (see discussions in section \ref{sec:inf-sup}). 

\noindent This example seeks to fulfill the following goals
\begin{enumerate}
	\item Incompressible elasticity problem: we first examine the stability of finite element spaces of $(u,p)$ in an incompressible elasticity problem $(\kh=0)$. We consider two elements, namely, the linear displacement-constant pressure $\linElem$ and the quadratic displacement-constant pressure $\quadElem$, and study each one on three mesh patterns as shown in Figure \ref{Meshpatterns_problem}.
		\begin{itemize}
			\item Spurious modes and inf-sup test: the performances of various finite elements are explored under different scenarios: element level (local), assemblage level (with boundary conditions applied), and the uniform stability (assessed with increasing mesh refinement).
		\end{itemize}
	\item Poroelasticity problem under very low hydraulic conductivity ($\kh$ very small)
		\begin{itemize}
			\item Checkerboard pattern in the pressure field: in a porous medium with very low hydraulic conductivity, numerical methods often produce checkerboard patterns and nonphysical oscillations in the pressure field. We examine the performance of our scheme over the time to address this issue from a dynamic point of view. 
		\end{itemize}
\end{enumerate}
\subsubsection{Incompressible elasticity problem} \label{sec:inf-sup}
\mbox{}\\The full discretized matrix form of the incompressible elasticity problem can be extracted from \eqref{Matrix form} by setting $\bld q=0$ as
\begin{equation}    \label{incompressElasticity} 
   \begin{bmatrix}   \bar{M}  & -\bar{Q} \\ 
			-\bar{Q}^\top & 0 \end{bmatrix}
			 \begin{bmatrix} \bld {\dot u}_{\text{n+1}} \\  \bld{p}_{\text{n+1}} \end{bmatrix}
			= \begin{bmatrix} \bar{\mathbb{P}}\\  0 \end{bmatrix}
\end{equation}

\begin{figure} 
 	\centering
	\resizebox{0.5\textwidth}{!}{\includegraphics{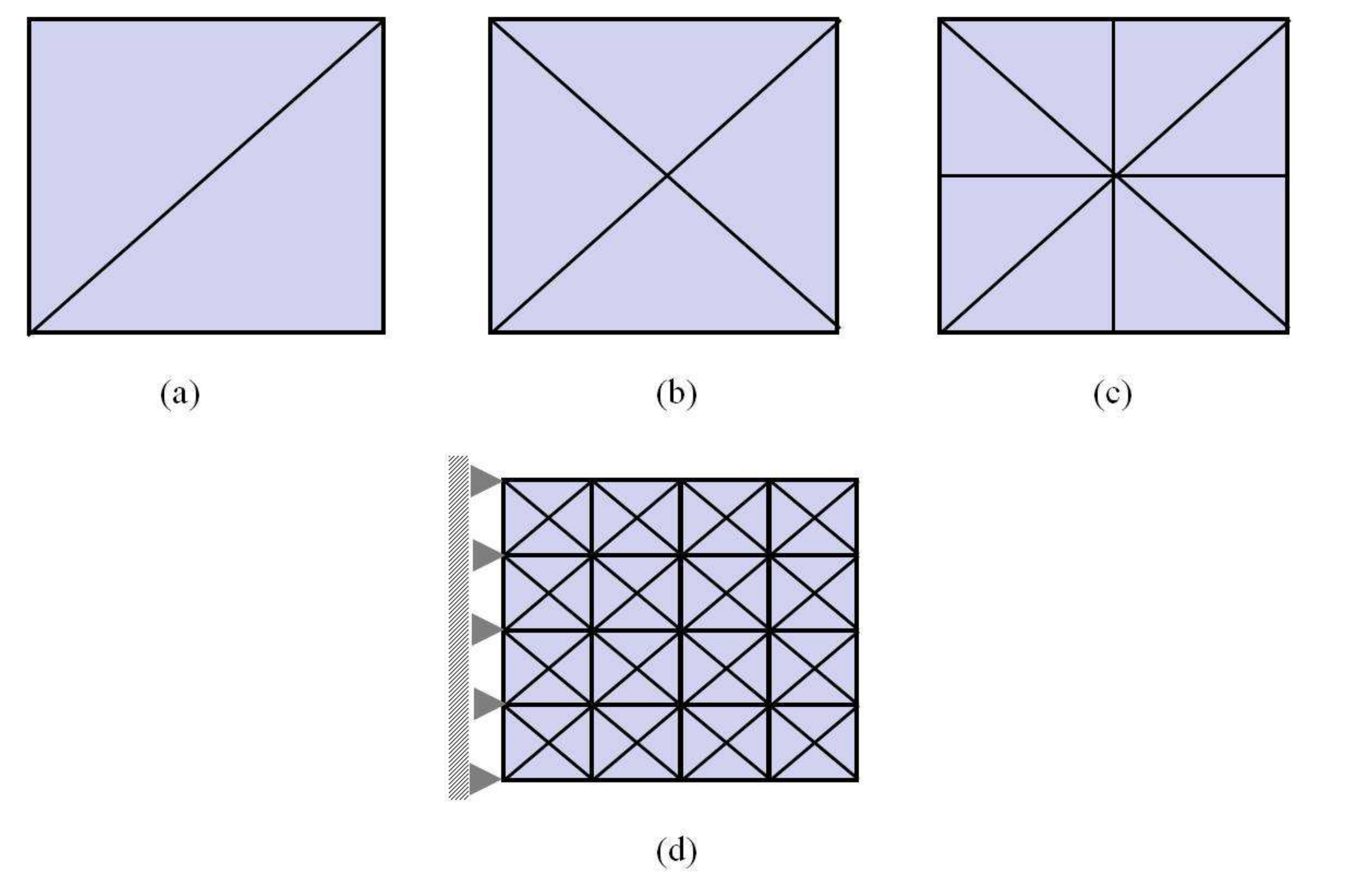}}  
	\caption[Mesh patterns and the problem considered for the inf-sup test] {Mesh patterns and the problem considered for the inf-sup test. (a) Criss pattern. (b) Crisscross pattern. (c) Union jack pattern. (d) Cantilever square bracket problem used for the inf-sup test, an instance of $N=4$ with a crisscross pattern is shown here.}
    \label{Meshpatterns_problem}
\end{figure}

\begin{figure}
 	\centering
	\resizebox{0.5\textwidth}{!}{\includegraphics{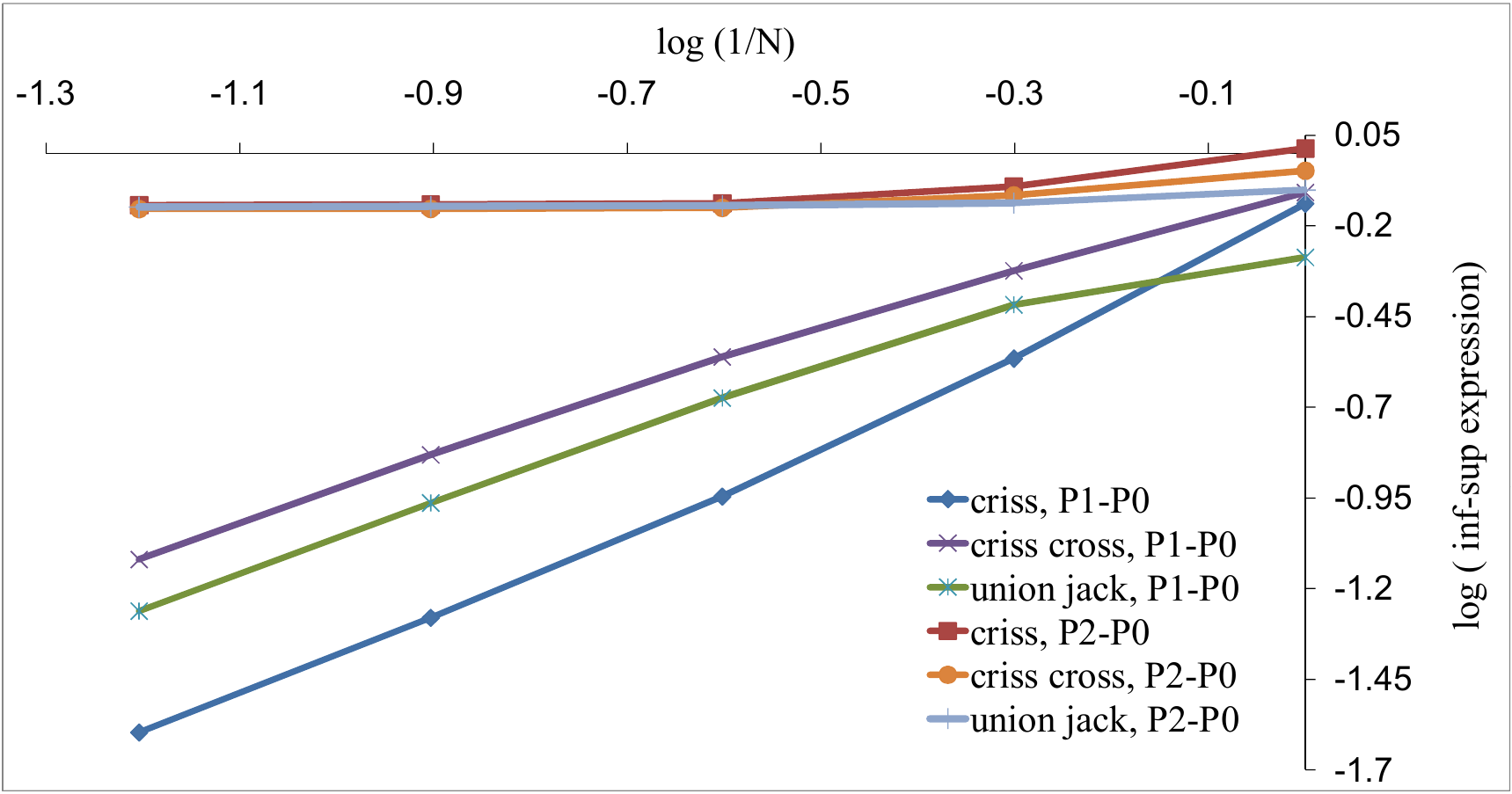}}	
	\caption[Numerical inf-sup test of triangular elements]{Numerical inf-sup test on finite elements. All $\quadElem$ elements satisfies the inf-sup condition, while the $\linElem$ elements do not.}
	\label{infSUp}   
\end{figure}

\noindent \textit{Spurious modes and inf-sup test}
\\We investigate on the stability of our finite elements from three perspective:
\begin{itemize}
\item Linear independence of constraint equation (with boundary conditions applied): one indicator of instability is the presence of spurious pressure modes, which occurs when the kernel of the $\bar Q$ is nontrivial. In this case, the spurious non-zero pressure mode satisfies 
$$\bar{Q} P_\text{s}=0$$
This is the result of constraint equations being linearly dependent. Clearly, in the presence of these spurious modes, the inf-sup expression yields exactly zero and the problem is not solvable. The existence of such non-zero $P_\text{s}$ can be crucially dependent on boundary conditions.
\item Linear independence of constraint equation (without boundary conditions applied): there is another type of spurious mode, which is the side effect of using certain mesh patterns.  These are element-wise (local) spurious modes due to redundant constraints at the macroelement level regardless of any imposed boundary conditions. An extensive discussion on this issue is found in \cite{hughes2012finite}.
\\We assess combination of several elements and mesh patterns in a two-dimensional plane strain problem and a cantilevered square block geometry (Figure \ref{Meshpatterns_problem}). Table \ref{tab:linearIndependency} shows the results of linear independence of the constraint block in equation \eqref{incompressElasticity} under various scenarios. Columns A and B report if the matrix of constraints, $\bar Q$, has full column-rank. It is observed that the $\linElem$ crisscross macroelement has one redundant constraint in each macroelement, which causes an element-wise spurious pressure; this was first discovered by Mercier \cite{mercier1979conforming}. We also note that the criss element does not contain any spurious pressure mode either in element level or in the assemblage level (when boundary conditions is applied). 
\item Uniform stability -- inf-sup expression: in order to fully verify the stability condition, the inf-sup expression must be bounded from below by a constant value independent of mesh size. The inf-sup condition is not simply equivalent to linear independence of constraint equations which is a purely algebraic property, it is rather an analytic one that ensure the \textit{uniform} linear independency as the mesh is increasingly refined \cite{braess2001finite}.
\\To perform the inf-sup test of Chapelle and Bathe \cite{chapelle1993inf} on our elements, we construct the matrix form of the inf-sup condition and normalize it. The corresponding eigenvalue problem is then solved and the smallest eigenvalue is considered as the inf-sup expression. We observe that the normalization by the displacement norm is inevitable otherwise the smallest singular values of the constraint matrix alone may never approach to any limit (as occurred for the $\quadElem$ element without normalization). It is noted that in the presence of spurious modes the inf-sup expression is exactly zero, however, we use the smallest non-zero eigenvalue to build the graph of Figure \ref{infSUp}. We display in this figure the inf-sup expression evaluated for a series of five mesh refinements, $N=1,2,4,8,16$ where $N$ is the number of macroelements along each side (an example of $N=4$ is shown in Figure \ref{Meshpatterns_problem}d). 
We observe that the $\linElem$ finite element violates the inf-sup condition no matter what mesh pattern is used as indicated by inf-sup values unbounded when $N$ increases. In contrast, the $\quadElem$ element satisfies the inf-sup condition, indicated by asymptotically approaching of the inf-sup expression to a value greater than zero. Column C of Table \ref{tab:linearIndependency} summarizes the results of these tests. 
\end{itemize}
\begin{table}
            \centering
            \begin{tabular}{clccc}
            	\toprule   
		$\text{element type}^{\ast}$ & $\text{mesh pattern}^{\dagger}$ & A: $\text{element-wise}$ & B: $\text{linear}$&C: $\text{inf-sup}$ \\
             &&  linear independent$^{\dagger}$& $\text{independent}^{\ddagger}$& $\text{condition} ^{\ddagger}$\\
             &&   (no BCs) & (N=1) & (N $\rightarrow \infty$)\\
\midrule
             & criss & $\checkmark$ & $\checkmark$ & X \\
		$\linprt$ & crisscross & X & X & X   \\
		& union jack & $\checkmark$ & X & X \\
		\midrule
		 & criss & $\checkmark$ & $\checkmark$ & $\checkmark$ \\
		$\quadElem$ & crisscross & $\checkmark$ & $\checkmark$ & $\checkmark$   \\
		& union jack & $\checkmark$ & $\checkmark$ & $\checkmark$ \\      
            	\bottomrule
            \end{tabular}
            \caption[Linear independency of the constraint in the element and global level as well as the numerical inf-sup test] {Linear independency of the constraint in the element and global level (assemblage) as well as the numerical inf-sup test. The $\linElem$ crisscross macroelement contains element-wise and global spurious modes. $\linElem$ union jack macroelement has global spurious mode. \\
			    $\text{}^\ast$ See Figure \ref{fig:mixedElement}, $\text{}^\dagger$ See Figure \ref{Meshpatterns_problem}{a-c}, $\text{}^\ddagger$ See Figure \ref{Meshpatterns_problem}d }
            \label{tab:linearIndependency}
\end{table}
\subsubsection{Poroelasticity under very low permeability} \label{PoroWLowPermeabil}
\begin{figure} 
 	\centering
	\resizebox{0.4\textwidth}{!}{\includegraphics{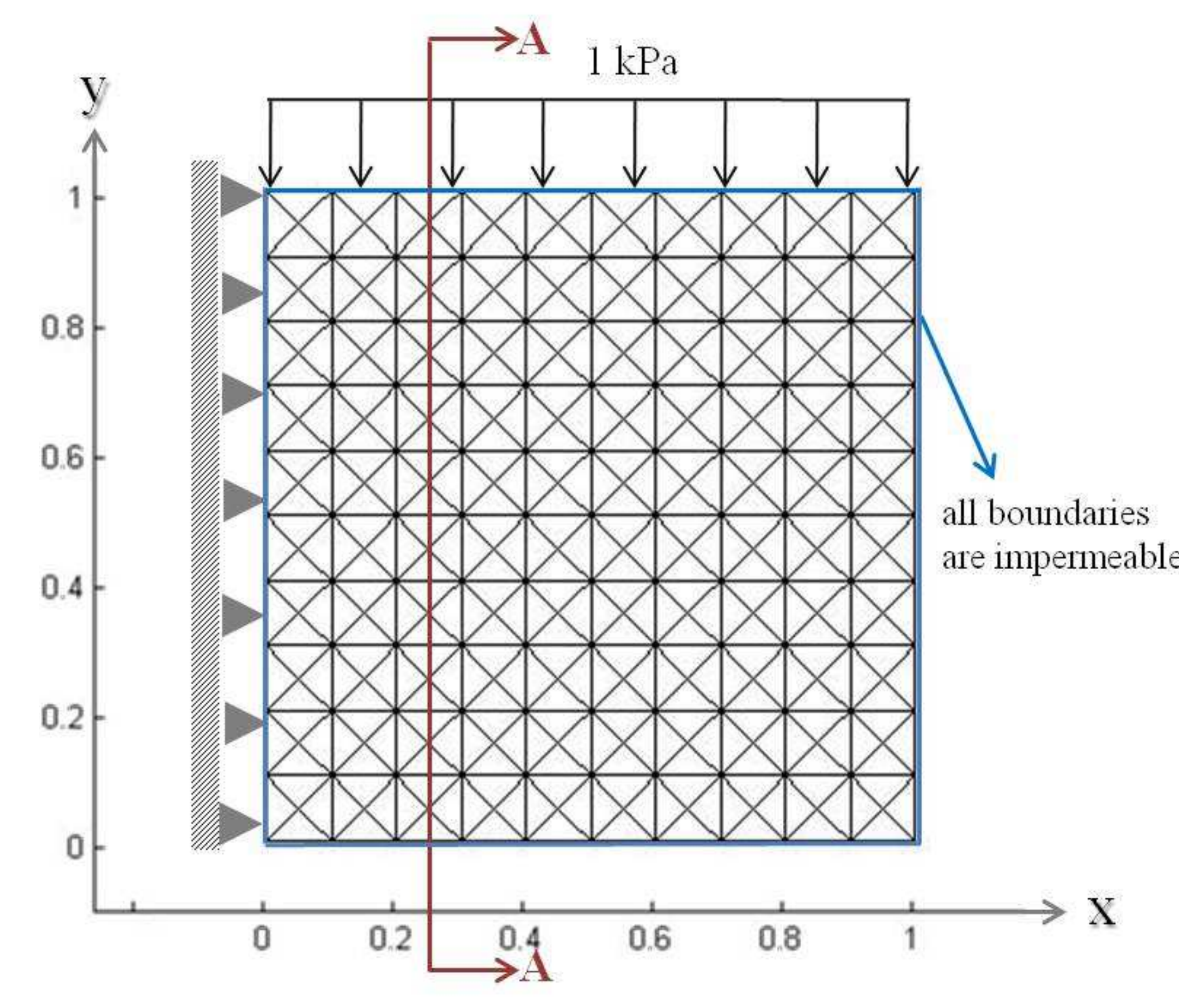}}
	\caption[Example 3 --- cantilever square bracket]{Example 3 --- cantilever square bracket: Geometry and finite element mesh}
\label{locking_geom}
\end{figure}

\mbox{}\\By passing the inf-sup test, the $\quadElem$ element is shown to be capable of overcoming the problem of locking in the incompressible elasticity problem. We wish to establish that consequently, it performs well in the low-permeability poroelasticity problems as well. We note that, in poroelasticity the $\linprt$ ($\quadprt$) element is the three field extension of $\linElem$ ($\quadElem$) elements introduced earlier for incompressible elasticity. 
In this section we study the configuration of Figure \ref{locking_geom} which is a 1 (m) $\times$ 1 (m) cantilever square bracket similar to the one
given in the previous section for inf-sup analysis. This problem was first used by Liu \cite{Liu_2004} and later by Phillips and Wheeler \cite{Phillips_Wheeler_2009} to demonstrate the problem of locking in poroelasticity. Plane strain behavior is assumed and the material properties, similar to those used by Liu \cite{Liu_2004}, are shown in Table \ref{tab:Material_properties}. To simulate the incompressible elasticity situation and induce locking, we study the case of very low permeability, $\kh=1 \times 10^{-1}$. A step load of 1 $\frac{\text{kN}}{\text{m}^2}$ is applied uniformly on top and all the boundaries are impermeable. Also on the fixed left boundary, the normal and tangential displacements to the surface are constrained.
\begin{figure}
    \centering
    \subfloat[Pressure along cross section A-A plotted at different \newline instances of time using the $\linprt$ element]
    {
        \resizebox{0.5\textwidth}{!}{\includegraphics{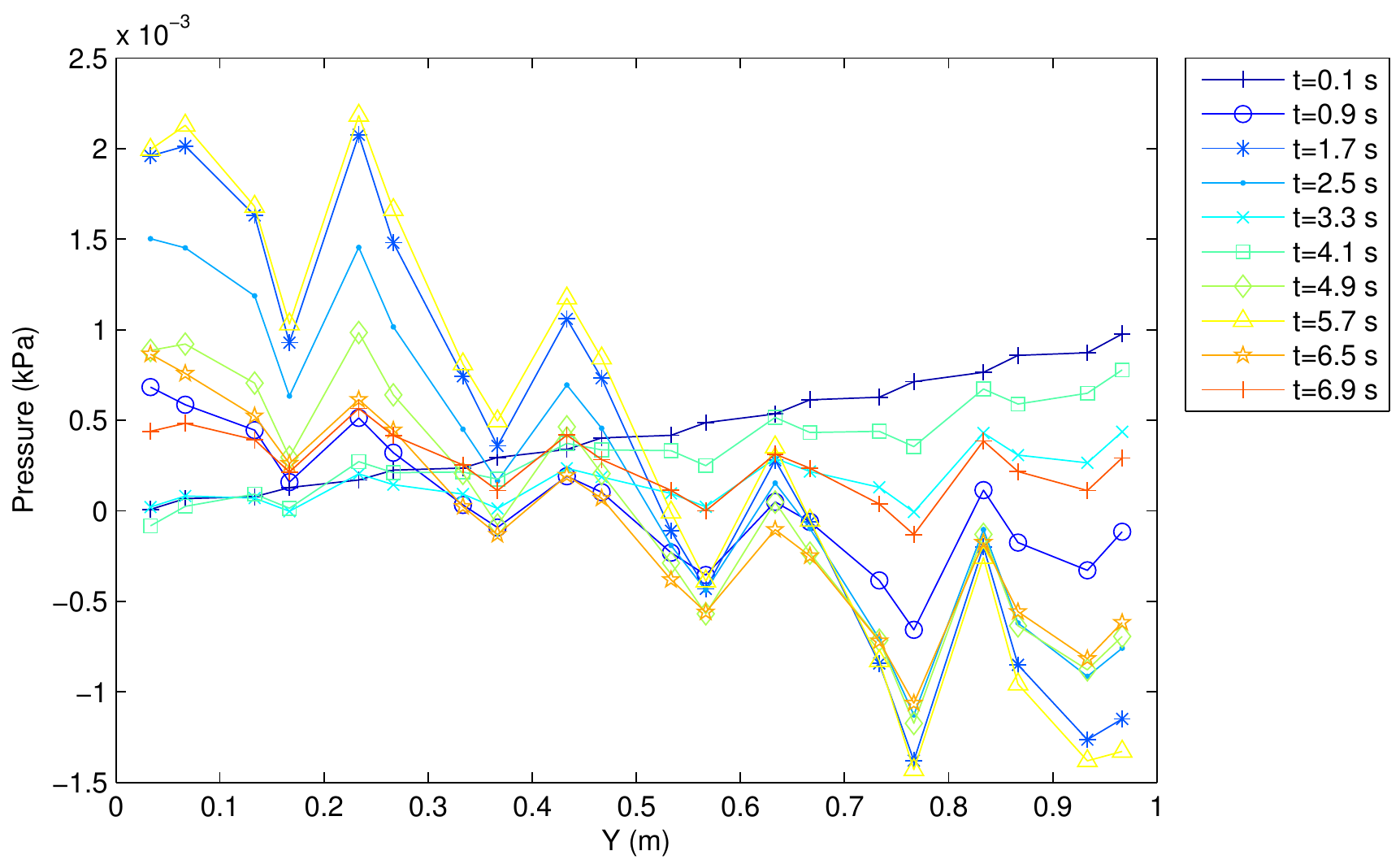}}
        \label{fig:CRS_time_Cst}
    }
    \subfloat[Pressure along cross section A-A plotted at different \newline instances of time using the $\quadprt$ element]
    {
        \resizebox{0.5\textwidth}{!}{\includegraphics{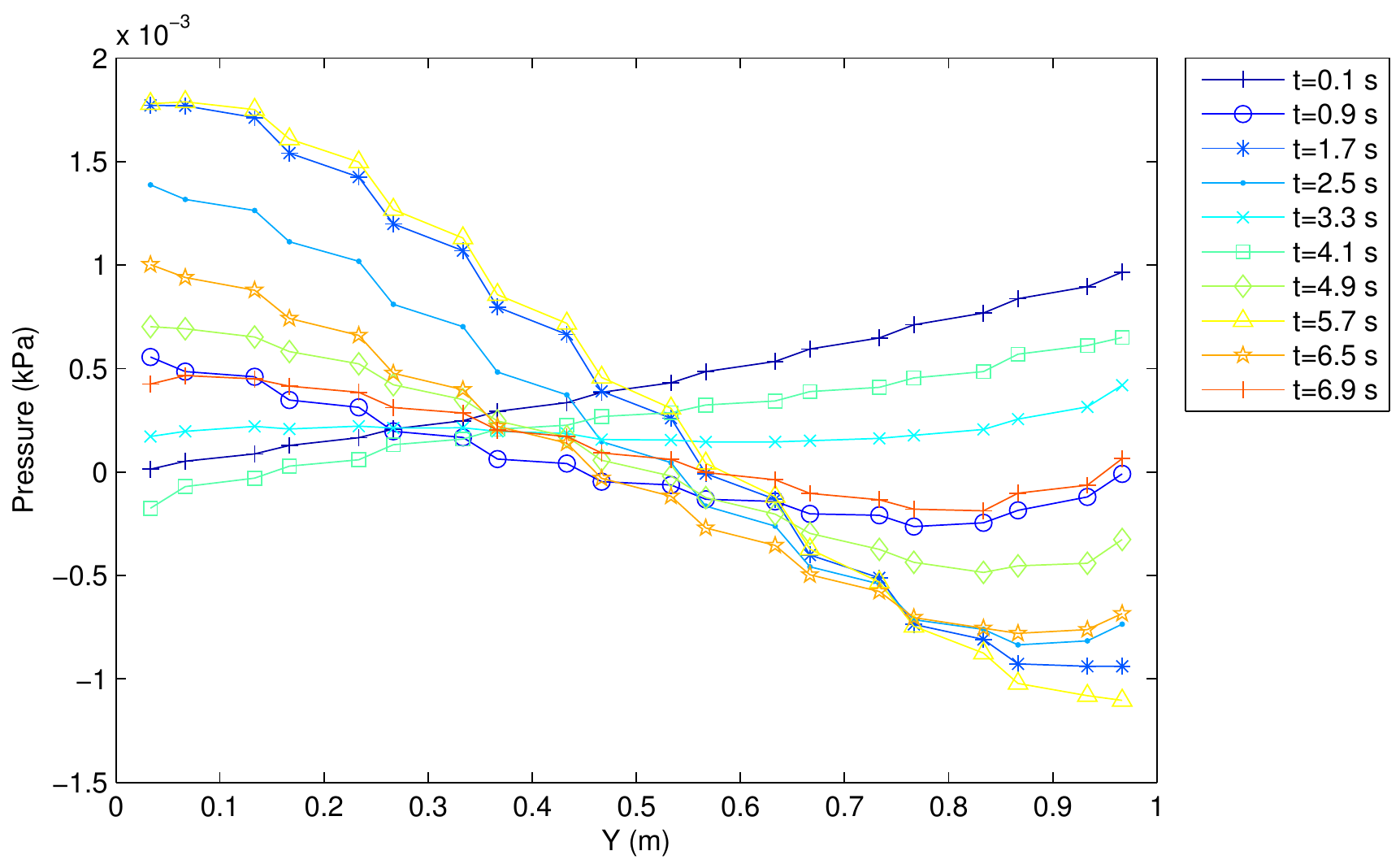}}
        \label{fig:CRS_time_Lst}
    }
\\
    \subfloat[Pressure along various cross sections \newline at time=2.5 s using the $\linprt$ element]
    {
        \resizebox{0.5\textwidth}{!}{\includegraphics{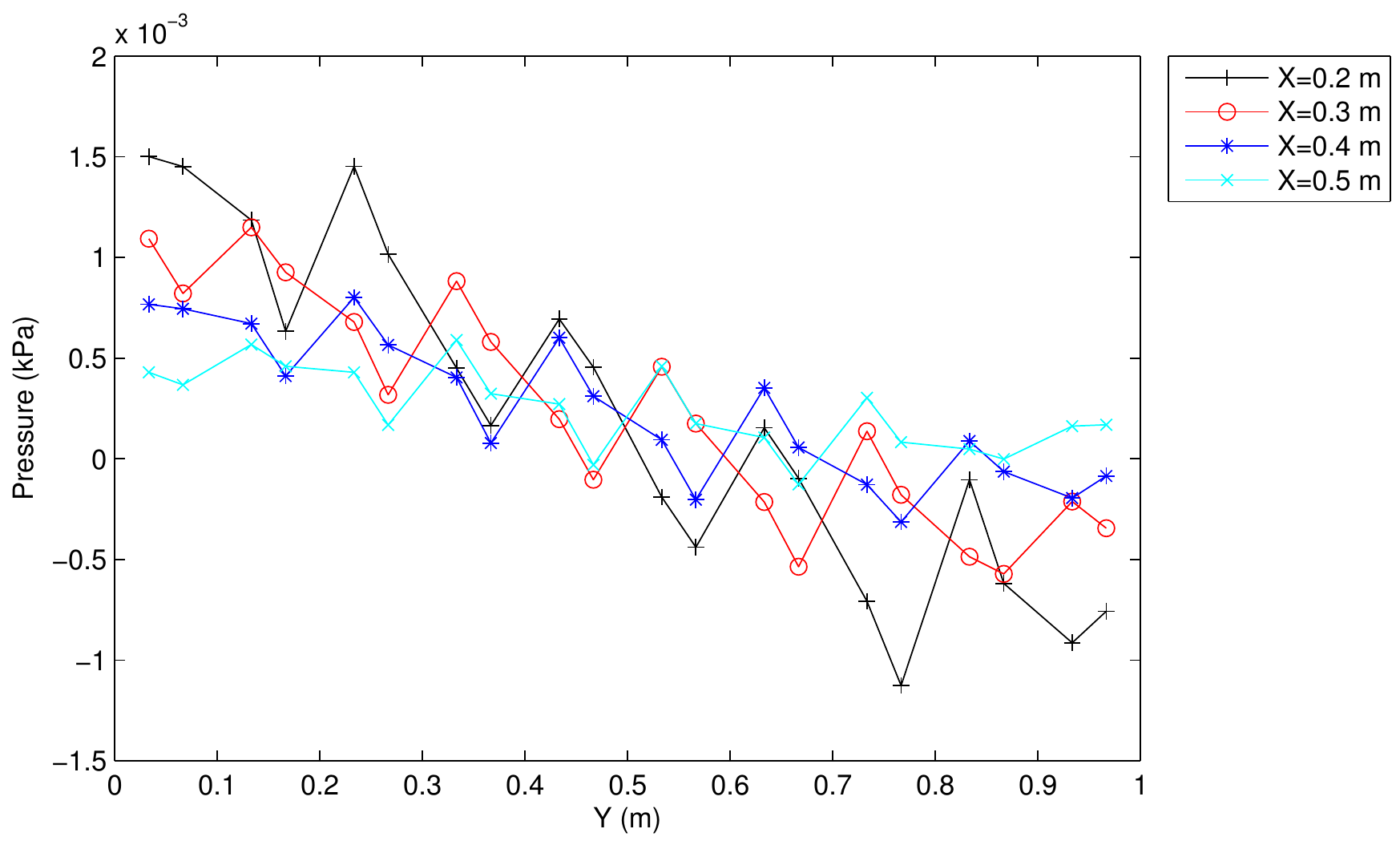}}
        \label{fig:CRS_sections_Cst}
    }    
    \subfloat[Pressure along various cross sections \newline at time=2.5 s using the $\quadprt$ element]
    {
        \resizebox{0.5\textwidth}{!}{\includegraphics{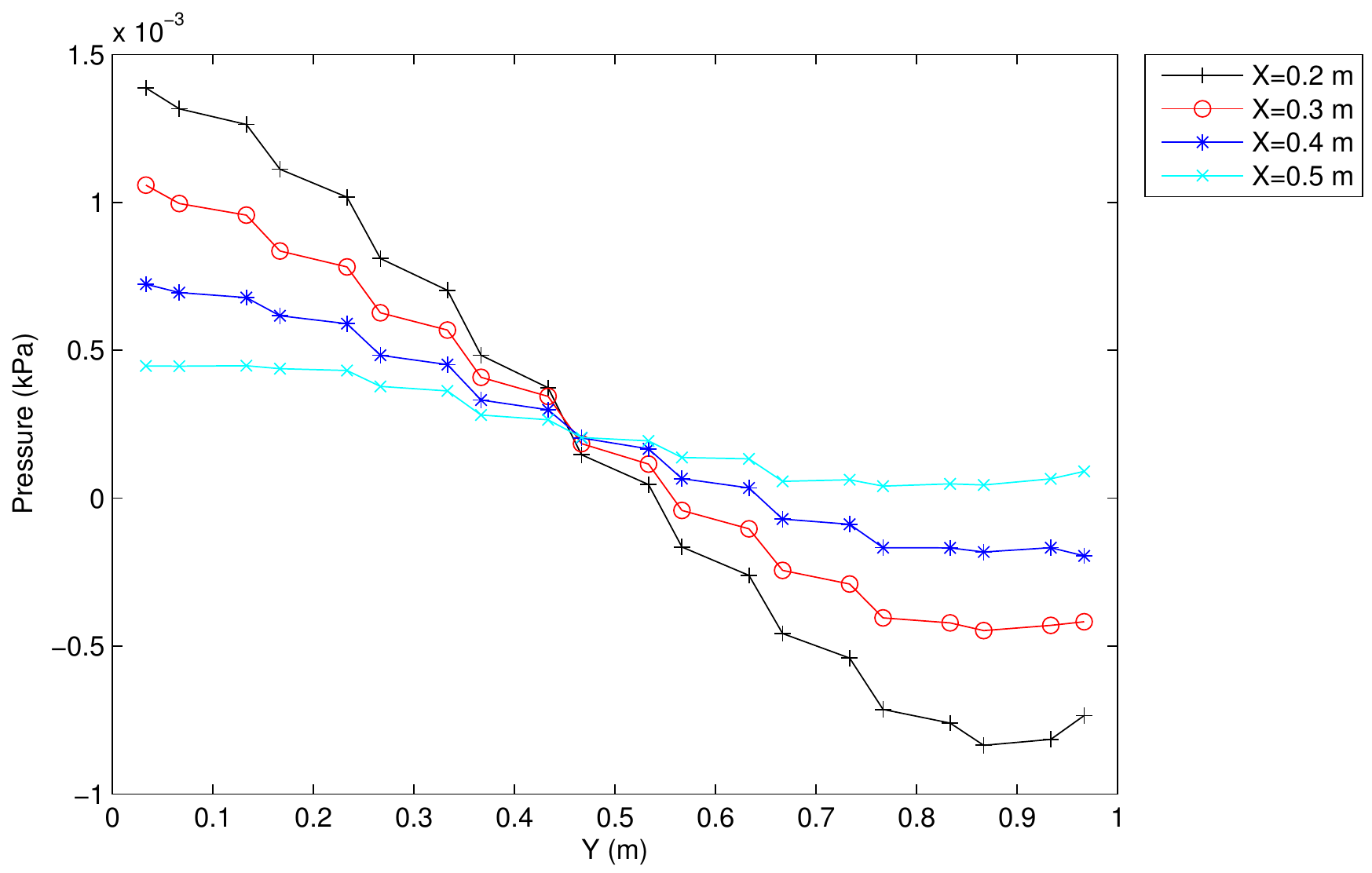}}
        \label{fig:CRS_sections_Lst}
    }
    \caption[Example 3--- nonphysical pressure oscillation] {Example 3 --- nonphysical pressure oscillation is visible at different moments of time and cross sections when using the $\linprt$ element, whereas a smooth pressure field is obtained when using the $\quadprt$ element}
    \label{fig:example3_CRS}
\end{figure}

\begin{figure}
    \centering
    \subfloat[time=0.1 s]
    {
        \resizebox{0.4\textwidth}{!}{\includegraphics{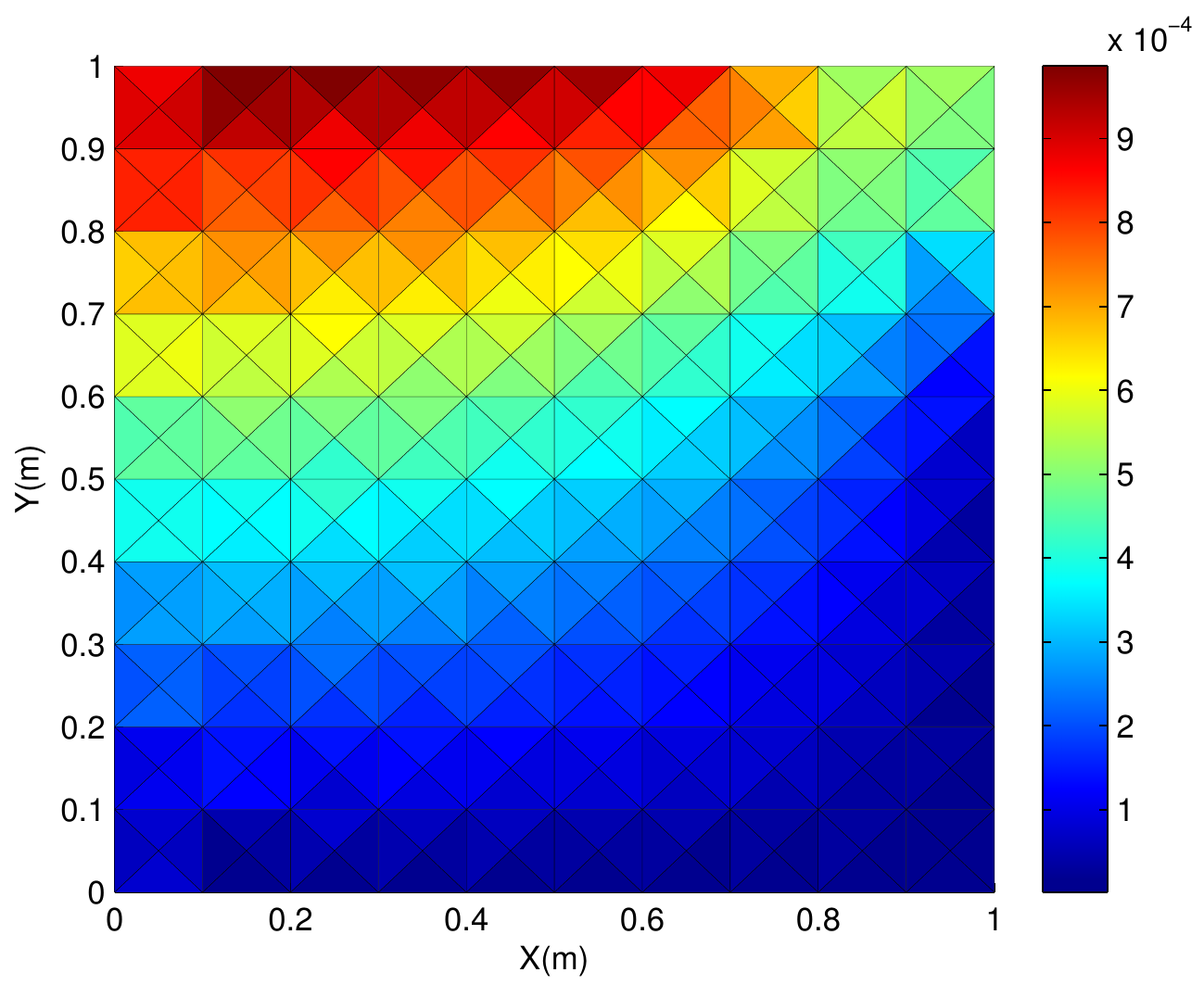}}
        \label{fig:patch1_CST}
    }
    \subfloat[time=0.1 s]
    {
        \resizebox{0.4\textwidth}{!}{\includegraphics{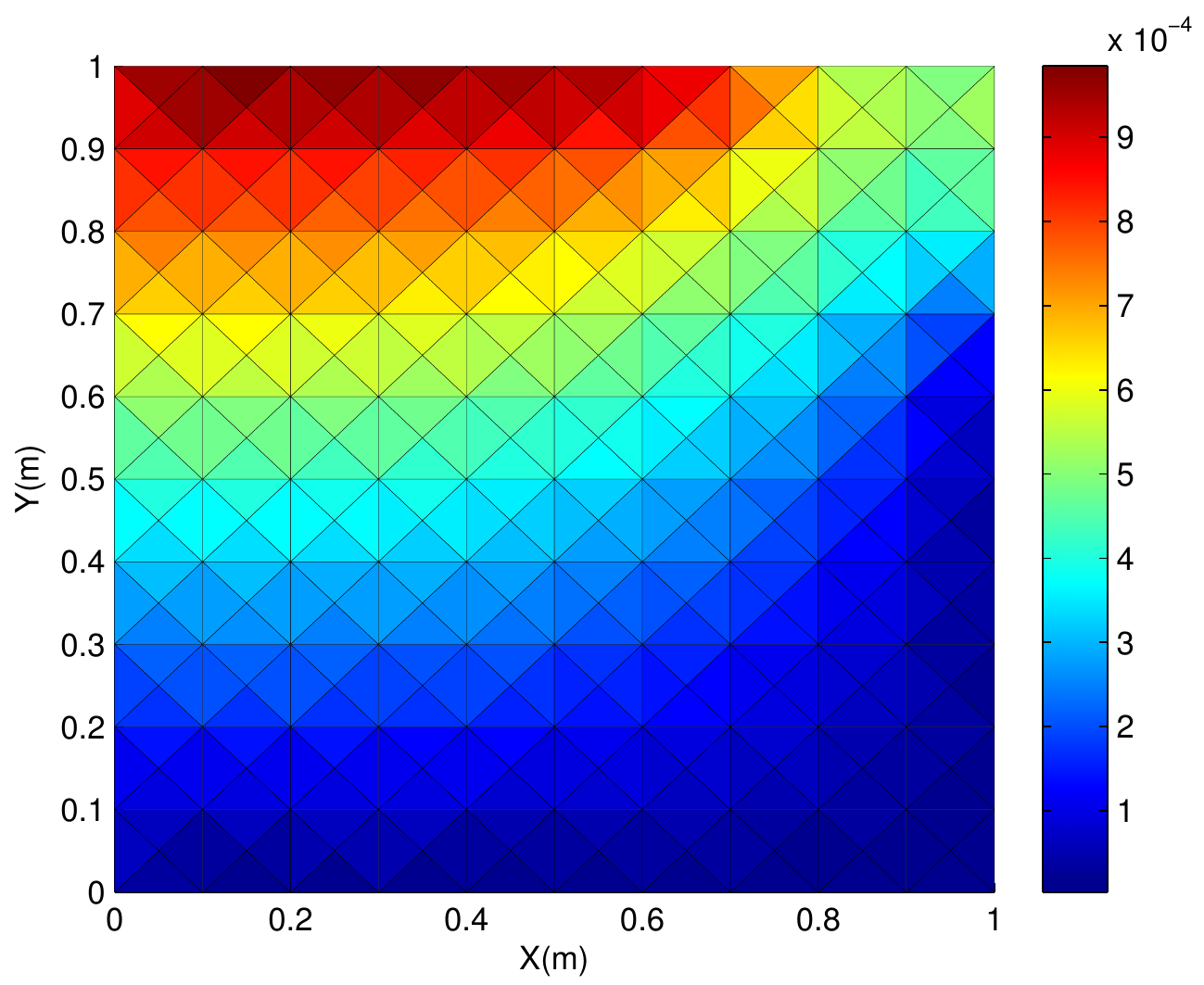}}
        \label{fig:patch1_LST}
    }
 \vspace{-1\baselineskip}
    \subfloat[time=0.3 s]
    {
        \resizebox{0.4\textwidth}{!}{\includegraphics{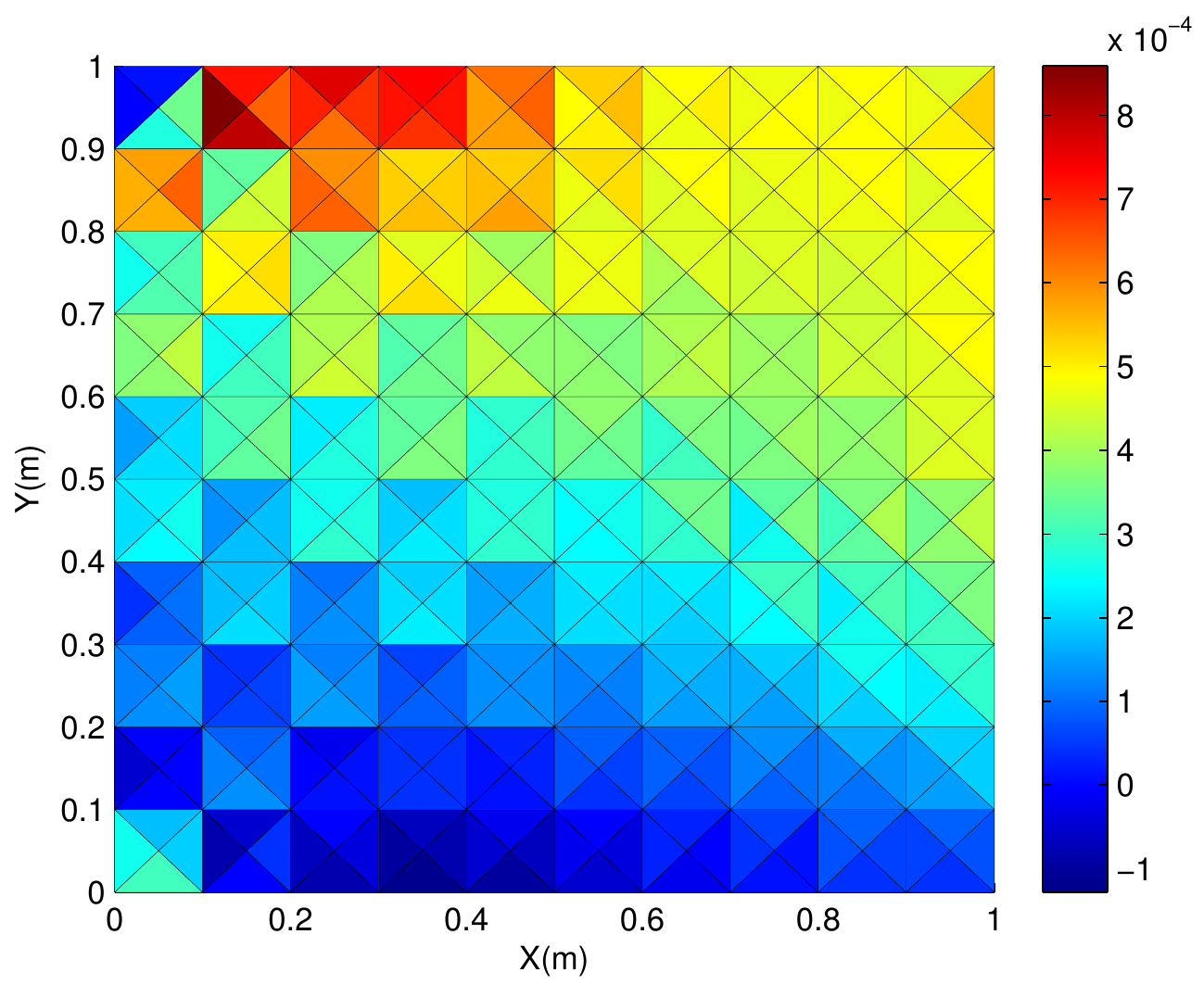}}
        \label{fig:patch2_CST}
    }   
    \subfloat[time=0.3 s]
    {
        \resizebox{0.4\textwidth}{!}{\includegraphics{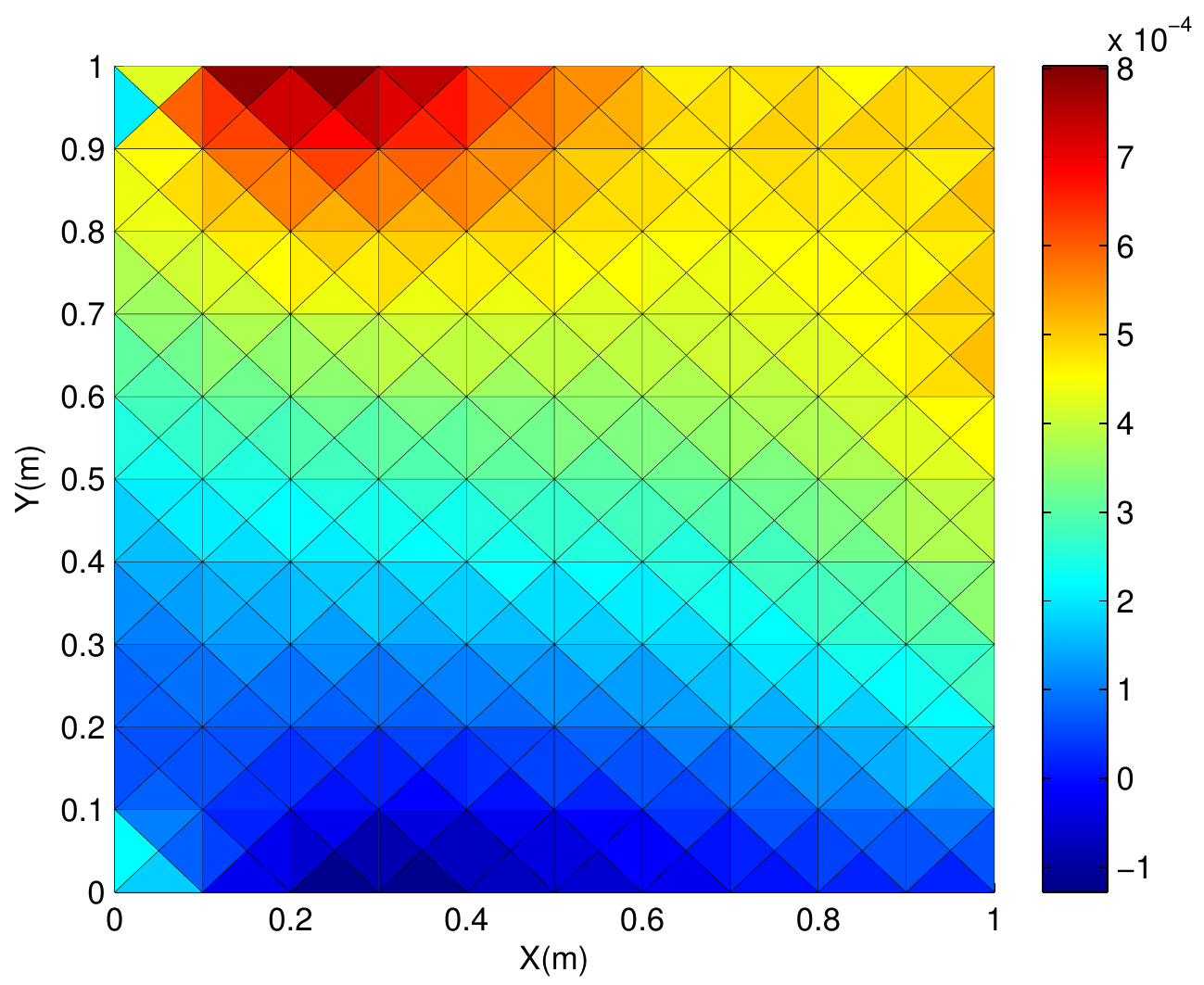}}
        \label{fig:patch2_LST}
    }
 \vspace{-1\baselineskip}
    \subfloat[time=4.1 s]
    {
        \resizebox{0.4\textwidth}{!}{\includegraphics{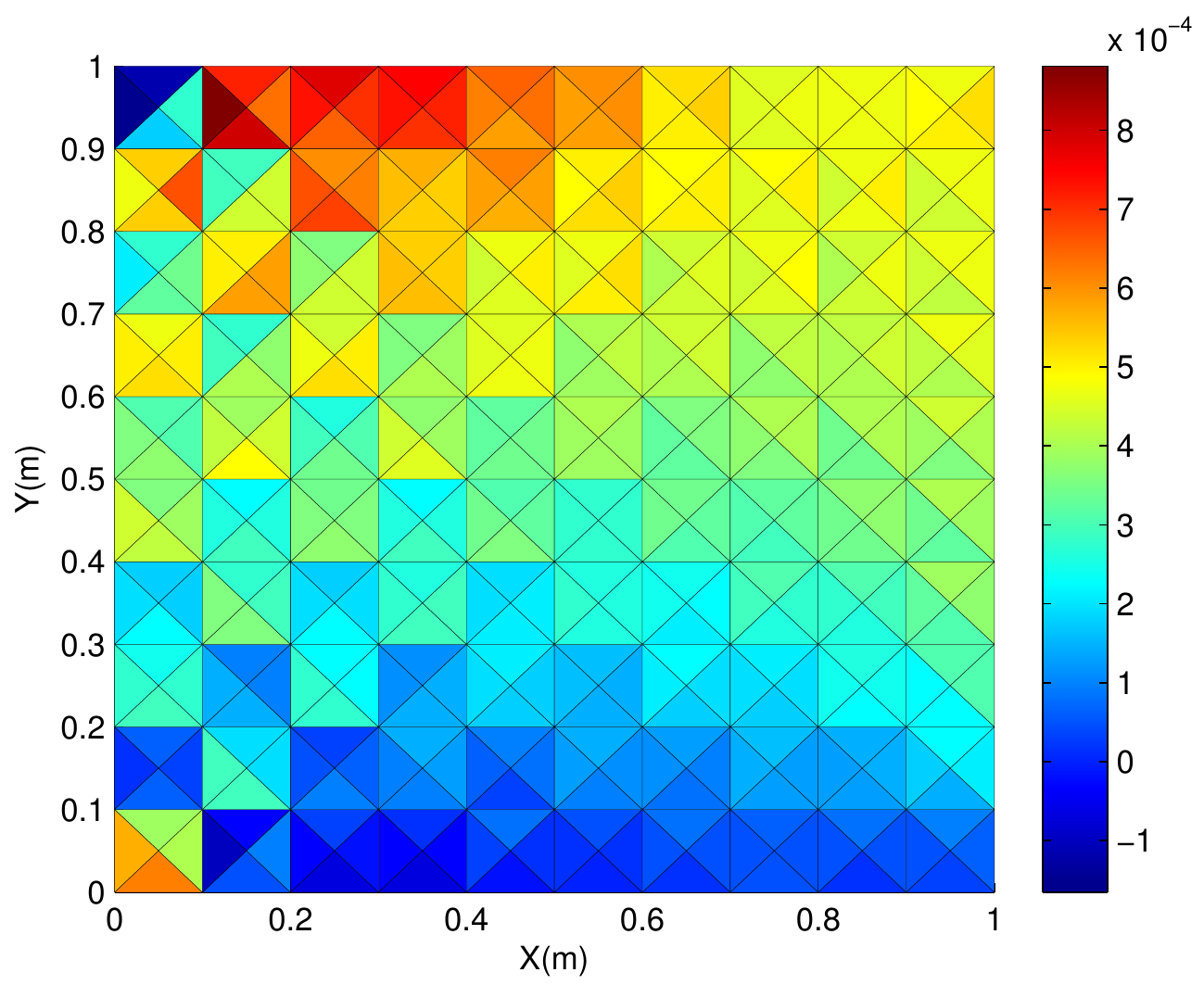}}
        \label{fig:patch4_CST}
    }    
    \subfloat[time=4.1 s]
    {
        \resizebox{0.4\textwidth}{!}{\includegraphics{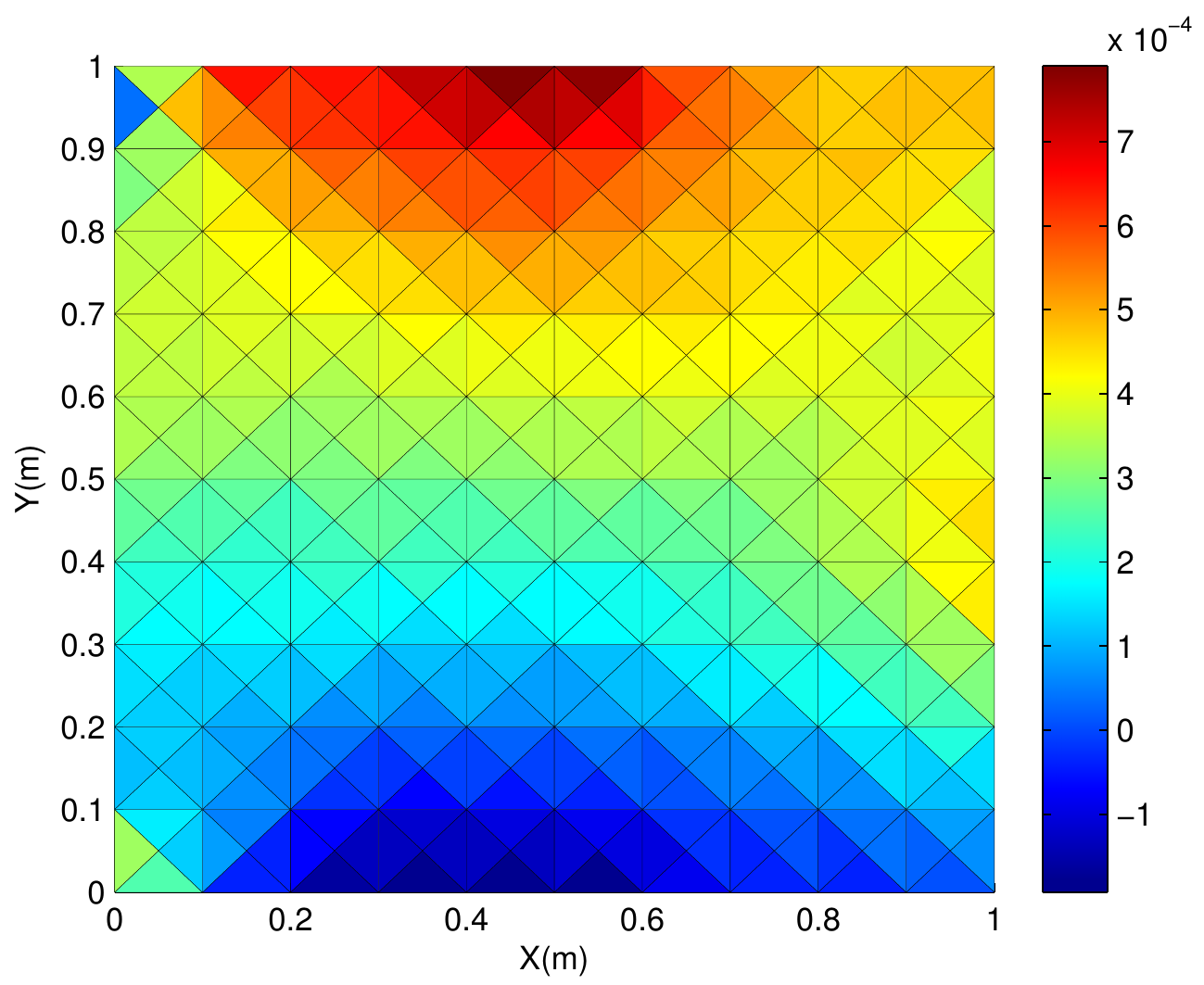}}
        \label{fig:patch4_LST}
    }
 \vspace{-1\baselineskip}
   \subfloat[time=6.3 s]
    {
        \resizebox{0.4\textwidth}{!}{\includegraphics{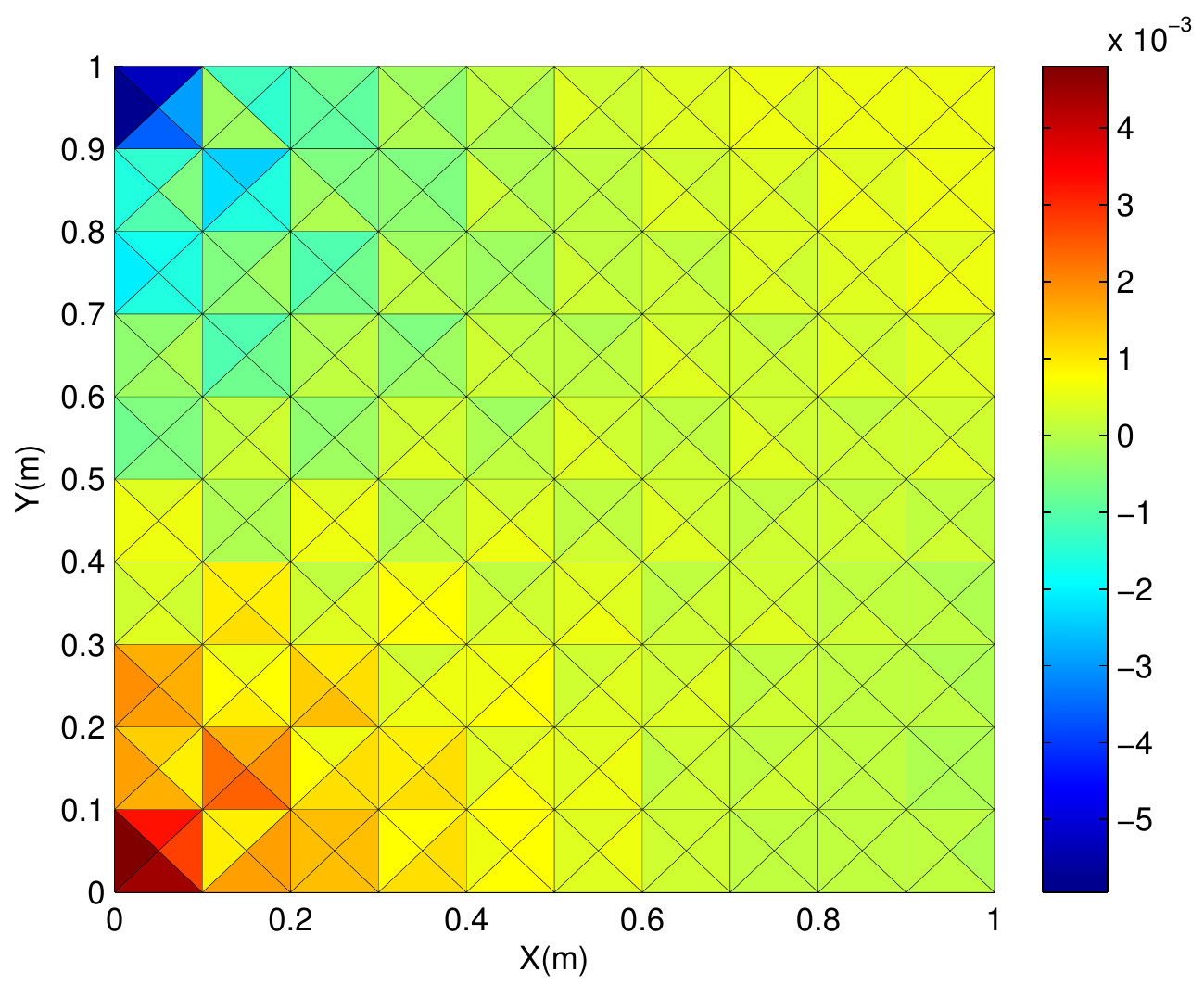}}
        \label{fig:patch5_CST}
    }    
    \subfloat[time=6.3 s]
    {
        \resizebox{0.4\textwidth}{!}{\includegraphics{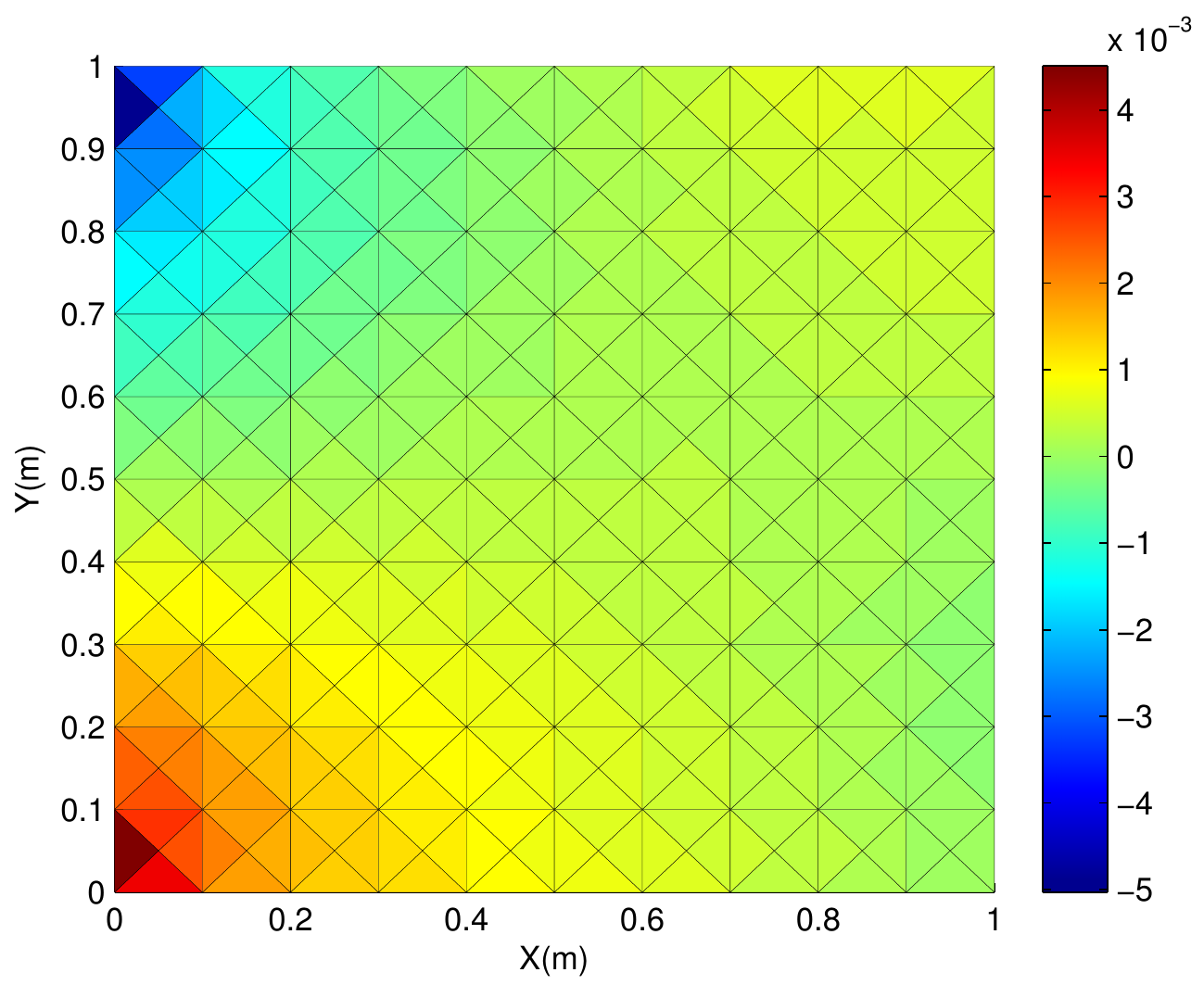}}
        \label{fig:patch5_LST}
    }
    \caption[Example 3 --- pressure contours at some representative moments of time] {Example 3 --- the left pressure contours illustrate the results of the $\linprt$ element and its associated locking response with the checkerboard pattern, while the right ones are for the $\quadprt$ stable and locking free elements. Note that the colorbar is not identical in all plots.}
   \label{fig:example3_patch}
\end{figure}

\noindent \textit{Checkerboard pattern in the pressure field}
\\In Figures \ref{fig:CRS_time_Cst} and \ref{fig:CRS_time_Lst}, we present pressure distributions along cross section A-A ($X=0.25$ m) for $\linprt$ and $\quadprt$ elements respectively. These are plotted at different instances of time. 
The nonphysical oscillations and checkerboard pattern in the pressure field in Figure \ref{fig:CRS_time_Cst} reveals the problem of the $\linprt$ element in handling low hydraulic conductivity as also reported in \cite{Liu_2004,Phillips_Wheeler_2009}. However, implementing a quadratic skeleton displacement field as in the $\quadprt$ element resolves the problem and results in a smooth pressure field (Figure \ref{fig:CRS_time_Lst}). The pressure field at other cross sections is plotted in Figures\ref{fig:CRS_sections_Cst} and \ref{fig:CRS_sections_Lst} at time ($t=2.5$), showing the same phenomena. 
To get further insight, the pressure over the whole domain is plotted in Figure \ref{fig:example3_patch} at different instances of time confirming similar observations. One indication of instability is the checkerboard patterns on the pressure field with the $\linprt$ element (left figures), whereas the $\quadprt$ element (right figures), which satisfies the inf-sup condition, effectively eliminates locking.

\begin{figure} 
 	\centering
	\resizebox{0.5\textwidth}{!}{\includegraphics{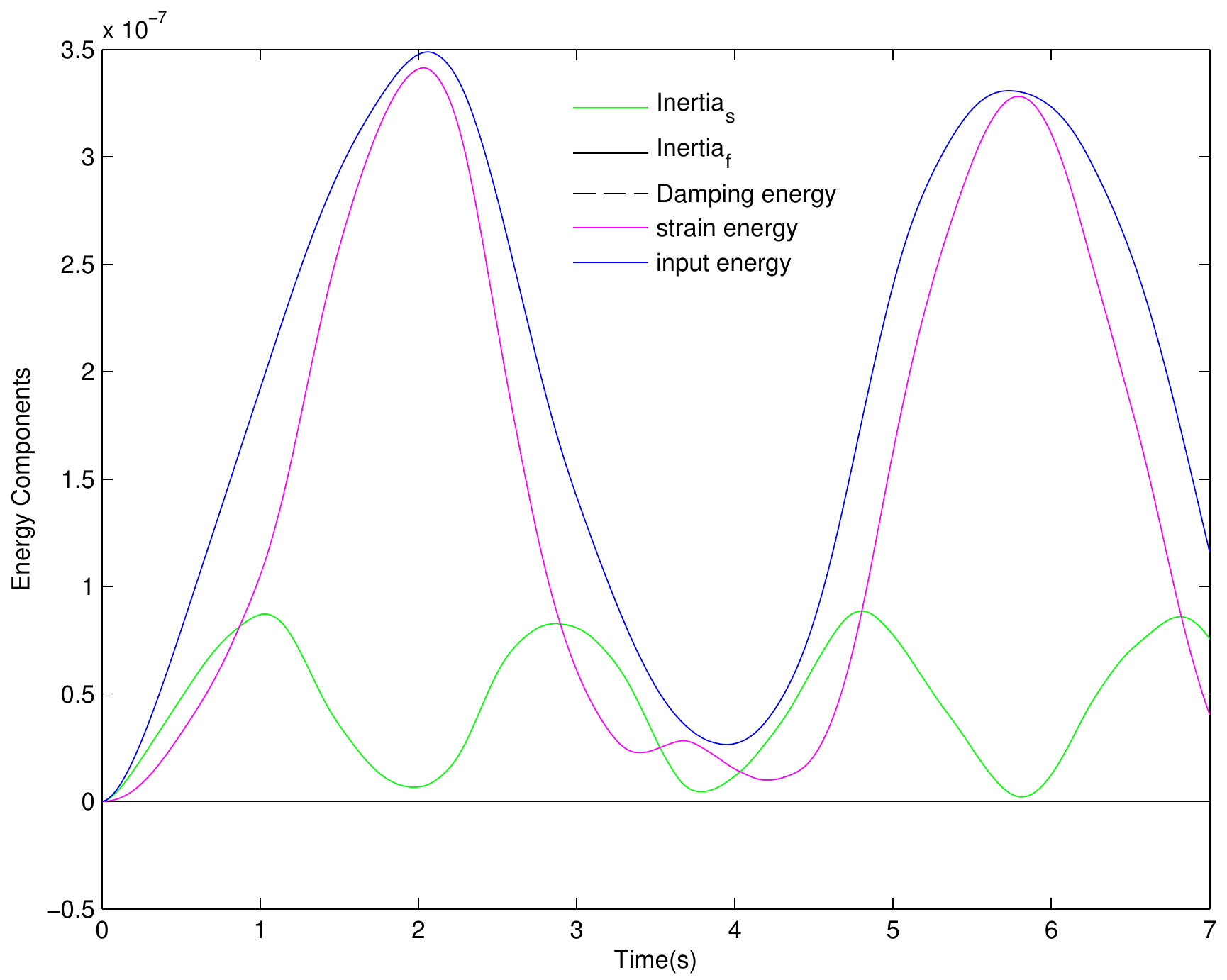}}
	\caption[Example 3 --- components of energy]{Example 3 --- components of energy. Initially, the kinetic energy of skeleton is dominant that prevents the locking even in unstable elements shown in Figure \ref{fig:patch1_CST} }
\label{locking_energy}
\end{figure}

\noindent \textit{Dynamic analysis affects oscillations}: 
\\We observe in Figure \ref{fig:CRS_time_Cst} the incapability of $\linprt$ element in computing a correct pressure field, with the checkerboard patter persisting over time. In contrast, in quasistatic poroelasticity, it has been observed that these instabilities are pronounced in the beginning, but tend to disappear as time progresses \cite{Phillips_Wheeler_2009,Murad_Loula_1994,Murad_etal_1996}. Such behavior in dynamic analysis with low hydraulic conductivity could be explained by the periodic deformation and velocity time histories; in such periodic condition, one encounter the static problem over and over, leading to a frequent instability issue.
One a different note, by comparing Figure \ref{fig:patch1_CST} with \ref{fig:patch1_LST}, we see that in early moments of time (up to about $t=0.1$ s) both $\linprt$ and $\quadprt$ elements predict similar and smooth pressure fields. One possible reason could be the sudden application of step load, which initially creates a dominant skeleton kinetic energy as illustrated in Figure \ref{locking_energy}. 
\subsubsection{Summary of findings}
\begin{itemize}
\item The stability of our finite elements in terms of local and global spurious modes as well as the inf-sup condition is studied (Table \ref{tab:linearIndependency}). It is found that the crisscross mesh pattern on the $\linElem$ element contains a redundant constraint, which causes an element-wise spurious pressure. Also, results indicate that the $\linElem$ element violates the inf-sup condition, whereas the $\quadElem$ element satisfies it.
\item The checkerboard pattern and nonphysical oscillation in the pressure field, shown in Figures \ref{fig:CRS_time_Cst} and left Figures of \ref{fig:example3_patch}, reveals the problem of the $\linprt$ element in handling dynamic analysis of porous medium with low hydraulic conductivity. Therefore, we introduce a quadratic skeleton displacement field in the $\quadprt$ element. This resolves the problem and results in a smooth pressure field (Figure \ref{fig:CRS_time_Lst} and right Figures in \ref{fig:example3_patch}). 
\item  In dynamic analysis, the checkerboard patterns produced by an unstable finite element is shown to persist over time and is not eliminated, in contrast to static analysis, where these instabilities appear early but disappear at later time.
\end{itemize}
\section{Concluding remarks}
An accurate and stable numerical framework for dynamics of incompressible saturated porous media is presented and verified over a variety of benchmark numerical studies. Our topology-motivated finite element is based on coupling of the mixed $\rt$ and the nodal Galerkin elements, implemented based on a three-field $(u-w-p)$ formulation. As a result, we do not neglect fluid acceleration terms, which enables applying the method to problems with considerable fluid dynamics such as soil liquefaction \cite{Jeremić_etal_2008} and biomechanics of porous tissues under rapid external loading \cite{Yang_2006,Levenston_et al_1998}. We also demonstrate stability of the method regarding the LBB condition in the limiting cases (rigid skeleton and very low hydraulic conductivity). In particular, the $\quadprt$ element, which uses a quadratic Galerkin discretization for the skeleton displacement, overcomes locking, which otherwise manifests in the form of checkerboard patterns in the pressure field. We also observe that with an unstable finite element, checkerboard patterns can persist over time. Given that we directly calculate the physical damping in the form of fluid diffusion, our numerical results suggest a decrease in damping with decreasing hydraulic conductivity, while frequency content is not affected by hydraulic conductivity. It is noteworthy that an appropriate prediction of dynamic behavior of porous with small hydraulic conductivity requires specific attention to the choice of finite element, mesh pattern, and mesh size.
\section{Acknowledgment}
The authors wish to thank Professors Gary F. Dargush for providing his comments to various aspects of this paper, in particular on the boundary
element method results. 

\printnomenclature[0.18\textwidth]

\newpage
\bibliography{Dynamic_poroElasticity}
\end{document}